\documentclass[a4paper,10pt]{book}
\usepackage[cp1251]{inputenc}
\usepackage{amsmath, amsfonts, amssymb}
\usepackage[russian]{babel}
\usepackage{makeidx}
\usepackage{calrsfs}
\usepackage{graphicx}
\usepackage{epsfig}
\usepackage {curves, eepic, epic}
\usepackage{wrapfig}
 \makeindex

 \DeclareMathAlphabet{\mathbbt}{U}{bbold}{m}{n}

\begin{document}

%%%%%%%%%%%%%%%%%%%%%%%%%%%%%%%%%%%%%%%%%%%%%%%%%%%%%%%%
\makeatletter%{@}
\newcommand{\Author}[2]{%\renewcommand{\thefootnote}{\fnsymbol{}}
\begin{center}\textbf{#1}%\footnote{$\copyright$\ 2004\ #2}
\end{center}\medskip\renewcommand{\@evenhead}{\raisebox{-3mm}[\headheight][0pt]%
{\vbox{\hbox to\textwidth{\thepage \hfill\strut {\small\sl
#2}}\hrule}}} }

\newcommand{\shorttitle}[1]{\renewcommand{\@oddhead}{\raisebox{-3mm}[\headheight][0pt]%
{\vbox{\hbox to\textwidth{\strut {\small\sl #1}\rightmark
\hfill\thepage}\hrule}}} }

\makeatother%@
%\newpage\pagestyle{plain}

 \newcommand{\UDC}[1]{\begingroup\newpage\thispagestyle{empty}\begin{flushleft}УДК #1\end{flushleft}}
 \newcommand{\Title}[1]{\begin{center}\uppercase{#1}\end{center}}
 \def\Adress#1#2{\bigskip{\footnotesize#1\par E-mail:\ #2\par} \newpage\endgroup}

  \renewcommand{\section}[1]{\par\medskip\begin{center}{\bf#1}\end{center}
  \addcontentsline{toc}{contsec}{\protect  #1}}

 \def\subsec#1{\smallskip\textbf{#1}}
 \def\subsubsec#1{\qquad{\bf #1}}
 \newcommand{\Lit}{\par\bigskip\centerline{\bf Литература}\smallskip
  \addcontentsline{toc}{contsec}{\protect  Литература}}
 \newcommand{\bib}[2]{{\leftskip-10pt\baselineskip=11pt\footnotesize\item{}\textsl{#1}~#2\par}}

\newcommand{\Def}[1]{\smallskip\par{\sc Определение#1.}}
\newcommand{\Zam}[1]{\smallskip\par{\sc Замечание#1.}}
\newcommand{\zam}[1]{{\sc Замечание#1.}}
\newcommand{\Primer}[1]{\smallskip\par{\sc Пример#1.}}
\newcommand{\primer}[1]{{\sc Пример#1.}}

\newcommand{\Proclaim}[1]{\smallskip{\bf#1}\sl}
\newcommand{\proclaim}[1]{{\bf#1}\sl}
\newcommand{\Theorem}[1]{\smallskip\par\textbf{Теорема#1.}\sl}
\newcommand{\theorem}[1]{\textbf{Теорема#1.\/}\sl}
\newcommand{\Sledstvie}[1]{\smallskip\par\textbf{Следствие#1.}\sl}
\newcommand{\sledstvie}[1]{\textbf{Следствие#1.}\sl}
\newcommand{\Endproc}{\rm}
\newcommand{\Lemma}[1]{\smallskip\par\textbf{Лемма#1.}\sl}
\newcommand{\lemma}[1]{\textbf{Лемма#1.}\sl}
\newcommand{\Predl}[1]{\smallskip\textbf{Предложение#1.}\sl}
\newcommand{\Utv}[1]{\smallskip\textbf{Утверждение#1.}\sl}
\newcommand{\utv}[1]{\textbf{Утверждение#1.}\sl}

\def\beginproof{\par\mbox{$\vartriangleleft$}}
\def\Beginproof{\smallskip\par\text{$\vartriangleleft$}}
\def\endproof{\text{$\vartriangleright$}}
\def\Endproof{\text{$\vartriangleright$}\smallskip}

%%%%%%%%%%%%%    неравенства
\renewcommand{\leq}{\leqslant}
\renewcommand{\geq}{\geqslant}
\renewcommand{\le}{\leqslant}
\renewcommand{\ge}{\geqslant}
%\newcommand{\loc}{\textrm{loc}\,}

%%%%%%%%%%%%%% o-lim  и др.
\def\osum{\mathop{o\text{\/-}\!\sum}}  %\def\osum{o\text{-}\!\sum}
\def\bosum{bo\text{-}\!\sum}
\def\olim{\mathop{o\text{-}{\fam0 lim}}} %\def\olim{\mathop{o\text{-}\!\lim}}
\def\bolim{\mathop{bo\text{-}\!\lim}}
\def\rlim{\mathop{r\text{-}{\fam0 lim}}}

%%%%%%%%%%%%%%%%% Операторы
\def\Ann{\mathop{\fam0 Ann}}
\def\Ban{\mathop{\fam0 Ban}}
\def\card{\mathop{\fam0 card}}
\def\Clop{\mathop{\fam0 Clop}}
\def\Dom{\mathop{\fam0 Dom}}
\def\dom{\mathop{\fam0 dom}}
\def\End{\mathop{\fam0 End}}
\def\Fin{\mathop{\fam0 fin}\nolimits}
\def\fin{\mathop{\fam0 fin}\nolimits}
\def\Id{\mathop{\fam0 Id}}
\def\ind{\mathop{\fam0 ind}}
\def\rot{\mathop{\fam0 rot}}
\def\loc{\mathop{\fam0 loc}}
\def\const{\mathop{\fam0 const}}
\def\span{\mathop{\fam0 span}}
\def\Exp{\mathop{\fam0 Exp}}
\def\grad{\mathop{\fam0 grad}}
\def\proj{\mathop{\fam0 proj}}
\def\im{\mathop{\fam0 im}}
\def\Lip{\mathop{\fam0 Lip}}
\def\mes{\mathop{\fam0 mes}}
\def\mix{\mathop{\fam0 mix}}

\def\On{\mathop{\fam0 On}}
\def\Orth{\mathop{\fam0 Orth}}
\def\sa{\mathop{\fam0 sa}}
\def\SC{\mathop{\fam0 SC}}
\def\sign{\mathop{\fam0 sign}}
\def\St{\mathop{\fam0 St}}
\def\supp{\mathop{\fam0 supp}}
\def\Re{\mathop{\fam0 Re}}
\def\ZFC{\mathop{\fam0 ZFC}}

%%%%%%%%% математические шрифты
\def\Cal#1{\mathcal{#1}}   %калиграфия (подключить пакет calrsfs)
\def\goth#1{\mathfrak{#1}} %готика
%\mathbbm --- полые буквы (\mathbb другой полый шрифт)
\def\bbold#1{{\mathbbt #1}} %полые цифры

%%%%%%%%%%%%%% жирная норма
% \[  \]
\def\[{\mathopen{\kern1pt \vrule height7.5pt depth2.4pt width1pt\kern1.5pt}}
\def\]{\mathclose{\kern1.5pt\vrule height7.5pt depth2.4pt width1pt\kern1pt}}
\def\sfatop{\mathopen{\kern1pt \vrule height5.5pt depth1.5pt width1pt\kern1pt}}
\def\sfatcl{\mathclose{\kern1pt\vrule height5.5pt depth1.5pt width1pt\kern1pt}}
%%%%%% большие
\def\leftbbr{\mathopen{\kern1pt \vrule height14pt depth8pt width1pt\kern1.5pt}}
\def\leftbbrr{\mathopen{\kern1pt\vrule height16pt depth13pt width1pt\kern1.5pt}}
\def\rightbbrr{\mathclose{\kern1.5pt\vrule height16pt depth13pt width1pt\kern1pt}}
\def\rightbbr{\mathclose{\kern1.5pt\vrule height14pt depth8pt width1pt\kern1pt}}
%%%%%%%%%%%%%%
\def\tvert{|\!|\!|} %тройная норма

\newcommand{\shortpage}{\enlargethispage{-\baselineskip}}

 \newcommand{\signatura}[2]{\enlargethispage{2\baselineskip}\par\begin{picture}(10,22)%
 \put(-12,#2){\makebox{#1}}
 \end{picture}}

 \newcommand{\Signatura}[3]{\begin{picture}(0,0)\put(#2,#3){\makebox{{#1}}}\end{picture}}

\newcommand{\be}[2]{\begin{equation}\label{#1} {#2}\end{equation}}
\newcommand{\bn}[1]{\[#1\]}
\newcommand{\ba}[2]{\begin{eqnarray}\label{#1} {#2}\end{eqnarray}}
\newcommand{\ban}[2]{\begin{eqnarray*}\label{#1} {#2}\end{eqnarray*}}
\newcommand{\lr}[1]{\left(#1\right)}
\newcommand{\sq}{\sqrt}
\newcommand{\ov}{\overline}
\newcommand{\pr}{\partial}

\Title{
Transmutations and Applications:\\
a survey.
}

\noindent\rule{0.4\linewidth}{0.3pt}

\noindent Originally published  in the book:

"Advances in Modern Analysis and Mathematical Modeling."\\

Editors: Yu.F.\,Korobeinik,  A.G.\,Kusraev.\\

Vladikavkaz: Vladikavkaz Scientific Center\\

of the Russian Academy of Sciences\\

and Republic of North Ossetia--Alania,\\

2008.\  P. 226--293.

\shorttitle{Transmutations and Applications.}

\Author{S.~M.~Sitnik}{Sitnik~S.~M.}

\addcontentsline{toc}{contpart}{\protect  Sitnik~S.~M.
{\rm Transmutations and Applications}}
\baselineskip=1.05\baselineskip

\begin{center}
Chair on   Mathematics,\\
Department of Radioengineering,\\
Voronezh Institute
of the Ministry of Internal Affairs of Russia,\\
Voronezh, Russia.\\

e-mail: pochtasms@gmail.com\\

postal: P.O.Box 12, Voronezh-5, \        Voronezh, \ 394005, \ Russia.
\end{center}

\medskip

In this survey we consider main   transmutation theory topics with many applications, including author's own results.
The topics covered are: transmutations for Sturm--Liouville operators,  Vekua--Erdelyi--Lowndes transmutations,
 transmutations for general differential  operators with variable coefficients, Sonine and Poisson transmutations,
transmutations and fractional integrals, Buschman--Erdelyi transmutations,
in the search for Volterra unitary transmutations, transmutations for singular differential  operators with variable coefficients,
composition method for transmutations, some applications and open problems.

\newpage
\centerline{\Large {Contents:}}

\begin{enumerate}
\item Introduction.
\item Transmutations for Sturm--Liouville operators.
\item Vekua--Erdelyi--Lowndes transmutations.
\item Transmutations for general differential  operators with variable coefficients.
\item Sonine and Poisson transmutations.
\item Transmutations and fractional integrals.
\item Buschman--Erdelyi transmutations.
\item In the search for Volterra unitary transmutations.
\item Transmutations for singular differential  operators with variable coefficients.
\item Composition method for transmutations.
\item Some applications for transmutations.
\item Problems.
\item References.
\end{enumerate}

\newpage

\UDC{517.444}
\begin{center}
Операторы преобразования и их приложения.\\
-----------------------------------------------------\\
ОБЗОР\\

В книге:\\ "Исследования по современному анализу и математическому моделированию"\\

Отв. ред. Коробейник Ю.Ф., Кусраев А.Г.\\

Владикавказ: Владикавказский научный центр\\

 Российской  Академии Наук\\ и Республики Северная Осетия---Алания,\\

 2008.\  Стр. 226--293.\\

\baselineskip=1.05\baselineskip
\medskip

В обзоре рассматриваются основные понятия и задачи теории операторов преобразования. Затем перечислены наиболее известные приложения операторов преобразования к обратным задачам, теории рассеяния, спектральной теории, нелинейным дифференциальным уравнениям и построению солитонов, обобщённым аналитическим функциям, сингулярным краевым задачам, теории дробного интегродифференцирования, вложениям некоторых функциональных пространств. Обзор заканчивается достаточно подробным списком литературы.
\end{center}
\newpage

\begin{center}
\vskip 3cm
ВАЖНАЯ ИНФОРМАЦИЯ :\\

Это версия опубликованного в 2007 обзора с существенными добавлениями 2008--2010 годов.
\end{center}

\newpage
\vskip 5.5cm
\begin{center}
ИНФОРМАЦИЯ ДЛЯ ЧИТАТЕЛЕЙ:\\
\end{center}

Уважаемый читатель! Это версия обзора, вышедшего в 2007 году. К сожалению, любой текст не может быть полным. Данный текст включает то, что было мне известно на момент написания в 2007 году, с минимальными последующими вставками. О некоторых результатах написано слишком кратко или поверхностно, хотя они этого не заслуживают из--за своей важности; о многих результатах я знал недостаточно или не знал вовсе. Приношу свои извинения авторам, чьи результаты не вошли в обзор или изложены недостаточно подробно.

Поздние вставки в первоначальный текст обзора привели к некоторому нарушению логики изложения и последовательности материала, сумбуру в библиографии, повторам и перекрытиям в тексте. Надеюсь, что это искупается включением дополнительной информации, а недостатки по возможности будут устранены в окончательном тексте.

Автор готовит обновлённую версию данного обзора для опубликования в 2010 году в виде монографии, она будет существенно расширена и дополнена. Многие результаты, которые не нашли отражения в первой версии, будут изложены в следующей.

Автор будет благодарен коллегам за любую критику и даже ругань (по делу), исправления, добавления, комментарии. Уже настоящий текст появился благодаря благожелательной помощи большого числа коллег, без которых данная работа никогда  не была бы выполнена.\\
\\

Автор: Ситник Сергей Михайлович,\\
доцент кафедры высшей математики\\
Воронежского института МВД, Воронеж, Россия.
\vskip 0.3cm

Электронный адрес: mathsms@yandex.ru
\vskip 0.3cm

Почтовый адрес: Ситник С.М.,\\
а.я. 12, Воронеж--5, Воронеж, 394005, Россия.
\newpage
\vskip 0.5cm

\centerline{\Large {Содержание обзора:}}

\begin{enumerate}
\item Введение.
\item Операторы преобразования для операторов Штурма--Лиувилля.
\item Операторы преобразования Векуа--Эрдейи--Лаундеса.
\item Операторы преобразования для общих дифференциальных операторов с переменными коэффициентами.
\item Операторы преобразования типа Сонина и Пуассона.
\item Связь операторов преобразования с дробным интегродифференцированием.
\item Операторы преобразования Бушмана--Эрдейи.
\item В поисках вольтерровых унитарных операторов преобразования.
\item Операторы преобразования для некоторых сингулярных дифференциальных операторов с переменными коэффициентами.
\item Композиционный метод  построения операторов преобразования различных классов.
\item Некоторые приложения метода операторов преобразования.
\item Задачи.
\item Литература.
\end{enumerate}

\newpage
\section{1. Введение.}
\textsc{Определение 1.} Пусть дана пара операторов $(A,B)$. Оператор
$T$ называется \textit{оператором преобразования} (ОП, transmutation), если
выполняется соотношение
\begin{equation}
\label{1.1}
\Large {T\,A=B\,T.}
\end{equation}

Соотношение (\ref{1.1}) называется иначе \textit{сплетающим свойством}, тогда
говорят, что ОП $T$ \textit{сплетает} операторы $A$ и $B$ (intertwining operator). Для превращения
(\ref{1.1}) в строгое определение необходимо задать пространства или
множества функций, на которых действуют операторы $A$, $B$, и,
следовательно, $T$. Иногда в определение ОП закладывают и
требование обратимости, что является желательным, но не обязательным
свойством. В конкретных реализациях операторы $A$ и $B$ обычно
являются дифференциальными, $T$ --- линейный оператор на стандартных
пространствах. Ясно, что понятие ОП является прямым и далеко идущим обобщением понятия подобия матриц из линейной алгебры [1--3]. Но ОП \textit{не сводятся к подобным (или эквивалентным) операторам},  так как сплетаемые операторы как правило являются неограниченными в естественных пространствах, к тому же обратный к ОП не обязан существовать и  действовать в том же пространстве. Так что спектры операторов, сплетаемых ОП, могут не совпадать.

Как же обычно используются операторы преобразования? Пусть, например,  мы изучаем некоторый достаточно сложно устроенный оператор $A$. При этом нужные свойства уже известны для модельного более простого оператора $B$. Тогда, если существует ОП (\ref{1.1}), то часто удаётся перенести свойства модельного оператора  $B$ и на $A$. Такова в нескольких словах примерная схема типичного использования ОП в конкретных задачах.

В частности, если рассматривается уравнение $Au=f$ с оператором $A$, то применяя к нему ОП $T$ со сплетающим свойством (\ref{1.1}), получаем уравнение с оператором $B$ вида $Bv=g$, где обозначено $v=Tu, g=Tf$. Поэтому, если второе уравнение с оператором $B$ является более простым, и для него уже известны формулы для решений, то мы получаем и представления для решений первого уравнения $u=T^{-1}v$. Разумеется, при этом обратный оператор преобразования должен существовать и действовать в рассматриваемых пространствах, а для получения явных представлений решений д
 лжно быть получено и явное представление этого обратного оператора. Таково одно из основных применений техники ОП в теории дифференциальных уравнений с частными производными.

Изложению теории ОП и их приложениям посвящены существенные части монографий [4--9], кроме того различные вопросы ОП рассматриваются также в [10--41]. К сожалению, на русском языке нет книг, полностью посвящённых  ОП, таких, как замечательные книги Роберта Кэррола на английском [4--6]. По--видимому, данный обзор также является первым на русском языке по теории ОП. Не претендуя на полноту, я включил в обзор только те результаты, которые представляются основными, а также связанные с собственными работами.

Сделаем одно терминологическое замечание. В западной литературе принят для ОП термин 'transmutation', восходящий к Ж.~Дельсарту. Как отмечает Р.~Кэрролл, похожий термин 'transformation' при этом закрепляется за классическими интегральными преобразованиями Фурье, Лапласа, Меллина и другими подобными им. Кроме того, термин 'transmutation' имеет в романских языках дополнительный оттенок 'волшебного превращения', что довольно точно характеризует действие ОП. Приведём точную цитату из [6]: "Such operators are often called transformation operators by the Russian school (Levitan, Naimark, Marchenko et. al.), but transformat
 ion seems too broad a term, and, since some of the machinery seems "magical"\, at times, we have followed Lions and Delsarte in using the word transmutation".

Необходимость теории ОП доказана большим числом её приложений. Методы ОП применяются в теории обратных задач, определяя обобщённое преобразование Фурье, спектральную функцию и решения знаменитого уравнения Гельфанда--Левитана; в теории рассеяния через ОП выписывается не менее знаменитое уравнение Марченко; в  спектральной теории получаются известные формулы следов и асимптотика спектральной функции; оценки ядер ОП отвечают за устойчивость обратных задач и задач рассеяния; в теории нелинейных дифференциальных уравнений метод Ла
 са использует ОП для доказательства существования решений и построения солитонов. Определёнными разновидностями ОП являются части теорий  обобщённых аналитических функций и операторов обобщённого сдвига. В теории уравнений с частными производными методы ОП применяются для построения явных решений некоторых  задач, изучении сингулярных и вырождающихся краевых задач, псевдодифференциальных операторов,  задач для решений с особенностями на части границы, оценки скорости убывания решений некоторых эллиптических и ультраэллиптич
 еских уравнений. Теория ОП позволяет дать некоторую новую классификацию специальных функций и интегральных операторов со специальными функциями в ядрах, в том числе различных операторов дробного интегродифференцирования. В теории функций найдены приложения ОП к вложениям функциональных пространств и обобщению операторов Харди, расширению теории Пэли--Винера.  Методы ОП с успехом применяются во многих прикладных задачах: оценках решений Йоста и квантовой теории рассеяния, исследовании системы Дирака, теории вероятностей и случа
 ных процессов, линейном стохастическом оценивании, фильтрации, стохастических случайных уравнениях, обратных задачах геофизики и трансзвуковой газодинамики.

Методы теории ОП и связанные с ними задачи в той или иной степени применялись в работах многих математиков. Перечислим некоторых из них:  A.~Boumenir, H.~Begehr, J.~Betancor, B.~Braaksma, L.~Bragg, R.~Carroll, M.~Chao, H.~Chebli, I.~Dimovski, C.~Dunkl, J.~Delsarte,
R.~Gilbert, V.~Hutson, R.~Hersh, V.Kiryakova, J.~L\H{o}ffstr\H{o}m, J.~Lions, A.~Mercer, J.~Peetre, J~.S~Pym, F.~Santosa, J.~Siersma, H.~deSnoo, K.~Stempak, V.~Thyssen,  K.~Trimeche, Vu Kim Tuan,
Агранович~З.~С., Андрощук~А.~А., Бритвина~Л.~Е., Буслаев~В.~С., Валицкий~Ю.~Н., Волк~В.~Я., Глушак~А.~В., Горбачук~М.~Л.,
Гаджиев~А.~Д., Гринив~Р.~О., Гулиев~В.~С., Житомирский~Я.~И.,
Иванов~Л.~А., Ерёмин~М.~С., Карп~Д.~Б., Катрахов~В.~В., Качалов~А.~П., Килбас~А.~А,
Киприянов~И.~А., Ключанцев~М.~И., Кононенко~В.~И., Коробейник~Ю.~Ф., Кулиш~П.~П., Курылёв~Я.~В.,
Лаптев~Г.~И, Левин~Б.~Я., Левитан~Б.~М., Ляхов~Л.~Н., Ляховецкий~Г.~В.,
Маламуд~М.~М., Марченко~В.~А., Мацаев~В.~И., Микитюк~Я.~В., Муравник~А.~Б., Наймарк~М.~А., Нагнибида~Н.~И., Олевский~М.~Н., Платонов~С.~С.,
Повзнер~А.~Я., Поляцкий~В.~Т., Сахнович~Л.~А., Ситник~С.~М., Сохин~А.~С., Сташевская~В.~В., Фаддеев~Л.~Д., Фаге~Д.~К., Фомин~В.~Л., Хачатрян~И.~Г., Хромов~А.~П., Шмулевич~С.~Д., Ярославцева~В.~Я.
Разумеется, этот список не полон и может быть существенно расширен.

Специально отметим  вклад математиков Харьковской школы в развитие теории операторов преобразования, который очень значителен: В.А.~Марченко, Н.И.~Ахиезер, Б.М.~Левитан, В.В.~Сташевская, А.Я.~Повзнер, Я.И.~Житомирский, А.С.~Сохин, Л.А.~Сахнович, В.Я.~Волк и другие.

\newpage

\section{2. Операторы преобразования для операторов Штурма-Лиувилля.}

Рассмотрим задачу построения различных ОП, сплетающих простейший оператор Штурма-Лиувилля
\begin{equation}
\label{2.1}
y''(x)+\lambda^2y(x)=(L_0y)(x)
\end{equation}
c оператором Штурма-Лиувилля общего вида
\begin{equation}
\label{2.2}
y''(x)+q(x)y(x)+\lambda^2y(x)=(Ly)(x).
\end{equation}
Функция $q(x)$ в (\ref{2.2}) называется потенциальной функцией.
Эта функция может быть комплексной, число $\lambda\in \mathbb{C}$, переменная $x\in \mathbb{R}$.
Мы ищем ОП, удовлетворяющий тождеству
\begin{equation}
\label{2.3}
SLf=L_0Sf
\end{equation}
на подходящих функциях $f(x)$ и при данных $q(x)$ и $\lambda$.
Естественные требования к ОП $S$ --- линейность и обратимость в стандартных
пространствах. Требование линейности $S$ после подстановки (\ref{2.1})
и (\ref{2.2}) в (\ref{2.3}) приводит к соотношению, которое не зависит
от $\lambda$
\begin{equation}
\label{2.4}
S\left(D^2+q(x)\right)f=D^2Sf,
\end{equation}
где обозначено $D=d/dx$. Требование обратимости естественным образом
приводит к поиску ОП $S$ в виде интегрального оператора
\begin{equation}
\label{2.5}
(Sf)(x)=\int\limits_{a(x)}^{b(x)}K(x,t)f(t)\,dt
\end{equation}
причем ядро $K(x,t)$ в общем случае может быть распределением
(например, $K(x,t)=\delta(x-t)+G(x,t)$, $G$ --- гладкая функция).
В формуле (\ref{2.5}) $a(x)$, $b(x)$ --- некоторые функции $\mathbb{R}\to \mathbb{R}$.

Существуют два несколько различных подхода к построению ОП. В первом они строятся
лишь на решениях уравнений $L_0y=0$, $Ly=0$ с операторами (\ref{2.1})--(\ref{2.2}).
Такой метод был принят в пионерских классических работах по ОП, он и излагается во всех книгах. Похоже, что во многом ему следуют только благодаря традиции. Но возможен и второй подход, который мне кажется более естественным и общим, когда ОП строятся на произвольных функциях, возможно с некоторыми ограничениями роста в фиксированных точках. Мы кратко изложим именно этот подход, следуя [41], он не является общепринятым. Ясно, что при нашем подходе класс исследуемых ОП побогаче. С другой стороны, оператор, построенный на собственных функци
 ях
при любых $\lambda$, может быть распространён на достаточно широкие классы  благодаря полноте систем собственных функций во всех основных пространствах. Поэтому за исключением сужения  множества рассматриваемых ОП и первый общепринятый подход не создаёт серьёзных ограничений в приложениях.

Заметим, что мы будем для краткости с некоторой вольностью называть операторами
дифференциальные или интегральные выражения, не всегда указывая пространства, в которых действуют операторы. При построении ОП
предполагается, что функции $f(x)$ принадлежат некоторому классу $\Phi$,
и что ядра $K(x,t)$ обладают определенной
гладкостью по обеим переменным.

Из вышеизложенного ясно, что наиболее целесообразно искать ОП $S$ как интегральные операторы Вольтерра второго рода:
\begin{equation}
\label{tag8}
(Sf)(x)=f(x)+\int\limits_{c}^{x}K(x,t)f(t)\,dt,
\end{equation}
Здесь $K(x,t)$---гладкая функция, число $c$ принадлежит
расширенной числовой оси $\mathbb{\overline{R}}$. Вместе с тем в [41] рассмотрены и другие возможные операторы.
Операторы Вольтерра (\ref{tag8}) легко обращаются в стандартных
пространствах. Свобода выбора предела интегрирования позволяет в каждом
конкретном случае использовать те ОП, которые наилучшим образом подходят для
данной задачи. В частности, при $c=0$ получаются ОП, впервые построенные
А.~Я.~Повзнером  и Б.~М.~Левитаном [42--43] (см. также [44]), а при $c=+\infty$ получаются ОП, впервые построенные
Б.~Я.~Левиным [45], последние  сохраняют асимптотику на бесконечности и используются
при решении обратных задач квантовой теории рассеяния.

Для ядра ОП (\ref{tag8}) в результате несложных вычислений получаем
\begin{eqnarray}
\frac{\partial^2 K}{\partial t^2}+q(t)K=
\frac{\partial^2}{\partial x^2} K,\label{tag9}\\
\frac{d}{dx}K(x,x)=\frac12q(x),\quad
(x,t)\in\Omega\subset \mathbb{R}^2,\label{tag10}\\
\lim_{t\to c}W_t\left(f,K(x,t)\right)=\lim_{t\to c}W(x)=0.\label{tag11}
\end{eqnarray}
(Через $\Omega$ обозначена область определения функции $K(x,t)$,
замыкание которой содержит часть диагонали $t=x, W_t$---определитель Вронского с производными по t.)

Уравнение (\ref{tag9}) является стандартным в той задаче, которую мы
рассматриваем. Условие (\ref{tag11}) показывает выделенность точки $c$.
Оно сводится к весовым граничным условиям на функцию $f(x)$ и ее
первую производную $f'(x)$ при $x\to c$ и диктует
выбор класса функций $\Phi$.
Важнейшее равенство (\ref{tag10}) устанавливает, что значение ядра ОП $K(x,t)$ на
диагонали $t=x$ позволяет восстановить потенциал $q(x)$. Этот факт является
основным в теории обратных задач. Поэтому наиболее распространенные методы
решения обратных задач сводятся именно к нахождению ядра ОП по спектральной функции (как в уравнении Гельфанда--Левитана, полученном Б.М.~Левитаном)
или по данным рассеяния (как в уравнении Марченко).

Чтобы построить решения системы (\ref{tag9}--\ref{tag11}), сначала нужно зафиксировать класс  $\Phi$ так, чтобы выполнялось условие (\ref{tag11}), например : $f(x)\in C^2, f(c)=f'(c)=0$.
Далее обычно реализуется следующий план действий. На первом шаге переходят от
уравнения в частных производных к интегральному. Эти уравнения не
эквивалентны, однако каждое решение интегрального уравнения удовлетворяет
исходному гиперболическому. На этом шаге доказывается существование
некоторого ядра ОП и его определенная гладкость. На втором шаге
проверяются  дополнительные условия для ядра и выбирается подходящий
класс функций $\Phi$. Этим заканчивается построение ОП.

Итак, перейдем к решению уравнения
\begin{equation}
\label{tag12}
\frac{\partial^2 K}{\partial x^2}=\frac{\partial^2 K}{\partial t^2}+q(t)K
\end{equation}
с дополнительным условием на диагонали $x=t$ (см. (\ref{tag10})),
\begin{equation}
\label{tag13}
K(x,x)=\frac12q(x).
\end{equation}
Выполним стандартную замену переменных по формулам
\begin{equation}
\label{tag14}
u=\frac12(x+t),\quad v=\frac12(x-t).
\end{equation}
При этом уравнение диагонали $x=t$ в новых переменных примет вид $v=0$.

Введем обозначение для ядра в новых переменных
\begin{equation}
\label{tag15}
H(u,v)=K(u+v,u-v)=K(x,t).
\end{equation}
Для функции $H$ перейдем от соотношений (\ref{tag12})--(\ref{tag13}) к новым
\begin{eqnarray}
&\frac{\partial^2 H}{\partial u\partial v}=q(u-v)H,\label{tag16}\\
&H(u,0)=\frac12\int\limits_c^uq(s)\,ds.\label{tag17}
\end{eqnarray}
Здесь функция $q(u-v)$ должна быть задана, $c$ --- произвольное число,
возможно $c=\pm\infty$.

Необходимо сделать важное замечание. Из формул замены (\ref{tag14})
следует, что возможно как $u>0$, $v>0$, так и $u<0$, $v<0$.

Система (\ref{tag16})--(\ref{tag17})---это задача Коши лишь с одним начальным
условием. Поэтому при естественных предположениях на $q$ ($q(x)\in C^1$) эта система
имеет бесконечно много решений. Следовательно, при каждом потенциале $q(x)$
существует бесконечно много ОП, например, вида (\ref{tag8}).
Это чрезвычайно удобно в приложениях, где можно выбирать при одном и том же
потенциале различные ОП, наиболее подходящие для каждой конкретной задачи.

Приведем 'нечестный' путь построения некоторого класса ядер, удовлетворяющих
системе (\ref{tag16})--(\ref{tag17}). Для этого достаточно убедиться, что
каждое $C^2$-решение интегрального уравнения
\begin{equation}
\label{tag18}
H(u,v)=\frac12\int\limits_c^uq(s)\,ds\,+\int\limits_d^u\,d\alpha\,
\int\limits_0^vq(\alpha-\beta)H(\alpha,\beta)\,d\beta
\end{equation}
удовлетворяет системе (\ref{tag16})--(\ref{tag17}) для
произвольных чисел $c$, $d\in \overline{R}$.

'Честный' вывод уравнения (\ref{tag18}) проводится по известной стандартной схеме с использованием
функции Римана [46--50]. Поэтому для теории ОП важное значение имеет нахождение функций Римана для конкретных уравнений в явном виде [47--66]. В этом направлении выделим работы Куйбышевско--Самарской математической школы под руководством С.~П.~Пулькина и В.~Ф.~Волкодавова [51--57].

Важность исследования уравнения (\ref{tag18}) с разными $c$, $d$ в том,
что мы можем одновременно изучать и случай $c=d=0$, возникающий при
построении ОП типа Повзнера--Левитана, и случай $c=d=+\infty$, возникающий
при построении ОП типа Левина, и случай произвольных различных $c$, $d$.
Обычно эти типы ОП изучались раздельно, а случай произвольных $c$, $d$
($c,d\neq 0$, $c,d\neq\pm\infty$) не рассматривался.

Доказательство существования решения для интегрального уравнения (\ref{tag18}) не представляет труда, оно проводится методом последовательных приближений. При этом решение представляется рядом Неймана [67]. Например, для случая $c=d=0$ и суммируемого потенциала получается такой типичный результат [4--6, 22--23].

\Theorem{ 1} Пусть функция $q(x)\in C^k$, $c=d=0$. Тогда уравнение  (\ref{tag18}) имеет единственное решение, удовлетворяющее оценке
\begin{equation}
\label{tag19}
\vert K(x,t)\vert \leq \frac{1}{2}w(\frac{x+t}{2})\exp\left( \sigma_1(x)-\sigma_1(\frac{x+t}{2})-\sigma_1(\frac{x-t}{2})\right),
\end{equation}
где введены обозначения
\begin{equation}
\label{tag20}
w(u)=\max_{0\leq s\leq u}\vert \int_0^s q(y)\,dy\vert, \sigma_0(x)=\int_0^x \vert q(t)\vert\,dt,\sigma_1(x)=\int_0^x \sigma_0(t)\,dt.
\end{equation}
Если функция $q(x)$ имеет $k\geq 0$ непрерывных производных, то ядро ОП $K(x,t)$ имеет $k+1$ непрерывную производную по обеим переменным.
\Endproc

Аналогично в этом случае доказывается  существование и единственность обратного оператора к ОП (\ref{tag8}) в пространствах $C^2$ или $L^2$. Для ядра обратного оператора доказываются аналогичные теореме 1 оценки, гладкость ядра, соотношение между ядрами прямого и обратного операторов. Заметим, что обратный к (\ref{tag8}) $P=S^{-1}$ также является ОП и действует по правилу $PL_0f=LPf$, см. [22--23].

Общий случай произвольных $c,d$ разобран в [41]. Там же явно построены ядра для некоторых простейших потенциалов, в основном выражающиеся через функции Бесселя и Макдональда. Получены оценки ядер для потенциалов, удовлетворяющих различным равномерным неравенствам.

Оценку (\ref{tag19}) можно уточнять разными способами, при этом ослабляя условия на потенциал, см. [4--6, 21--23]. Наиболее сильные результаты в этом направлении при минимальных требованиях к потенциалу можно получить, используя методы очень интересной книги [69]. В ней В.~А.~Чернятиным получены окончательные необходимые и достаточные условия на потенциал для разрешимости смешанной задачи для гиперболического уравнения (\ref{tag9}), что завершило исследования В.~А.~Стеклова, И.~Г.~Петровского, Б.~М.~Левитана, В.~П.~Михайлова, В.~А.~Ильина по этой задаче; при э
 ом использовался метод Фурье и принадлежащая А.~Н.~Крылову идея выделения особенностей из ряда. Метод работы [69] очень близок идейно методам теории ОП, он и был приспособлен для доказательства существования ОП с не очень хорошим потенциалом в [41].  Кроме того, в [41] получена оценка невязки
$SLf-L_0Sf$ в соотношении (\ref{2.3}) для того случая, когда ОП $S$ задан с ошибкой, например, образовавшейся от замены ряда Неймана при решении интегрального уравнения конечной суммой. Результаты В.~А.~Чернятина другим способом, который представляется более простым, были впоследствии получены В.Л.Прядиевым [389].

Вместе с тем, интересным как для теории, так и для практических приложений является рассмотрение операторов Штурма--Лиувилля с совсем плохими потенциалами, вплоть до распределений. Результаты в этом направлении с приложениями к спектральной теории и формулам следов получены в [70--74] А.~А.~Шкаликовым и его учениками . Некоторые другие результаты для подобных ОП см. в [75--77]. Разумеется, что это направление также содержит работы, в которых существенно используются методы теории ОП; так в [391--393] Гринив и Микитюк изучают построение ОП для плохи
 х потенциалов и их последующее применения для решения обратных задач.

Вместе с тем ОП определённого вида существуют не всегда. Ряд отрицательных результатов получен в [41]. Так не существует ОП в форме оператора Фредгольма второго рода с гладким ненулевым потенциалом. Не существует линейного ОП, сплетающего нелинейный оператор $D^2 y +y^2$ со второй производной. Для ОП (\ref{tag8}) в случае простейших потенциалов $q(x)=x^2, \cos x$ ядра не могут быть аналитическими функциями. Этим объясняется, например, тот факт, что для функций Матье (решений уравнения $D^2y +\cos(x) y=\lambda y$) нет простых интегральных выражений через тригонометри
 еские или экспоненциальные функции.

В заключение этого пункта отметим такой удивительный факт: во всех известных книгах и работах, начиная с возникновения теории операторов преобразования, сразу рассматривались и строились ОП для дифференциальных операторов \textit{второго порядка}. Но наиболее логично начать построение теории ОП с рассмотрения операторов \textit{первого порядка} вида $D-q(x)$. В этом случае соответствующие ОП  должны  удовлетворять  сплетающим соотношениям
$$
T_1 \lr{D-q(x)}=DT_1, \ \ \ \ \ T_2 D=\lr{D-q(x)} T_2 .
$$
Этот пробел в литературе устранён в [41], где построены соответствующие ОП. Открытым остаётся такой интересный вопрос: можно ли из этих ОП для дифференциальных операторов \textit{первого порядка} сконструировать ОП для дифференциальных операторов \textit{второго порядка}? Вопрос остаётся открытым даже в случае простейшего постоянного потенциала $q(x)=a^2$, хотя оператор Штурма--Лиувилля при помощи известных замен разлагается на множители первого порядка.

\newpage

\section{3. Операторы преобразования Векуа--Эрдейи--Лаундеса.}

Частный случай ОП для оператора Штурма--Лиувилля (\ref{2.4}) получается при выборе  постоянного потенциала. При этом возникает важный класс ОП Векуа--Эрдейи--Лаундеса, которые осуществляют сдвиг по спектральному параметру.

\textsc{Определение 2.} Обобщённым оператором преобразования Векуа
-- Эрдейи -- Лаундеса (ВЭЛ) называется сплетающий оператор для пары
$(A+\lambda_1,A+\lambda_2)$, где $A$ -- некоторый базовый оператор,
$\lambda_1,\lambda_2$ -- комплексные числа.

Иными словами
\begin{equation}
\label{31}
T(A+\lambda_1)=(A+\lambda_2)T.
\end{equation}

ОП ВЭЛ были по разным поводам введены и изучены в работах И.Н.~Векуа [78--81],
А.~Эрдейи [82--84] и Дж.~С.~Лаундеса [85--87]. В их работах
рассматривались такие базовые дифференциальные операторы:
\begin{equation}
A=D^2=\frac{d^2}{dx^2},\  A=B_{\nu}=D^2+\frac{2\nu +1}{x}D,
A=x^{\beta}B_{\nu} \ .
\end{equation}

Первый ОП ВЭЛ был построен Ильёй Несторовичем Векуа в виде
\begin{equation}
\label{33}
J_{\lambda}f(x)=f(x)-\int_0^x t\frac{J_1(\lambda
\sqrt{x^2-t^2})}{\sqrt{x^2-t^2}}\ f(t)\ dt,
J_{\lambda}(D^2+\lambda)=D^2 J_{\lambda},
\end{equation}
где $J_1(\cdot)$ --- функция Бесселя. Такой ОП может быть
использован, например, для представления решений телеграфного
уравнения через решения волнового. В этом проявляется основная роль, которую ОП ВЭЛ играют в теории дифференциальных уравнений в частных производных---они осуществляют сдвиг по спектральному параметру. Поэтому с их помощью получаются представления решений для задач 'с лямбда' через решения задач 'без лямбда', многочисленные примеры можно найти в энциклопедической монографии [34] и ссылках в ней (см. особенно главу 8).

Для А.~Эрдейи основным было изучение свойств самих операторов вида (\ref{33}), их композиций и оценок норм, тогда как Дж.~Лаундес в основном изучал представление решений различных задач для уравнений с частными производными при помощи подобных операторов. Мотивацией И.~Н.~Векуа для введения ОП было выразить решения уравнения $\Delta u+\lambda u=0$, к которому сводится простейшая задача теории упругости, через гармонические функции. Затем он рассмотрел и более сложные уравнения, в основном возникающие в механике. Простейший подобный интегральный оп
 ратор, переводящий аналитические функции в гармонические, был уже давно к тому времени известен, его ввёл Э.~Уиттекер.

Так возникла новая теория обобщённых аналитических функций, начинающаяся с построения операторов, выражающих решения различных сложных уравнений через аналитические или гармонические функции. Определение таких обобщённых аналитических функций может быть дано на языке ОП: в действительном случае  они возникают при рассмотрении ОП, которые сплетают некоторый дифференциальный оператор с частными производными с оператором Лапласа или его степенями, а в комплексном случае они возникают при рассмотрении ОП, которые сплетают уравнен
 я (оператор) Коши--Римана (=Эйлера--Деламбера) с его обобщёнными аналогами. После Векуа, Эрдейи и Лаундеса теория обобщённых аналитических функций именно в разрезе явного построения операторов преобразования достраивалась в работах известного американского математика Стефана Бергмана и его учеников (операторы Бергмана--Гилберта) [88, 7, 10, 14]. При этом для трёх и более переменных всё становится намного сложнее и интереснее. Отметим, что работы Бергмана существенно используют теорию специальных функций, в частности, гипергеометрических.
  Приложения к физическим задачам рассматривались Л.~Берсом, А.~Гельбартом, А.~Вайнштейном (GASPT---теория обобщённого осесимметрического потенциала) [89--93]. Следует отметить также работы [94, 37], в которых Г.~Н.~Положий  рассмотрел обобщения уравнений Коши--Римана, решения которых названы им $p$--аналитическими и $(p,q)$--аналитическими функциями. Принадлежность обобщённых аналитических функций различным пространствам (Харди, Смирнова, ВМО) изучается в работах С.~Б.~Климентова [95--98].
Существенные результаты для теории обобщённых аналитических функций были также получены в работах Б.~Боярского, Н.К.~Блиева, Л.Г.~Михайлова, А.~П.~Солдатова и других.

ОП ВЭЛ изучались также в [99--101], перечислим некоторые полученные там результаты. Вначале выделено семейство из восьми основных операторов ВЭЛ. Изучены факторизации этих операторов через более простые: Фурье,
Ханкеля, дробные интегралы Римана--Лиувилля, Эрдейи--Кобера и
другие [34]. Изучены полугрупповые свойства введённых операторов ВЭЛ по
параметру. Описаны общие методы построения ОП ВЭЛ из уже известных. На этом
пути получены новые операторы ВЭЛ, ядра которых выражаются через
гипергеометрические функции от нескольких переменных: Райта, Фокса, Гумберта, Аппеля,
Кампе де Ферье и другие [102--120]. В [41, 99--101] найдены новые семейства ОП, зависящие от произвольных параметров.

Также построено взаимно однозначное соответствие между ОП ВЭЛ и ОП, сплетающими весовые операторы Бесселя с разными индексами:
\begin{equation}
\label{32}
T(x^2 B_\nu)=(x^2 B_\mu)T, \  B_{\nu}=D^2+\frac{2\nu +1}{x}D,
\end{equation}
Такие ОП изучались в [121], где получен один  класс подобных операторов. В [99] найдены формулы, решающие задачу об описании ОП со свойством (\ref{32}) в общем виде. Оказывается, что по каждому ОП ВЭЛ со свойством  (\ref{31}) можно явно выписать ОП со свойством  (\ref{32}).

Таким образом, теория операторов преобразования Векуа---Эрдейи---Лаундеса представляется сформировавшимся разделом общей теории ОП с
достаточно большим набором собственных результатов и важными
применениями в теории уравнений с частными производными, теории функций, комплексном анализе, теории специальных функций. Замечательно,
что начало этой плодотворной тематики было заложено в трудах академика Ильи
Несторовича Векуа, столетие со дня рождения которого отмечалось в 2007 году.
\newpage

\section{4. Операторы преобразования для общих дифференциальных операторов с переменными коэффициентами.}

Теперь перейдём к рассмотрению вопроса о существовании ОП $T$, который сплетает два произвольных линейных дифференциальных оператора

\begin{eqnarray}
\label{41}
T \left( a_n(x)y^{(n)}+\cdots +a_1(x)y+a_0(x)\right)= \\
=\left(b_m(x)y^{(m)}+\cdots +b_1(x)y+b_0(x)\right)T.\nonumber
\end{eqnarray}

В этом случае задача существенно усложняется, не все построения оказываются теоретически возможными, что напоминает теорию разрешимости алгебраических уравнений высших порядков.

По--видимому, впервые  задачу о построении ОП вида (\ref{41}) в том случае, когда в (\ref{41}) справа стоит оператор дифференцирования $n$--го порядка,  в общих чертах сформулировал Жан Дельсарт (член группы Бурбаки) в 1938 г [122](обратим внимание, что во всех русских книгах в ссылках на эту статью неверно указаны номера страниц, то есть скорее всего ни один автор её не читал!). Затем общее обсуждение продолжилось, в том числе совместно с Ж.~Л.~Лионсом [123--126]. При этом роль этих двух французских математиков для теории ОП на мой взгляд совершенно различна. Раб
 оты Дельсарта содержат большое число новых интересных результатов (см. ниже). Работы же Лионса по данной тематике содержат набор очевидных рассуждений, тривиальных выкладок общего плана или неверных результатов. (Вспоминаю, как мучил нас, аспирантов, И.~А.~Киприянов, заставляя разбирать книги Лионса по квазиобращению и Лионса--Мадженеса по краевым задачам, а также Лере по преобразованию Лапласа. Основные черты этого стиля таковы: написано плохо, наскоро, без всякого уважения к читателю; доказательства содержат только очевидное, в перв
 м более--менее сложном месте следует ссылка на другой источник, часто неверная; цитируются в основном только свои работы, в крайнем случае других французов. Я не стал бы приводить эти злобные выпады здесь, но, к сожалению, манера написания работ нынешних французских математиков, как следует из отзывов живущего и работающего там В.~И.~Арнольда [127], нисколько не изменилась. Казалось бы, что мешает им брать пример  с Лагранжа, Пуассона, Пуанкаре, прекрасные работы которых пока что еще хранятся в библиотеках французских университетов, не гов
 ря уже об авторах  великой французской художественной литературы! Зато какое удовольствие после всей этой дребедени было прорабатывать наши учебники И.~Г.~Петровского, С.~Л.~Соболева, С.~Г.~Михлина, Ю.~В.~Егорова по УЧП~!).

Вернёмся к ОП. Оказалось, что для построения ОП в случае дифференциальных операторов порядков выше двух существенную роль неожиданно начинает играть аналитичность коэффициентов. Отметим, что для случая дифференциального оператора второго порядка аналитичность коэффициентов практически вообще не играет никакой роли в задачах о построении ОП. Кроме того, следует обратить внимание, что на самом деле существуют два разных подхода к построению ОП вида (\ref{41}). При первом подходе рассматриваются \textit{локальные ОП}, которые в разных точках
  области определения задаются разными выражениями, они определены локально и меняются от точки к точке. При втором подходе строятся \textit{глобальные ОП}, которые на всей области определения задаются одним выражением, например, вида (\ref{2.5}) с некоторым ядром. Нужно понимать, что это две совершенно разные задачи, по методам их решения как оказалось существенно различные и схожие лишь формулировкой. Например, ответ на вопрос о самом существовании глобальных или локальных ОП принципиально отличается.
Данный факт часто недооценивается или игнорируется, что приводит к принципиальным ошибкам.

Вопрос о существовании локальных ОП вида (\ref{41}) был впервые полностью решён М.~К.~Фаге в 1957--58 гг. ([128--132], подробное изложение в монографии [9]). В литературе часто цитируют как авторов этого результата Дельсарта и Лионса, однако их результат является неверным и содержит принципиальные неустранимые ошибки, (см. ниже). Основной вывод таков: дифференциальные операторы вида (\ref{41}) с непрерывными коэффициентами одного порядка с единичными старшими коэффициентами \textit{локально эквивалентны} на множестве аналитических функций. Эквивалентность
 означает, что на аналитических пространствах функций было доказано не только существование локальных ОП, но и их непрерывная обратимость.  Основные результаты были сформулированы и доказаны М.~К.~Фаге в рамках созданной им теории $(M,L)$--аналитических функций. Интересно отметить, что при этом по существу даже не использовались методы теории уравнений с частными производными. Существенные результаты в данном направлении были получены В.А.~Марченко.

Теперь рассмотрим \textit{другую задачу} о построении \textit{глобального} ОП в виде интегрального оператора Вольтерра второго рода с некоторым ядром
\begin{equation}
\label{42}
(Tf)(x)=f(x)+\int\limits_{0}^{x}K(x,t)f(t)\,dt.
\end{equation}
Как показали впервые Л.~А.~Сахнович и В.~И.~Мацаев, это не всегда возможно уже для простейших дифференциальных операторов третьего порядка [133--136].

\Theorem{ 2} Если существует ОП вида (\ref{42}), сплетающий операторы третьего порядка по формуле
\begin{equation}
\label{43}
T(D^3-q(x))f=(D^3)Tf, f\in C^2(0,a), f'(0)=f''(0)=0,
\end{equation}
то потенциал $q(x)$ является аналитической функцией на $(0,a)$.
\Endproc

Аналитичность в этом результате строго по существу. С одной стороны, Л.~А.~Сахновичем было показано, что для аналитических потенциалов  $q(x)$ задача о построении искомого ОП всегда разрешима [136]. С другой стороны, было доказано, что для чуть худших не бесконечно дифференцируемых $q(x)$, искомый ОП для операторов четвёртого порядка $ T(D^4-q(x))f=(D^4)Tf$  может не существовать [135--136]; переработанный подобный пример для оператора третьего порядка с кусочно--постоянным потенциалом приведён в [40]. Таким образом, получается, что свойство наличия ОП может бы
 ть устранено или добавлено малым шевелением коэффициентов, это даёт ответ на вопрос В.~Б.~Лидского [142]. О необходимых условиях существования ОП см. также [137--139].

В этом направлении для дифференциальных операторов общего вида окончательные результаты были получены в работах М.~М.~Маламуда [140--144]. В них рассмотрен вопрос о существовании ОП на решениях простейшей задачи Коши
\begin{eqnarray}
\label{44}
u^{(n)}(x,\lambda) - \lambda^n u(x,\lambda)=0, u^{(k)}(0,\lambda)=a_k,0\leq k\leq n-1,
\end{eqnarray}
и  общей задачи Коши с переменными коэффициентами
\begin{eqnarray}
\label{441}
v^{(n)}+q_1(x) v^{(n-2)}+\cdots+q_{n-1}(x) v-\lambda^n v=0,\\ \nonumber
v^{(k)}(0,\lambda)=b_k,0\leq k\leq n-1.
\end{eqnarray}

Будем кратко говорить, что существует ОП для уравнения (\ref{441}), если для его решений  существует представление вида (\ref{42}) через решения (\ref{45}).

Сначала приведём результаты М.~М.~Маламуда для двучленного уравнения, а затем для общего случая.

\Theorem{ 3} Если для уравнения $v^{(n)}+ q(x)v-\lambda^n v=0$ при $n\geq 3$ существует ОП, то функция  $q(x)$ является аналитической.
\Endproc

\Theorem{ 4} Пусть $a\leq \infty$ и $q_k (x)$ аналитична на $(0,a)$ при $1\leq k \leq [\frac{n}{2}]$. Тогда, если для уравнения (\ref{441}) существует ОП, то функции  $q_k(x)$ являются аналитическими и при $[\frac{n}{2}]+1\leq k \leq n-1$.
\Endproc

Эти результаты уже используют многие факты теории уравнений с частными производными: метод функций Римана, разрешимость строго эллиптических по Лопатинскому уравнений, результаты Агмона по разрешимости нелинейных  уравнений с эллиптическим оператором в главной части, уточнённую форму теоремы Коши--Ковалевской (с использованием леммы Розенблюма вместо мажорант Коши [145]) и т.~д. Внутренний смысл теоремы 4 в том, что вторая группа коэффициентов при наличии ОП выражается через первую. Рассмотрено также аналогичное (\ref{42}) представление,
 но с интегрированием по промежутку $(x,\infty)$.
Существенные результаты по построению ОП высокого порядка в уточнённых областях аналитичности коэффициентов были получены И.Г.~Хачатряном.

Теперь перейдём к более подробному рассмотрению результатов из замечательной монографии [9]. А их довольно много. В этой монографии введены и изучены обобщения аналитических функций---$L$ и ($L,M$)--аналитические функции, что позволило в том числе объединить многие ранее полученные результаты. Для  таких функций построена теория рядов Тэйлора, введённых ранее в одном  частном случае Дельсартом. Наиболее принципиальным результатом представляется получение глобального ОП для случая переменных аналитических коэффициентов в пространстве а
 алитических функций.

\Theorem{ 5} Пусть все коэффициенты дифференциальных операторов $P,Q$ являются аналитическими функциями. Тогда в некотором пространстве аналитических функций (см. точные формулировки в [9]) существует ОП, представимый в интегральном виде и сплетающий  $P,Q$:
\begin{eqnarray}
\label{45}
P=D^n + p_{n-1}(x) D^{n-1}+\cdots+p_1(x)D+p_0(x),\\
Q=D^n + q_{n-1}(x) D^{n-1}+\cdots+q_1(x)D+q_0(x), D=\frac{d}{dx}.\nonumber
\end{eqnarray}
\Endproc

Существенные результаты получены в [9] и для представления ОП (\ref{45}) в интегральном виде. В этом случае для ядра ОП вместо стандартного двумерного гиперболического уравнения  (\ref{tag9}) получается так называемое уравнение Бианки. В [9] построена полная теория разрешимости задачи Коши для уравнения Бианки. Построение проведено методом Римана с использованием так называемых характеристических пирамид. Формулы получаются достаточно сложные, на некотором этапе, например, в них входят семикратные суммы. Поэтому для общего случая (\ref{45}) они но
 ят качественный характер, просто указывая на возможность реализации ОП в некотором интегральном виде. А вот для случая $Q=D^n$ интегральные представления уже существенно упрощены и получены в явном виде. В частности показано, что для модельных случаев операторов Штурма--Лиувилля эти представления сводятся к известным в теории ОП.

Разобран также важный вопрос об операторах, коммутирующих с производной или её степенью в пространствах аналитических функций. Коммутирующие операторы---это частный случай ОП. Их важность для теории ОП заключается в следующем. Предположим, что мы имеем ОП со свойством $TA=B\,T$, а также некоторый оператор $K$, коммутирующий с $B$. Тогда можно построить новый ОП с тем же свойством:
\begin{eqnarray}
\label{46}
T_1=KT, KB=BK, T_1A=KTA=KBT=BKT=BT_1.
\end{eqnarray}
Или аналогично для операторов $K$, коммутирующих с $A$:
\begin{eqnarray}
\label{47}
T_2=TK, KA=AK, T_2A=TKA=TAK=BTK=BT_2.
\end{eqnarray}

Таким образом, задача об описании всех ОП эквивалентна задаче об описании всех коммутирующих операторов с любым из двух сплетаемых, так как между этими двумя семействами есть взаимно--однозначное соответствие. Можно сказать, что для описания \textit{всех} ОП, сплетающих $A$ и $B$, достаточно знать \textit{всего один любой} такой ОП $T_1$ и множество \textit{всех} операторов $K$, коммутирующих с $B$ (с $A$). В этом случае все ОП получаются по формуле $T=KT_1$ ($T=T_1K$). С другой стороны, для описания \textit{всех} операторов, коммутирующих с $B$ (с $A$), достаточно знать  множес
 тво \textit{всех} ОП $T$ , сплетающих $A$ и $B$. В этом случае все коммутирующие с $B$ (с $A$) получаются по формуле $K=T_1T_2^{-1}$ ($K=T_1^{-1}T_2$) при дополнительном условии обратимости ОП. Таким образом, эти две задачи практически эквивалентны. Вот почему в теории ОП такую важную роль играют операторы дробного интегродифференцирования различных видов---при естественных условиях они коммутируют со 'своей родной' производной и служат материалом для конструирования ОП.

Считалось, что задача об описании операторов, коммутирующих с производной, давно решена в процитированных выше работах Лионса и Дельсарта, поэтому считалась решённой и задача об описании всех ОП для дифференциальных операторов с аналитическими коэффициентами.  Но как указано в [9], соответствующие работы содержат грубые ошибки, на самом деле там описано не всё множество коммутирующих операторов, а только его некоторая собственная часть. Правильный результат также получен М.~К.~Фаге и приведён в [9], с его помощью можно описать все возмо
 ные ОП, сплетающие данный дифференциальный оператор с производной того же порядка.

Отметим, что теория операторов в аналитических пространствах функций,  важная и для приложений к ОП, рассматривалась Ю.~Ф.~Коробейником в его  монографиях [146--147].
Теория интегральных операторов в ещё более общих пространствах рассмотрена, например,  в книгах В.~Г.~Фетисова [148], А.~Г.~Кусраева [149], А.~Э.~Пасенчука [150].

Укажем на существенность для перечисленных результатов равенства единице старших коэффициентов дифференциальных операторов. Так, в [151] доказано, что операторы $D$ и $\exp(-z)D$ не эквивалентны в пространстве целых аналитических функций. Многочисленные примеры эквивалентности или её отсутствия для дифференциальных операторов  в различных аналитических пространствах рассмотрены во второй части монографии [9]. Отметим, что существенную роль в теории общих ОП сыграли  также результаты В.~А.~Марченко [152--156].

О некоторых других результатах этого направления  см.  [157--162].

Таким образом, мне представляется более справедливым называть ОП со свойством (\ref{41}) ОП Дельсарта--Фаге, (или Дельсарта--Марченко--Фаге), а не как принято на Западе ОП Лионса--Дельсарта. О вкладе отечественных математиков не скажешь лучше ведущего американского специалиста по ОП Роберта Кэррола: "Идея ОП, возникшая в начале 50--х,  восходит к Гельфанду, Левитану, Марченко, Наймарку и другим. Она была подхвачена вновь Дельсартом и Лионсом..." [163].
\newpage

\section{5. Операторы преобразования типа Сонина и Пуассона.}

Перейдём к рассмотрению, наверное, самого известного класса ОП, сплетающих дифференциальный оператор Бесселя со второй производной:
\begin{equation}
\label{51}
T(B_\nu)f=(D^2)Tf, B_{\nu}=D^2+\frac{2\nu +1}{x}D, D^2=\frac{d^2}{dx^2}, \nu \in \mathbb{C}.
\end{equation}

Как было отмечено в п.~2, одним из способов построения ОП является установление соответствий между решениями соответствующих дифференциальных уравнений. Решениями уравнения вида $B_\nu f=\lambda f$ являются функции Бесселя, а уравнения $D^2f=\lambda f$---тригонометрические функции или экспонента. Поэтому прообразами ОП вида (\ref{51}) были формулы Пуассона и Сонина:
\begin{eqnarray}
\label{52}
J_{\nu} (x)=\frac{1}{\sqrt{\pi}\Gamma({\nu+\frac{1}{2}})2^{\nu-1}x^\nu}
\int_0^x \left( x^2-t^2\right)^{\nu-\frac{1}{2}}\cos(t)\,dt,\Re \nu> \frac{1}{2}\\
\label{53}
J_{\nu} (x)=\frac{2^{\nu+1}x^{\nu}}{\sqrt{\pi}\Gamma({\frac{1}{2}-\nu})}
\int_x^\infty \left( t^2-x^2\right)^{-\nu-\frac{1}{2}}\sin(t)\,dt,-\frac{1}{2}<\Re \nu<\frac{1}{2}
\end{eqnarray}

Как обычно, именные названия формул носят условный характер. Начнём с того, что функции Бесселя для произвольного индекса были введены и исследованы великим Леонардом Эйлером, когда Бессель ещё и не родился. Вообще эти функции возникли и постепенно вводились в период 1690--1770 гг. для решения дифференциального уравнения, которое мы сейчас называем именем Риккати. Участвовали в этом члены семьи Бернулли, венецианский граф Джаккопо Франческо Риккати и Эйлер [164--166]. Вместе с тем вклад самого Фридриха Бесселя в изучение названных его именем
 ункций достаточно велик.  Интеграл (\ref{52})  начал изучать Эйлер в 1769 г. Затем Парсеваль посчитал интеграл при $\nu=0$ в 1805 г.,  для целых $\nu$ формулу (\ref{52}) получил Плана в 1821 г., Пуассон вывел её для полуцелых $\nu$ в 1823 г., его метод применим и для целых $\nu$, но он этого не заметил.  Далее этот интеграл встречался в работах Куммера, Лоббато и Дюамеля. Окончательно формулу (\ref{52}), которую мы приписываем Пуассону, установил в общем случае Ломмель в 1868 г., а Сонин вывел формулу (\ref{53}) в 1880 г. Формулу Сонина, в существование которой все верят и считают прообр
 зом ОП (\ref{55}) , и которая бы выражала экспоненту или тригонометрические функции через интеграл по \textit{конечному промежутку} от функции Бесселя, я в литературе не нашёл.

\textsc{Определение 3.} ОП Пуассона называется выражение
\begin{equation}
\label{54}
P_{\nu}f=\frac{1}{\Gamma(\nu+1)2^{\nu}x^{2\nu}}
\int_0^x \left( x^2-t^2\right)^{\nu-\frac{1}{2}}f(t)\,dt,\Re \nu> -\frac{1}{2}.
\end{equation}
ОП Сонина называется выражение
\begin{equation}
\label{55}
S_{\nu}f=\frac{2^{\nu+\frac{1}{2}}}{\Gamma(\frac{1}{2}-\nu)}\frac{d}{dx}
\int_0^x \left( x^2-t^2\right)^{-\nu-\frac{1}{2}}t^{2\nu+1}f(t)\,dt,\Re \nu< \frac{1}{2}.
\end{equation}

Операторы (\ref{54})--(\ref{55}) действуют как ОП по формулам
\begin{equation}
\label{56}
S_\nu B_\nu=D^2 S_\nu, P_\nu D^2=B_\nu P_\nu.
\end{equation}
Их можно доопределить на все значения $\nu\in\mathbb {C}$.

Идею изучения операторов подобных (\ref{54})--(\ref{55}) высказывал ещё Лиувилль, их реальное использование в контексте теории функций Бесселя начал Николай Яковлевич Сонин [167]. Как ОП эти операторы были введены в работах Дельсарта [122, 168--170] и затем изучены в работах Дельсарта и Лионса [123--126, 17]. Поэтому мы будем называть (\ref{54})--(\ref{55}) ОП Сонина--Пуассона--Дельсарта (СПД). В нашей стране об операторах СПД в основном узнали из великолепно написанной статьи Б.~М.~Левитана [171].

Не будет преувеличением сказать, что  операторы СПД (\ref{54})--(\ref{55}) являются самыми знаменитыми объектами всей теории ОП, их изучению, приложениям и обобщениям посвящены сотни работ. Перечислим очень кратко только основные направления.

Дельсартом на базе ОП СПД было введено фундаментальное понятие обобщённого сдвига.

\textsc{Определение 4.} Оператором обобщённого сдвига (ООС) называется решение $u(x,y)=T_x^yf(x)$ задачи
\begin{eqnarray}
\label{57}
(B_\nu)_y u(x,y)=(\frac{\partial^2}{\partial y^2}+\frac{2\nu +1}{y}\frac{\partial}{\partial y})\  u(x,y)= \frac{\partial^2}{\partial x^2} \ u(x,y),\\ \nonumber
u(x,0)=f(x), u_y (x,0)=0.
\end{eqnarray}

Название объясняется тем, что в частном случае $\nu=-1/2$ ООС сводится к почти обычному сдвигу $T_x^yf(x)=\frac{1}{2}\left( f(x+y)+f(x-y)\right)$. Для ООС (\ref{57}) Дельсартом была получена явная формула
\begin{equation}
\label{58}
T_x^y f(x)=\frac{\Gamma(\nu+1)}{\sqrt{\pi}\Gamma(\nu+\frac{1}{2})} \int_0^{\pi}f(\sqrt{x^2+y^2-2xy \cos(t)}) \sin^{2\nu}t\,dt.
\end{equation}
Можно рассматривать в определении (\ref{57}) и произвольные пары дифференциальных (или даже любых) операторов. Например, при таком определении $\frac{\partial u}{\partial x}=\frac{\partial u}{\partial y},u(x,0)=f(x)$ получаем привычный сдвиг $T_x^y f(x)=f(x+y)$. Отметим, что ООС (\ref{57})--(\ref{58}) явно выражаются через ОП СПД (\ref{54})--(\ref{55}) [4--6, 22--29, 171].

Теперь, когда есть новый сдвиг, то можно обобщать прежние теории, основанные на обычном сдвиге. Так возникли теория обобщённых почти--периодических функций [169,  25--27, 18, 20, 126](отметим, что определение и основные свойства почти--периодических функций были впервые даны профессором Юрьевского ( Тартусского (Дерптского)) университета Пирсом Георгиевичем Болем задолго до Г.~Бора, Боль также получил знаменитую топологическую теорему о неподвижных точках непрерывного отображения сферы в себя до Брауэра; из эвакуированных за время войны и револ
 юции сотрудников Юрьевского университета начал в 1918 г. создаваться Воронежский университет), разложения по обобщённым рядам Тэйлора, справедливо названных рядами Тэйлора--Дельсарта [168, 9, 24--26], обобщённая свёртка и её приложения [172--174].

Огромную роль конструкции ОП и ООС сыграли в теории уравнений с частными производными. ОП позволяют преобразовывать более сложные уравнения в более простые, ООС помогают в сингулярных уравнениях переносить особенность из начала в произвольную точку, а также с помощью обобщённой свёртки находить фундаментальные решения.

Особо отметим один класс уравнений с частными производными с особенностями, типичным представителем которого является $B$--эллиптическое уравнение с операторами Бесселя по каждой переменной вида
\begin{equation}
\label{59}
\sum_{k=1}^{n}B_{\nu,x_k}u(x_1,\dots, x_n)=f,
\end{equation}
аналогично рассматриваются $B$--гиперболические и $B$--параболические уравнения. Изучение этого класса уравнений было начато в работах Эйлера, Пуассона, Дарбу, продолжено в теории обобщённого осесимметрического потенциала А.~Вайнштейна и в трудах отечественных математиков И.~Е.~Егорова, Я.~И.~Житомирского, Л.~Д.~Кудрявцева, П.~И.~Лизоркина, М.~И.~Матийчука, Л.~Г.~Михайлова, М.~Н.~Олевского, М.~М.~Смирнова, С.~А.~Терсенова,   Хе Кан Чера, А.~И.~Янушаускаса   и других.

Наиболее полно весь круг вопросов для  уравнений с операторами Бесселя был изучен воронежским математиком Иваном Александровичем Киприяновым и его учениками Л.~А.~Ивановым, А.~В.~Рыжковым, В.~В.~Катраховым, В.~П.~Архиповым, А.~Н.~Байдаковым, Б.~М.~Богачёвым, А.~Л.~Бродским, Г.~А.~Виноградовой, В.~А.~Зайцевым, Ю.~В.~Засориным, Г.~М.~Каганом,  А.~А.~Катраховой, Н.~И.~Киприяновой, В.~И.~Кононенко, М.~И.~Ключанцевым, А.~А.~Куликовым, А.~А.~Лариным, М.~А.~Лейзиным, Л.~Н.~Ляховым, А.~Б.~Муравником, И.~П.~Половинкиным, А.~Ю.~Сазоновым, С.~М.~Ситником, В.~П.~Шацким, В.~Я.~Ярославце
 вой;   основные результаты этого направления представлены в [30]. Для описания классов решений соответствующих уравнений И.~А.~Киприяновым были введены и изучены функциональные пространства [175],  позднее названные его именем (см. монографии Х.~Трибеля, Л.~Д.~Кудрявцева и С.~М.~Никольского [176--177], в которых пространствам Киприянова посвящены отдельные параграфы).

В  этом направлении работал В.~В.~Катрахов (Валерий Вячеславович Катрахов скончался в 2010 году, вечная ему память). Сейчас уравнения с оператором Бесселя и связанные с ними вопросы изучают  А.~В.~Глушак, В.~С.~Гулиев,  Л.~Н.~Ляхов   со своими коллегами и учениками. Задачи для операторных уравнений вида (\ref{59}), берущие начало в известной монографии [11], рассматривали А.~В.~Глушак, С.~Б.~Шмулевич, В.~Д.~Репников и другие.

Здесь следует отметить, что первой фундаментальной работой, с которой начался отсчёт  изучения вырождающихся и сингулярных дифференциальных уравнений в частных производных с переменными коэффициентами, является статья М.~В.~Келдыша [178] (задача Е). Общая теория вырождающихся уравнений эллиптического типа впоследствии разрабатывалась Г.~Фикерой, О.~А.~Олейник, Е.~А.~Радкевичем, В.~Н.~Враговым, В.~П.~Глушко и многими другими.

Интересные результаты по изучению ОП СПД были получены В.~В.~Катраховым, они также изложены в специально переработанном Р.~Кэрролом виде в отдельной главе в [30].
Рассматривались также ОП для многочисленных обобщений оператора Бесселя [179--190]. Вместе с тем были построены аналогичные теории и для некоторых других модельных операторов, например таких [4--6, 189--190]:
\begin{eqnarray}
\label{510}
A=\frac{1}{v(x)}\frac{d}{dx}v(x)\frac{d}{dx},v(x)=\sin^{2\nu+1}x,\sh^{2\nu+1}x,\\ (e^x-e^{-x})^{2\nu+1}(e^x+e^{-x})^{2\mu+1}.\nonumber
\end{eqnarray}
Важность операторов $A$ вида (\ref{510}) для теории заключается в том, что по знаменитой формуле Гельфанда они представляют радиальную часть оператора Лапласа на симметрических пространствах [191]. При этом оператор Бесселя получается при выборе в (\ref{510}) $v(x)=x^{2\nu+1}$. Подход через ОП позволяет с единых позиций изучать многие свойства специальных функций, через которые выражаются собственные функции (\ref{510}), этот подход последовательно проводится в [4--6].

Другим модельным оператором, для которого построены ОП, является оператор Эйри $D^2+x$, см. [387]. Рассматривался также его возмущённый вариант, связанный с эффектом Штарка из квантовой механики.

Важным обобщением операторов Сонина--Пуассона--Дельсарта являются ОП для гипербесселевых функций, которые были подробно изучены в работах Ивана Димовски и его учеников. Соответствующие ОП заслуженно  получили в литературе названия ОП Сонина--Димовски и Пуассона--Димовски, они также изучались в работах ученицы Ивана Димовски --- Виржинии Киряковой.

Другой важный класс обобщений операторов СПД возник при изучении операторов Дункла (я благодарен профессору В.П.~Лексину, который научил меня правильно произносить фамилию этого австрийского математика). Операторы Дункла включаются в общую схему ОП как специальный случай, когда один из сплетаемых операторов является дифференциально--разностным, а второй по--прежнему дифференциальным. Построению ОП типа Сонина и Пуассона для операторов Дункла в последнее время посвящено много работ, см. [].

Отметим, что некоторый способ построения ОП со свойствами (\ref{76}) на решениях соответствующих уравнений для частного случая целых $\nu$ упоминается в литературе под названием метода Крама--Крейна [21--23], но он приводит к длинной последовательности усложняющихся подстановок, а результирующий ОП не может быть выписан в явном виде. На самом деле, и подстановки Крама--Крейна могут быть включены в общую схему ОП. Они возникают, когда сплетаемые операторы остаются дифференциальными, а ОП для этой пары ищется уже не в интегральном, а в дифференциа
 льном представлении. В наиболее общем виде эта теория известна как метод преобразований Дарбу [Матвеев], к которым и относятся в частном случае  подстановки Крама--Крейна. К таким дифференциальным ОП формально могут быть отнесены и многие известные замены переменных типа преобразования Миуры, которые сплетают нелинейный и линейный операторы, решая задачу линеаризации.
\newpage

\section{6. Связь операторов преобразования с дробным интегродифференцированием.}

Операторы дробного интегродифференцирования играют важную роль во многих современных разделах математики: уравнениях с частными производными, функциональном анализе и теории функций, специальных функциях, многочисленных приложениях.
Число книг по дробному исчислению всё время возрастает, процесс их увеличения  недавно даже был изображен графически на красивом цветном плакате, опубликованном в специализированном журнале Fractional Calculus and Applied Analysis; этот международный журнал издаётся Болгарской АН, его главным редактором является Виржиния Кирякова. В виду вышесказанного мы ограничимся из всей литературы лишь упоминанием энциклопедической монографии [34] ( её в России иногда называют Красной Книгой по цвету обложки),  диссертации [353], которую можно рассматривать как со
 ременное продолжение [34], а также известных книг
[192--196, 355].

Для теории специальных функций важность дробного интегродифференцирования отражена в названии известной статьи Виржинии Киряковой [197]: "\textbf{\textit{Все специальные функции получаются дробным интегродифференцированием элементарных функций}}"!\ (кроме $H$---функции Фокса --- замечание профессора А.А.~Килбаса (Анатолий Александрович Килбас погиб в 2010 году, вечная ему память)).

Приведём список основных операторов дробного интегродифференцирования: Римана--Лиувилля, Эрдейи--Кобера, дробного интеграла по произвольной функции $g(x)$
\begin{eqnarray}
\label{61}
I_{0+,x}^{\alpha}f=\frac{1}{\Gamma(\alpha)}\int_0^x \left( x-t\right)^{\alpha-1}f(t)d\,t,\\ \nonumber
I_{-,x}^{\alpha}f=\frac{1}{\Gamma(\alpha)}\int_x^\infty \left( t-x\right)^{\alpha-1}f(t)d\,t,
\end{eqnarray}
\begin{eqnarray}
\label{62}
I_{0+,x^2}^{\alpha}f=\frac{1}{\Gamma(\alpha)}\int_0^x \left( x^2-t^2\right)^{\alpha-1}2tf(t)d\,t,\\ \nonumber
I_{-,x^2}^{\alpha}f=\frac{1}{\Gamma(\alpha)}\int_x^\infty \left( t^2-x^2\right)^{\alpha-1}2tf(t)d\,t,
\end{eqnarray}
\begin{eqnarray}
\label{63}
I_{0+,g}^{\alpha}f=\frac{1}{\Gamma(\alpha)}\int_0^x \left( g(x)-g(t)\right)^{\alpha-1}g'(t)f(t)d\,t,\\ \nonumber
I_{-,g}^{\alpha}f=\frac{1}{\Gamma(\alpha)}\int_x^\infty \left( g(t)-g(x)\right)^{\alpha-1}g'(t)f(t)d\,t,
\end{eqnarray}
во всех случаях предполагается, что $\Re\alpha>0$, на оставшиеся значения
$\alpha$ формулы также без труда продолжаются [34].
При этом обычные дробные интегралы получаются при выборе в (\ref{63}) $g(x)=x$, Эрдейи--Кобера при $g(x)=x^2$, Адамара при $g(x)=\ln x$.

Связь с ОП проявляется в том, что, как видно из (\ref{54})--(\ref{55}), ОП СПД с точностью до множителей как раз и являются операторами Эрдейи--Кобера, то есть дробными степенями $(\frac{d}{dx^2})^{-\alpha}$. Так же как ОП для перечисленных в (\ref{510}) модельных операторов являются дробными степенями $(\frac{d}{dg})^{-\alpha}$ при $g(x)=\sin x,\sh x$. Поэтому основные свойства этих ОП можно получить из теории операторов дробного интегродифференцирования, а не изобретать заново, что нередко и делалось.

Операторы дробного интегрирования в комплексной плоскости, к которым сводится понятие дробного интегрирования по произвольной функции (\ref{63}), были введены в работах М.М.~Джрбашяна. На это указал мне Армен Джрбашян --- сын Мхитара Мкртичевича Джрбашяна. Соответствующие операторы, определяемые в виде ряда, первоначально назывались операторами Гельфонда--Леонтьева (по предложению В.Ф.~Коробейника, который был рецензентом оригинальной статьи --- личное сообщение), а сейчас чаще называются операторами дробного интегродифференцирования Дж
 рбашяна--Гельфонда--Леонтьева [353].

Наиболее широкие обобщения операторов дробного интегродифференцирования были предложены Виржинией Киряковой [353] в виде
$$
x^{-\alpha}\int_0^x t^{\alpha-1} K\left(\frac{x}{t}\right)f(t)\,dt\ , x^{-\alpha}\int_x^\infty t^{\alpha-1} K\left(\frac{x}{t}\right)f(t)\,dt,
$$
где $K\left(\frac{x}{t}\right)$ --- подходящее ядро, $f(t)$ --- достаточно произвольная функция.
Такие объекты сами по себе не являются новыми, это одна из форм записи интегрального оператора с однородным ядром степени минус единица
$$
\int_0^\infty K(x,t) f(t) \,dt\ , K(\lambda x,\lambda y)= \frac{1}{\lambda} K( x,y),
$$
или свёртки Меллина
$$
\int_0^\infty K\left(\frac{x}{t}\right) f(t) \,\frac{dt}{t}
$$
при дополнительном условии, что ядро интегрального оператора или свёртки обращается в ноль при $t<x$ или при $t>x$.
Для этой формы обобщённого интегродифференцирования в [353] построена достаточно полная теория.

Теперь вспомним связь ОП с коммутированием (\ref{46})--(\ref{47}). Получаем, что домножение с нужной стороны на различные формы дробных интегралов приводит к новым ОП. Для ОП Векуа--Эрдейи--Лаундеса эта идея использована в [99--101] для построения новых преобразований с различными экзотическими ядрами, для ОП СПД см. ниже. Новая возможность---это использование операторов дробного интегродифференцирования, коммутирующих с оператором Бесселя. Для этого требуется, например, явное описание дробных степеней оператора Бесселя. И оно есть, хоть об этом ма
 о кто знает!

В этом пункте построены в явном интегральном виде  дробные степени дифференциального оператора Бесселя, заданного на подходящих гладких функциях,  который мы для удобства переопределим в виде:
$$B_{\nu}=D^2+\frac{\nu}{x}D, \ Re~\nu\ge 0.$$

Безусловно, по аналогии с обычными производными, можно определить такие дробные степени на полуоси при помощи естественного свойства, что они действуют как умножение на степень в образах преобразования Ханкеля. Такой подход оправдан за неимением лучшего и позволяет получить ряд интересных результатов, хотя на этом пути невозможно, по--видимому, получить явные представления дробных степеней. Но представим на минуту, что мы умеем определять обычные операторы дробного интегрирования Римана--Лиувилля только через их действие в образах
  преобразований Лапласа или Меллина, а интегральные формулы для этих операторов нам неизвестны. Тогда сразу становится ясным, что подобная теория будет достаточно бедной и лишится большинства своих наиболее полезных и красивых результатов. Примерно в таком состоянии сейчас находится теория дробных степеней оператора Бесселя, поэтому их построение в явном интегральном виде является актуальной и интересной задачей.

Несмотря на то, что операторы дробного интегродифференцирования Бесселя являются  пока менее известными  по сравнению с перечисленными выше, уже созданы основы их теории. Формальные определения даны в [198--201].

\textsc{Определение 5.} Пусть $f(x)\in C^{2k}(0,b]$. Определим правосторонний оператор дробного интегрирования Бесселя при условии $f^{(i)}(b)=0, 0\leq i \leq 2k-1, k \in N$  по формуле
\begin{eqnarray}
\label{651}
(B_{b-}^{{\nu},k}f)(x)=\frac{1}{{\Gamma (2k)}}\int_{x}^{b}
\left(\frac{y^{2}-x^{2}}{2y}\right)^{2k-1} {_2}F_1(k+\frac{{\nu}-1}{2},k;2k;1-\frac{x^{2}}{y^{2}})\cdot
\nonumber\\
\cdot f(y)\,dy=
\frac{\sqrt{{\pi}}}{2^{2k-1}{\Gamma (k)}}
\int_{x}^{b}(y^{2}-x^{2})^{k-\frac{1}{2}}
\left(\frac{y}{x}\right)^{\frac{{\nu}}{2}}
P_{\frac{{\nu}}{2}-1}^{\frac{1}{2}-k}
\left(\frac{1}{2}\left(\frac{x}{y}+\frac{y}{x}\right)\right)\cdot\\
\cdot f(y)\,dy.\phantom{111111}\nonumber
\end{eqnarray}

Определим левосторонний оператор дробного интегрирования Бесселя при условии $f^{(i)}(a)=0, 0\leq i \leq 2k-1, k \in N$  по формуле

\begin{eqnarray}
\label{65}
(B_{a+}^{{\nu},k}f)(x)=\frac{1}{{\Gamma (2k)}}\int_{a}^{x}
\Bigl(\frac{x^{2}-y^{2}}{2x}\Bigr)^{2k-1} {_2}F_1(k+\frac{{\nu}-1}{2},k;2k; ;1-\frac{y^{2}}{x^{2}})\cdot
\nonumber \\
\cdot f(y)\,dy=\frac{\sqrt{{\pi}}}{2^{2k-1}{\Gamma (k)}}
\int_{a}^{x}(x^{2}-y^{2})^{\left(k-\frac{1}{2}\right)}\left(\frac{x}{y}\right)^{\frac{{\nu}}{2}}
P_{\frac{{\nu}}{2}-1}^{\frac{1}{2}-k}\left(\frac{1}{2}\left(\frac{x}{y}+\frac{y}{x}\right)\right)\cdot\\
\cdot f(y)\,dy,\phantom{111111}\nonumber
\end{eqnarray}
где ${_2}F_{1}$ - гипергеометрическая функция Гаусса, $P_{\nu}^{\mu}(z)$ - функция Лежандра.

Введённые операторы и являются интегральными реализациями отрицательных целых степеней оператора Бесселя $(B_{\nu})^{-k}$. Их распространение на произвольные комплексные значения параметра $k$ проводится аналогично классическому случаю. При ${\nu}=0$ оператор Бесселя сводится ко второй производной, а введённые операторы  --- к дробным интегралам Римана-Лиувилля
\begin{equation*}
B_{b-}^{0,k}f=I_{b-}^{2k}f,\ \ B_{a+}^{0,k}f=I_{a+}^{2k}f.
\end{equation*}

Определённые выше дробные степени оператора Бесселя (\ref{651}--\ref{65}) в явном виде как одна из разновидностей операторов дробного интегродифференцирования были определены автором в [198--201]. При этом в качестве наводящих соображений были существенно использованы результаты работы [204], в которой рассматривалось  решение  в явном виде обыкновенных дифференциальных уравнений с целыми степенями операторов Бесселя.
При этом автором было замечено, что выражение гипергеометрических функций Гаусса в формулах (\ref{651}--\ref{65}) через функции Лежандра существенно упрощает вычисления.

При помощи достаточно простых выкладок, основанных на тождествах для функций Лежандра, получается первоначальный набор простейших свойств операторов (\ref{651}--\ref{65}). К ним относятся полугрупповое свойство по параметру $k$, действие как умножение на соответствующую степень в образах преобразования Ханкеля (при более ограничительных условиях на параметры) и преобразования Мейера (при менее ограничительных условиях на параметры).

Многочисленные приложения операторов дробного интегродифференцирования Римана--Лиувилля основаны на их вхождении в остаточный член формулы Тэйлора. Поэтому после определения дробных степеней оператора Бесселя сразу возникает задача о построении обобщённой формулы Тэйлора, в которой функция раскладывается по степеням оператора Бесселя. Эта задача возникла достаточно давно и имеет некоторую историю.

Впервые формулы  разложения по степеням оператора Бесселя были получены Жаном Дельсартом  (ряды Тэйлора--Дельсарта). Общий способ их построения изложен в [9] в терминах операторно--аналитических функций. Но ряды Тэйлора--Дельсарта позволяют разложить по степеням оператора Бесселя не обычный, а обобщённый сдвиг. По существу такие разложения являются операторными вариантами рядов для функции Бесселя, так же как обычные ряды Тэйлора являются операторными версиями разложения в ряд экспоненты. Разумеется, ряды Тэйлора--Дельсарта имеют сво
 ю область приложений. Но для численного решения дифференциальных уравнений в частных производных нужны обобщённые формулы и ряды Тэйлора
несколько другой природы. При пересчёте решения со слоя на слой, например, методом сеток формулы для обобщённого сдвига совершенно бесполезны, а нужны именно формулы для обычного сдвига, выражающие решение на очередном рассчитываемом слое через его значения на предыдущих слоях. Оказалось, что строить такие формулы для обычного сдвига намного труднее, чем для обобщённого, так как они уже не являются прямыми аналогами известных тождеств для специальных функций.

Впервые с указанной мотивацией для применения к численному решению уравнений с оператором Бесселя методом конечных элементов  формула Тэйлора нужного типа была рассмотрена в работе В.В.~Катрахова [203]. Но полученный там результат может рассматриваться только как первое приближение для желаемых формул в явном виде, так как коэффициенты выражались неопределёнными постоянными, задаваемыми системой рекуррентных соотношений, а ядро остаточного члена представлялось некоторым многократным интегралом. Это не случайно, угадать одновреме
 но явный вид ядер и остаточных членов невозможно, пока не известны конкретные выражения для остатка в виде дробных степеней оператора Бесселя.

Мы находимся в другой ситуации, когда нужные выражения для дробных степеней известны в явном виде (\ref{651}--\ref{65}), поэтому перейдём к описанию  полного решение задачи о явном построении обобщённой формулы Тэйлора для разложения обычного сдвига в ряд по степеням оператора Бесселя. Подобная формула в простейшем случае была анонсирована ранее в [198--201].

\Theorem { 7} Справедлива формула Тэйлора разложения произвольной достаточно гладкой функции по степеням дифференциального оператора Бесселя при $x=b$ с остаточным членом в интегральной форме
\begin{eqnarray}
\label{66}
f(x)=\sum_{i=1}^{k}
\Biggl\{
 \frac{1}{\Gamma(2i-1)}
\Bigl(\frac{b^{2}-x^{2}}{2b}\Bigr)^{2i-2}
{_2}F_{1}(i+\frac{{\nu}-1}{2},i-1;2i-1; \nonumber\\
1-\frac{x^{2}}{b^{2}})(B^{i-1}f)|_{b}- \frac{1}{\Gamma(2i)}
\Bigl(\frac{b^{2}-x^{2}}{2b}\Bigr)^{2i-1} {_2}F_{1}(i+\frac{{\nu}-1}{2},i;2i;1-\frac{x^{2}}{b^{2}})\cdot \\ \nonumber\cdot(DB^{i-1}f)|_{b}
\Biggr\} +B_{b-}^{{\nu},k}(B^{k}f),\phantom{11111111111111111}
\end{eqnarray}
где $B_{b-}^{{\nu},k}$ есть оператор левостороннего дробного интегрирования Бесселя (\ref{651}), ${_2}F_{1}$ --- гипергеометрическая функция Гаусса.
\Endproc

Также справедлива двойственная формула , использующая формально сопряжённый оператор
\begin{equation*}
C_{{\nu}}f=D^{2}f-D(\frac{{\nu}}{y}f)=D^{2}-\frac{{\nu}}{y}Df+\frac{{\nu}}{y^{2}}
\end{equation*}

\Theorem { 8}  Справедлива формула Тэйлора разложения произвольной достаточно гладкой функции по степеням дифференциального оператора  $C_{{\nu}}$  при $x=a$ с остаточным членом в интегральной форме
\begin{eqnarray}
\label{67}
f(x)=\sum_{i=1}^{k}
\Biggl\{\frac{1}{\Gamma(2i-1)}
\Bigl(\frac{x^{2}-a^{2}}{2x}\Bigr)^{2i-2} (\frac{a}{x})\  {_2}F_{1}
(i+\frac{{\nu}-1}{2},i;2i-1;\nonumber\\
1-\frac{a^{2}}{x^{2}})
(C_{{\nu}}^{i-1}f)|_{a}
+\frac{1}{\Gamma(2i)}
\Bigl(\frac{x^{2}-a^{2}}{2x}\Bigr)^{2i-1}\cdot \phantom{1111111}\\\nonumber \cdot  \ {_2}F_{1}(i+\frac{{\nu}-1}{2},i;2i; 1-\frac{a^{2}}{x^{2}}) \ a^{{\nu}} \ (Dx^{-{\nu}}C_{{\nu}}^{i-1}f)|_{a}
\Biggr\} +B_{a+}^{{\nu},k}(C_{{\nu}}^{k}f),
\end{eqnarray}
где $B_{a+}^{{\nu},k}$ есть оператор правостороннего дробного интегрирования Бесселя (\ref{65}), $_{2}F_{1}$ --- гипергеометрическая функция Гаусса.

\Endproc

Гипергеометрические функции в формулах Тэйлора могут быть выражены через функции Лежандра аналогично определениям. Из формул (\ref{66})--(\ref{67}) сразу следуют модификации операторов дробного интегродифференцирования Бесселя  по Герасимову--Капуто. Для этого нужно в определениях (\ref{651}--\ref{65}) вычесть из функций начальный отрезок их разложения по соответствующей обобщённой формуле Тэйлора.

Автором изучены и другие основные свойства операторов дробного интегродифференцирования Бесселя вида (\ref{651}--\ref{65}): действие в стандартных функциональных пространствах, интегральные уравнения в особых случаях
$ a=0 $ или $b=\infty$ на полуоси, модификации операторов по Маршо и Капуто, сразу вытекающие из формул Тэйлора (\ref{66})--(\ref{67}),  связь с преобразованиями Фурье--Бесселя, Лапласа, Меллина и Мейера,  действие на распределениях  в классах типа Лизоркина, связь с обобщённым преобразованием Стильтьеса, весовые и разносторонние формулы композиций, выражения левосторонних операторов через правосторонние и наоборот с помощью сингулярных интегралов, а также приложения в теории операторов преобразования.

Следует отметить, что введённые операторы дробного интегродифференцирования Бесселя родственны некоторым известным классам интегральных преобразований. К ним относятся  операторы Лава, подробно изученные О.~А.~Репиным операторы Сайго [205], гипербесселевы операторы И.~Димовски, операторы Бушмана--Эрдейи (о них см. далее). Кроме того, введённые операторы (\ref{66})--(\ref{67}) можно включить в общую канву интегральных операторов с $G$--функциями Мейера, функциями Райта   и $H$--функциями Фокса в ядрах.

Задачу о выводе формул Тэйлора (\ref{66})--(\ref{67}) по предложенному плану с использованием интегрирования по частям я поставил в 1994 г. для курсовой работы одной студентке третьего курса Дальневосточного политехнического института, в котором тогда преподавал во Владивостоке (при этом двадцать сотрудников двух академических институтов учили в течение всех пяти лет одну группу из двадцати студентов, читая все курсы). Она спросила, что такое гипергеометрическая функция, записала ссылку на книгу Бэйтмена--Эрдейи и ушла. Через две недели без еди
 ной консультации студентка принесла текст, в котором без единой помарки был записан подробный вывод теорем 7,8 и их обобщений (см ниже). Ничего подобного в своей жизни я больше не видел, а жаль, завидую коллегам, которые сталкивались с подобными случаями чаще. Сейчас, насколько я знаю, Дина Коновалова живёт и работает в Новосибирске.

Рассмотрены и более общие комбинированные дробные степени для пары операторов $(\frac{1}{x}D)^{m}(B_{{\nu}})^{k}$.
Это семейство операторов интересно тем, что содержит обычные операторы Римана --Лиувилля ($m=0,\ {\nu}=0$), дробное интегродифференцирование Бесселя (m=0), операторы Эрдейи--Кобера (k=0). В этом случае также доказаны соответствующие формулы Тэйлора.

Далее приведём выражение для резольвенты дробных степеней оператора Бесселя. Оно обобщает знаменитую формулу для дробных интегралов Римана--Лиувилля, написанную без доказательства Эйнаром Хилле и Яковом Давидовичем Тамаркиным в работе   1930 года [206]. В этой работе указывалось, что формула для резольвенты может быть выведена методом преобразования Лапласа с использованием нового на тот момент понятия свёртки в духе работ Дёйтча, но этот способ похоже не был никогда реализован. Формула Тамаркина--Хилле была на самом деле впервые доказ
 на в монографии М.М. Джрбашяна [207] обычным для теории интегральных уравнений методом последовательных приближений, хотя в монографии Мхитара  Мкртитевича нет упоминания, что доказательство даётся им впервые, что характеризует этого замечательного математика. Поэтому, возможно, исторически правильным было бы называть формулу для резольвенты операторов дробного интегрирования Римана--Лиувилля формулой \textsl{Тамаркина--Хилле--Джрбашяна}. Кроме того,  в [207] впервые были подробно изучены свойства функции Миттаг--Лефлёра, из этой книги от
 чественные математики узнали о существовании подобной функции. В дальнейшем свойства функций Миттаг--Лефлёра изучались в работах многих математиков, результаты дальнейших  исследований самого М.М. Джрбашяна просуммированы им в [352].

Формула  Тамаркина--Хилле--Джрбашяна  является самым известным применением функций Миттаг--Лефлёра, а также  основой огромного числа работ по приложениям дробного исчисления к дифференциальным уравнениям. Как отмечает в своём недавнем обзоре Виржиния Кирякова [354],  итальянские математики R. Gorenflo и F. Mainardi предложили называть функцию Миттаг--Лефлёра королевской функцией теории дробного исчисления (Queen function of the fractional calculus) [358]. Аналогично, самым известным применением $H$--функции Фокса по моему мнению является результат А.Н.~Кочубея [208--209], в
  работах которого через указанную $H$--функцию Фокса, не сводящуюся к функциям из более простых классов, была выражена функция Грина для уравнения дробной диффузии.

\Theorem { 9}  Для резольвенты оператора дробного интегрирования Бесселя (\ref{651}) при $a=0, \ 0\leq \nu <1$ на подходящих функциях справедливо интегральное представление
\begin{equation}
R_\lambda f=-\frac{1}{\lambda}f-\frac{1}{\lambda}\int_0^x K(x,y)f(y)\,dy,
\end{equation}
где ядро $K(x,y)$  выражается по формулам
\begin{eqnarray}
\label{900}
K(x,y)=\frac{2y}{x^2-y^2}\int_0^1 S_{\alpha,\nu}(z(t))
\frac{dt}
{
\left(
t\left(1-t\right)
\right)
^{\frac{\nu+1}{2}}
},\\ \nonumber
z(t)=\left(\frac{t(1-t)\left(x^2-y^2\right)^2}{\left(1-\left(1-\frac{x^2}{y^2}
\right)t\right)4y^2}\right)^\alpha,\\ \nonumber
S_{\alpha,\nu}(z)=\sum_{k=1}^{\infty}\frac{z^k}{\Gamma(\alpha k+\frac{\nu-1}{2})\Gamma(\alpha k-\frac{\nu-1}{2})}
\end{eqnarray}
---разновидность гипергеометрической функции Райта или Фокса.
\Endproc

Подобные  функции также встречались в работах А.~В.~Псху [196], это специальные случаи более общих функций Райта, $G$--функций Мейера и $H$--функций Фокса, которые первоначально вводились как обобщения функций Бесселя [102--120]. Отметим также такой бессмысленный термин, к сожалению использованный в [34], как функция Бесселя--Мэйтленда, который возник в результате недоразумения из расчленения имени Эдварда Мэйтленда Райта. Современное состояние теории обобщённых гипергеометрических функций типа Миттаг--Лефлёра, Мейера, Райта и Фокса изложено в [353
 --357].
Интересной задачей является дальнейшее упрощение представления ядра (\ref{900}), если оно возможно.

Основные интегральные представления для функции Миттаг--Лефлёра с учётом формулы Тамаркина--Хилле--Джрбашяна можно теперь безвозмездно, то есть даром \copyright, получить из известных формул Гильберта для резольвент, что, насколько известно автору, ранее не отмечалось.

Отметим, что резольвенты, как аналитические оператор--функции, коммутируют со своими операторами. Следовательно, их в принципе можно опять рассматривать как сырьё для получения новых ОП. Кроме того, из основного свойства ОП (\ref{1.1}) формально следует подобное соотношение и для резольвент $TR_\lambda(A)=R_\lambda(B)T$, но мне неизвестны применения этого соотношения.
Более того, опять же формально сразу получается и общая формула
\begin{equation}
\label{610}
Tf(A)=f(B)T
\end{equation}
для произвольного линейного ОП со свойством (\ref{1.1}) и аналитической функции $f$.

Полученная формула для резольвенты дробных степеней оператора Бесселя позволяет рассматривать  задачи для обыкновенного интегро--дифференциального уравнения  вида
\begin{equation}
\label{611}
B_\nu^\alpha u(x)-\lambda u(x)=f(x),
\end{equation}
при различных краевых условиях. По аналогии с известными результатами возможно также рассмотрение уравнений в частных производных с дробными степенями оператора Бесселя и их модификациями по Капуто.
К числу таких уравнений относятся обобщение $B$--эллиптического по терминологии И.А.Киприянова [30] дробного уравнения Лапласа--Бесселя
\begin{equation*}
\sum_{k=1}^n B_{\nu_k}^{\alpha_k} u(x_1,x_2,\ldots x_n)=f(x_1,x_2,\ldots x_n),
\end{equation*}
нестационарное уравнение вида
\begin{equation*}
B_{\nu,t}^\alpha u(x,t)=\Delta_x u(x,t)+f(x,t),
\end{equation*}
а также операторные уравнения в  пространствах Банаха, аналогичные изученным в [213].

Кроме того, рассмотрение  спектральных свойств уравнения (\ref{611}) нуждается в изучении асимтотики функции $K(x,y)$ из формулы (\ref{900}) в комплексной плоскости, а также распределения её корней.
\newpage

\section{7. Операторы преобразования Бушмана--Эрдейи.}

Рассмотрим важный класс ОП, который при определённом выборе параметров является одновременным обобщением ОП СПД и их сопряжённых, операторов дробного интегродифференцирования Римана--Лиувилля
и Эрдейи--Кобера, а также интегральных преобразований Мелера--Фока.

\textsc{Определение 6.} Операторами Бушмана--Эрдейи называются интегральные операторы
\begin{eqnarray}
\label{71}
B_{0+}^{\nu,\mu}f=\int_0^x \left( x^2-t^2\right)^{-\frac{\mu}{2}}P_\nu^\mu \left(\frac{x}{t}\right)f(t)d\,t,\\
E_{0+}^{\nu,\mu}f=\int_0^x \left( x^2-t^2\right)^{-\frac{\mu}{2}}\mathbb{P}_\nu^\mu \left(\frac{t}{x}\right)f(t)d\,t,\\
B_{-}^{\nu,\mu}f=\int_x^\infty \left( t^2-x^2\right)^{-\frac{\mu}{2}}P_\nu^\mu \left(\frac{t}{x}\right)f(t)d\,t,\\
\label{72}
E_{-}^{\nu,\mu}f=\int_x^\infty \left( t^2-x^2\right)^{-\frac{\mu}{2}}\mathbb{P}_\nu^\mu \left(\frac{x}{t}\right)f(t)d\,t.\\ \nonumber
\end{eqnarray}
Здесь $P_\nu^\mu(z)$---функция Лежандра первого рода [102], $\mathbb{P}_\nu^\mu(z)$---та же функция на разрезе $-1\leq  t \leq 1$, $f(x)$---локально суммируемая функция, удовлетворяющая некоторым ограничениям на рост при $x\to 0,x\to\infty$. Параметры $\mu,\nu$---комплексные числа, $\Re \mu <1$, можно ограничиться значениями $\Re \nu \geq -1/2$.

Интегральные операторы указанного вида с функциями Лежандра в ядрах впервые встретились в работах E.T.~Copson по уравнению Эйлера-Пуассона-Дарбу в конце 1950-х годов.
А именно, в работе [359] доказано следующее утверждение, которое я бы назвал

\textbf{Лемма Копсона.}

Рассмотрим дифференциальное уравнение в частных производных с двумя переменными:

$$
\frac{\pr^2 u(x,y)}{\pr x^2}+\frac{2\alpha}{x}\frac{\pr u(x,y)}{\pr x}=
\frac{\pr^2 u(x,y)}{\pr y^2}+\frac{2\beta}{y}\frac{\pr u(x,y)}{\pr y}
$$
(обобщённое уравнение Эйлера--Пуассона--Дарбу или В--гиперболическое уравнение по терминологии И.А.Киприянова) в открытой четверти плоскости $x>0, y>0$ при положительных параметрах $\beta>\alpha>0$ с краевыми условиями на осях координат (характеристиках)
$$
u(x,0)=f(x), u(0,y)=g(y), f(0)=g(0).
$$
Предполагается, что решение u(x,y) является непрерывно дифференцируемым в замкнутом квадранте, имеет непрерывные вторые производные в открытом квадранте, граничные функции $f(x), g(y)$ являются непрерывно дифференцируемыми.

Тогда, если решение поставленной задачи существует, то для него выполняются соотношения:

\begin{equation}
\label{Cop1}
\frac{\pr u}{\pr y}=0, y=0,  \frac{\pr u}{\pr x}=0, x=0,
\end{equation}
\begin{equation}
\label{Cop2}
2^\beta \Gamma(\beta+\frac{1}{2})\int_0^1 f(xt)t^{\alpha+\beta+1}
\lr{1-t^2}^{\frac{\beta -1}{2}}P_{-\alpha}^{1-\beta}{t}\,dt=
\end{equation}
\begin{equation*}
=2^\alpha \Gamma(\alpha+\frac{1}{2})\int_0^1 g(xt)t^{\alpha+\beta+1}
\lr{1-t^2}^{\frac{\alpha -1}{2}}P_{-\beta}^{1-\alpha}{t}\,dt,
\end{equation*}
$$
\Downarrow
$$
\begin{eqnarray}
\label{Cop3}
g(y)=\frac{2\Gamma(\beta+\frac{1}{2})}{\Gamma(\alpha+\frac{1}{2})
\Gamma(\beta-\alpha)}y^{1-2\beta}\cdot\nonumber\\
\cdot\int_0^y x^{2\alpha-1}f(x)
\lr{y^2-x^2}^{\beta-\alpha-1}x \,dx.
\end{eqnarray}

Соотношения (\ref{Cop1}) были известны ранее до Копсона, они очевидны. В работе приводится нестрогий вывод (\ref{Cop2}), то есть получено, что граничные функции (или значения решения на характеристиках) не могут быть произвольными, они связаны в современной терминологии операторами Бушмана--Эрдейи. Далее утверждается, что если две функции связаны операторами  Бушмана--Эрдейи указанных порядков, то на самом деле выполняется (\ref{Cop3})---то есть они связаны более простыми операторами Эрдейи--Кобера. В этом на мой взгляд основное содержание леммы Копсон
 а.

Но отсюда не следует, как иногда отмечается, что теперь можно сразу  получить обращение соответствующего оператора  Бушмана--Эрдейи, хотя бы формально. Для этого произвольную функцию в правой части соответствующего уравнения надо записать также в виде оператора  Бушмана--Эрдейи соответствующего порядка, чтобы подогнать под лемму Копсона. Но для этого уже надо уметь  оператор Бушмана--Эрдейи обращать--получается порочный круг. Таким образом, неверно приписывать Копсону первый результат по обращению операторов  Бушмана--Эрдейи, хотя на
 колько нам известно в его работе эти операторы действительно встречаются впервые.

Доказательства в работе [359] скорее являются нестрогими рассуждениями, намечающими что и в каком порядке надо делать, хотя всё понятно и видимо легко доводится до строгого. Видимо, у Копсона были некоторые сомнения в этих результатах, так как я не нашёл их в последующей его книге [360], что несколько странно. А может, он не придавал им особого значения. Отметим также, что в Красной Книге [34 --- Килбас, Маричев, Самко] и других работах даётся не совсем точная ссылка на работу [359], которую мы здесь исправили.

Эта первая работа Копсона нашла продолжение в совместной работе с Эрдейи [361]. Там даётся строгий вывод, вводятся подходящие классы функций, явно озвучена связь с дробными интегралами и операторами Кобера--Эрдейи.

Таким образом, можно сделать вывод, что впервые подробное изучение разрешимости и обратимости данных операторов было начато в 1960--х годах в работах    Р.~Бушмана и А.~Эрдейи [216--219]. Операторы Бушмана--Эрдейи изучались также в работах Higgins, Ta Li, Love, Habibullah, K.N.~Srivastava, Динь Хоанг Ань, Смирнова, Вирченко, Федотовой, Килбаса, Скоромник и др. При этом в основном изучались задачи о решении интегральных уравнений с этими операторами, их факторизации и обращения. Эти результаты изложены в монографии [34], хотя случай выбранных нами пределов интегрировани
  считается там особым и не рассматривается, за исключением одного набора формул композиции, см. также [220--225].

Термин   "операторы Бушмана--Эрдейи"\  как наиболее исторически оправданный был введён автором в [226--227], впоследствии он использовался и другими авторами. Ранее в [34] встречался предложенный О.И.Маричевым термин "операторы Бушмана". Наиболее полное изучение операторов Бушмана--Эрдейи на наш взгляд было проведено в работах автора в 1980--1990-е годы [229, 374, 226--228] и затем продолжено в [252--260, 362--372]. При этом необходимо отметить, что  роль операторов Бушмана--Эрдейи как ОП до работ [229, 374, 226--228] вообще ранее нигде не отмечалась и не рассматривалась.

Особо отметим работы Н.А.Вирченко [238--240], А.А.Килбаса  и их учеников. Так в работах А.А.Килбаса и О.В.Скоромник [375--376] рассматривается действие операторов Бушмана--Эрдейи в весовых пространствах Лебега, а также  многомерные обобщения в виде интегралов по пирамидальным областям.

Важность операторов Бушмана--Эрдейи во многом обусловлена их многочисленными приложениями. Например, они встречаются в следующих вопросах теории уравнений с частными производными [34]: при решении задачи Дирихле для уравнения Эйлера--Пуассона--Дарбу в четверти плоскости и установлении соотношений между значениями решений уравнения Эйлера--Пуассона--Дарбу на многообразии начальных данных и характеристике (см. лемму Копсона выше), теории преобразования Радона, так как в силу результатов Людвига [230--232] действие преобразования Радона пр
 и разложении по сферическим гармоникам сводится как раз к операторам  Бушмана--Эрдейи по радиальной переменной, при исследовании краевых задач для различных уравнений с существенными особенностями внутри области, доказательству вложения пространств И.А.~Киприянова  в весовые пространства С.Л.~Соболева, установлению связей между операторами преобразования и волновыми операторами теории рассеяния, обобщению классических интегральных представлений Сонина и Пуассона и операторов преобразования Сонина--Пуассона--Дельсарта.

Отметим, что полученные автором в [226--227] в 1990--1991 годах результаты по факторизации позволяют сразу  описать условия ограниченности общих операторов Бушмана--Эрдейи в весовых пространствах Лебега. Действительно, в этих работах найдены точные условия ограниченности операторов Бушмана--Эрдейи нулевого порядка гладкости в весовых пространствах Лебега. Как показано ниже, общие операторы являются композициями операторов нулевого порядка гладкости и операторов Римана--Лиувилля или Эрдейи--Кобера, а точные условия ограниченности последни
 х были известны уже давно [34]. На этом пути только не удаётся установить точные константы в неравенствах для норм операторов. Ниже приведены условия ограниченности операторов Бушмана--Эрдейи \textit{с точными константами} для невесового случая, и даже найдены условия на параметры, обеспечивающие их \textit{унитарность}. Другой подход к описанию действия операторов Бушмана--Эрдейи в весовых пространствах Лебега недавно был применён в [375--376].

Приведём полный список результатов, полученных автором, при изучении операторов Бушмана--Эрдейи.
\begin{itemize}
\item Введены операторы Бушмана--Эрдейи первого, второго и третьего родов, причём два последних класса ранее не вводились и не исследовались.
\item Впервые операторы Бушмана--Эрдейи изучены как операторы преобразования
\item Введено понятие гладкости для операторов Бушмана--Эрдейи.
\item Впервые изучены операторы Бушмана--Эрдейи с интегрированием по промежуткам $[0,x]$ и $[x,\infty]$.
\item Получены новые формулы факторизации.
\item Изучено действие операторов Бушмана--Эрдейи нулевого порядка гладкости в пространствах $L_2(0,\infty)$ и весовых пространствах Лебега. Получены точные значения норм, найдены  условия ограниченности и неограниченности. Ввиду полученных формул факторизации отсюда следуют результаты по действию в весовых пространствах Лебега общих операторов Бушмана--Эрдейи.
\item Найдены условия унитарности операторов Бушмана--Эрдейи нулевого порядка гладкости в пространствах $L_2(0,\infty)$.
\item Решена задача об унитарном обобщении операторов преобразования Сонина и Пуассона. Найдены ОП этого типа, унитарные при всех значениях параметра. Решение получено в рамках метода композиций, искомые ОП выражаются через  операторы Бушмана--Эрдейи третьего рода.
\item Описан спектр операторов Бушмана--Эрдейи нулевого порядка гладкости в весовых пространствах Лебега.
\item Найдены выражения операторов Бушмана--Эрдейи через классические интегральные преобразования Фурье, синус-- и косинус--преобразования Фурье, преобразование Ханкеля.
\item Изучены свойства мультипликаторов операторов Бушмана--Эрдейи при действии преобразования Меллина. Через мультипликаторы найдены условия, при которых произвольные операторы являются ОП типа Сонина и Пуассона, а также получены формулы связи разносторонних  операторов Бушмана--Эрдейи в терминах обобщённых преобразований Стильтьеса и Гильберта.
\item Введены и изучены операторы Бушмана--Эрдейи второго рода, с функциями Лежандра второго рода в ядрах.
\item Получены обобщения интегральных представлений Сонина и Пуассона для специальных функций.
\item Решено большое число интегро--дифференциальных уравнений с функциями Лежандра различных типов в ядрах и получены оценки решений в весовых пространствах Лебега.
\item С помощью операторов Бушмана--Эрдейи нулевого порядка гладкости в одномерном случае доказано вложение пространств И.А.~Киприянова в весовые пространства С.Л.~Соболева.
\item Предложен общий метод построения новых начальных и краевых задач с интегральными операторами типа дробного интегродифференцирования различных типов в краевых условиях для уравнения Эйлера--Пуассона--Дарбу и других уравнений $B$--эллиптического, $B$--гиперболического и $B$--параболического типов.
\item Изучено действие операторов Бушмана--Эрдейи в пространствах аналитических функций, доказаны коэффициентные оценки типа теорем искажения и решены интегро--дифференциальные уравнения дробного порядка в комплексной плоскости.
\end{itemize}

Приведём основные результаты, в основном следуя в изложении [226--227]. Все рассмотрения ведутся ниже на полуоси. Поэтому будем обозначать через $L_2$ пространство $L_2(0, \infty)$ и $L_{2, k}$ весовое пространство $L_{2, k}(0, \infty)$. $\mathbb{N}$ обозначает множество натуральных, $\mathbb{N}_0$ -- неотрицательных целых, $\mathbb{Z}$
-- целых и $\mathbb{R}$ -- действительных чисел.

Вначале распространим определение 6 на  важный не исследованный ранее случай $\mu =1$.

\textsc{Определение 7.} Введём при $\mu =1$ операторы  Бушмана--Эрдейи нулевого порядка гладкости по формулам
\begin{eqnarray}
\label{73}
B_{0+}^{\nu,1}f=\frac{d}{dx}\int_0^x P_\nu \left(\frac{x}{t}\right)f(t)\,dt,\\
\label{731}
E_{0+}^{\nu,1}f=\int_0^x P_\nu \left(\frac{t}{x}\right)\frac{df(t)}{dt}\,dt,\\
\label{732}
B_{-}^{\nu,1}f=\int_x^\infty P_\nu \left(\frac{t}{x}\right)(-\frac{df(t)}{dt})\,dt,\\
\label{733}
E_{-}^{\nu,1}f=(-\frac{d}{dx})\int_x^\infty P_\nu \left(\frac{x}{t}\right)f(t)\,dt,
\end{eqnarray}
где $P_\nu(z)=P_\nu^0(z)$---функция Лежандра.

Разумеется, при очевидных дополнительных условиях на функции в (\ref{73})--(\ref{733}) можно продифференцировать под знаком интеграла или проинтегрировать по частям.

\Theorem{ 10} Справедливы следующие формулы факторизации операторов Бушмана--Эрдейи на подходящих функциях через дробные интегралы Римана--Лиувилля, Эрдейи--Кобера и Бушмана--Эрдейи нулевого порядка гладкости:
\be{1.9}{B_{0+}^{\nu,\,\mu} f=I_{0+}^{1-\mu}~ {_1 S^{\nu}_{0+}f},~B_{-}^{\nu, \,\mu} f={_1 P^{\nu}_{-}}~ I_{-}^{1-\mu}f,}
\be{1.10}{E_{0+}^{\nu,\,\mu} f={_1 P^{\nu}_{0+}}~I_{0+}^{1-\mu}f,~E_{-}^{\nu, \, \mu} f= I_{-}^{1-\mu}~{_1 S^{\nu}_{-}}f.}
\Endproc

Эти формулы позволяют "разделить" параметры $\nu$ и $\mu$. Мы докажем, что операторы
\eqref{73}--\eqref{733} являются изоморфизмами пространств $L_2(0, \infty)$, если $\nu$
не равно некоторым исключительным значениям. Поэтому операторы \eqref{71}--\eqref{72}
по действию в пространствах типа $L_2$ в определённом смысле подобны операторам дробного интегродиффенцирования  $I^{1-\mu}$, с которыми они совпадают при $\nu=0$. Далее операторы Бушмана--Эрдейи будут доопределены при всех значениях $\mu$.
Исходя из этого, введём следующее

\textsc{Определение 8.} Число $\rho=1-Re\,\mu $ назовём порядком гладкости операторов Бушмана--Эрдейи \eqref{71}--\eqref{72}.

Таким образом, при $\rho > 0$ (то есть при $Re\, \mu > 1$) операторы Бушмана--Эрдейи
являются сглаживающими, а при $\rho < 0$ (то есть при $Re\, \mu < 1$) уменьшающими
гладкость в пространствах типа $L_2 (0, \infty)$. Операторы \eqref{73}--\eqref{733},
для которых $\rho = 0$, являются по данному определению операторами нулевого порядка гладкости.

Перечислим основные свойства операторов Бушмана--Эрдейи первого рода
\eqref{71}--\eqref{72} с функцией Лежандра I рода в ядре. Будем обозначать области определения операторов через $\mathfrak{D}(B_{0+}^{\nu , \, \mu })$,
$\mathfrak{D}(E_{0+}^{\nu , \, \mu })$ и т.д.

Простейшие свойства функций Лежандра приводят к тождествам
\begin{eqnarray}
& B_{0+}^{\nu, \, \mu} f = B_{0+}^{-\nu-1,\mu} f, & E_{0+}^{\nu,\, \mu} f=E_{0+}^{-\nu-1, \, \mu}f,
\nonumber \\
& B _{-}^{\nu, \, \mu} f=B_{-}^{-\nu-1, \, \mu}f, & E_{-}^{\nu, \,\mu} f = E_{-}^{-\nu-1, \, \mu} f,
\label{2.1}
\end{eqnarray}

\begin{eqnarray}
& & (2 \nu +1)x \, B_{0+}^{\nu,\, \mu} [ \frac{1}{x}f]=(\nu-\mu+1)B_{0+}^{\nu+1, \, \mu}
f +  (\nu+\mu)B_{0+}^{\nu-1, \, \mu}f ,\nonumber \\
& & (2 \nu +1) \frac{1}{x} \, B_{-}^{\nu,\mu} [x f]=(\nu-\mu+1)B_{-}^{\nu+1,\mu}
f +  (\nu+\mu)B_{-}^{\nu-1,\mu}f, \label{2.2}
\end{eqnarray}

\begin{eqnarray}
& & B_{0+}^{\nu-1, \, \mu}f - B_{0+}^{\nu+1, \, \mu}f = -(2 \nu +1)B_{0+}^{0, \, \mu-1} [ \frac{1}{x}f]  \nonumber ,\\
& & B_{-}^{\nu-1, \, \mu}f - B_{-}^{\nu+1, \, \mu}f = -(2 \nu +1)\frac{1}{x}B_{-}^{\nu, \, \mu-1} f. \label{2.3}
\end{eqnarray}
Из формул разложения $L_{\nu}$ на множители получаются тождества

\begin{eqnarray}
& & B_{0+}^{\nu, \, \mu -1}\left(\frac{d}{dx}- \frac{\nu}{x}\right) f = B_{0+}^{\nu-1, \, \mu } f,  \label{2.4}  ,\\
& & B_{0+}^{\nu, \, \mu -1}\left(\frac{d}{dx}+ \frac{\nu}{x}\right) f = B_{0+}^{\nu+1,  \mu } f, \label{2.5}
\end{eqnarray}
справедливые при условиях $Re\,\mu < 1$, $Re\,\nu > - \frac{1}{2}$
$$
\lim\limits_{y \to 0} f(y) / y^{\nu}=0, ~ f \in \mathfrak{D}(B_{0+}^{\nu \pm 1, \, \mu }),~ \left(\frac{d}{dx}\pm \frac{\nu}{x}\right)f \in \mathfrak{D}(B_{0+}^{\nu,  \, \mu -1})
$$

Формулы \eqref{2.1} позволяют ограничиться случаем $Re\,\nu \geq - \frac{1}{2}$.
Функции, на которые действуют операторы, должны принадлежать их областям определения. Для оператора $E_{0+}^{\nu,\, \mu}$ справедливы те же формулы, что и для $B_{0+}^{\nu,\, \mu}$

\Theorem{ 11} Операторы Бушмана-Эрдейи \eqref{71}-\eqref{72} определены,
если $Re ~ \mu < 1$ или $\mu \in \mathbb{N}$ и дополнительно выполнены условия:

а) для оператора  $B_{0+}^{\nu, \, \mu}$
$$
\int\limits_0^x \sq y \,|f(y) \ln y|\, dy~ < \infty,
$$
если $\nu=-\frac{1}{2}$, $\mu \neq \frac{1}{2}$, а во всех остальных случаях
$$
\int\limits_0^x y^{-Re\,\nu} \,|f(y)|\, dy~ < \infty,
$$

б) для оператора $E_{0+}^{\nu, \, \mu}$ дополнительные условия не требуются;

в) для оператора $E_{-}^{\nu, \, \mu}$
$$
\int\limits_x^{\infty} y^{-Re\,\nu} \,|f(y)|\, dy~ < \infty,
$$

г) для оператора $B_{-}^{\nu, \, \mu}$
$$
\int\limits_x^{\infty} y^{-\frac{1}{2}-Re\,\nu} \,|\ln y \cdot f(y)|\, dy~ < \infty,
$$
если $\nu=-\frac{1}{2}$, $\mu \neq \frac{1}{2}$, а во всех остальных случаях
$$
\int\limits_x^{\infty} y^{Re(\nu-\mu)} \,|f(y)|\, dy~ < \infty.
$$
\Endproc

В этой теореме предполагается, что функция $f(x)$ является локально суммируемой на
$(0, \infty)$, $x$ -- произвольное положительное число.

При некоторых специальных значениях параметров $\nu,~\mu$ операторы Бу\-шмана--Эрдейи
сводятся к более простым. Так при значениях $\mu=-\nu$ или $\mu=\nu+2$ они являются
операторами Эрдейи--Кобера; при $\nu = 0$ операторами дробного интегродифференцирования
$I_{0+}^{1-\mu}$ или $I_{-}^{1-\mu}$; при $\nu=-\frac{1}{2}$, $\mu=0$ или $\mu=1$
ядра выражаются через эллиптические интегралы; при  $\mu=0$,  $x=1$, $v=it-\frac{1}{2}$  оператор $B_{-}^{\nu, \, 0}$ лишь на постоянную отличаются от преобразования Малера--Фока.

\Theorem{ 12} Пусть или $Re \, \mu < 0$, или $\mu = m \in \mathbb{N}$,
$-m \leq \nu \leq m-1$, $\nu \in \mathbb{Z}$. Тогда справедливы тождества

\begin{eqnarray}
& \frac{d}{dx} B_{0+}^{\nu, \, \mu } f = B_{0+}^{\nu, \, \mu +1 } f, &
E_{0+}^{\nu, \, \mu }\frac{d \, f}{dx} = E_{0+}^{\nu, \, \mu +1 } f  \label{2.6}  ,\\
& B_{-}^{\nu, \, \mu } \lr{-\frac{d \, f}{dx}} = B_{-}^{\nu, \, \mu + 1}f, &
\lr{-\frac{d}{dx}} E_{-}^{\nu, \, \mu } f = E_{-}^{\nu, \, \mu +1 } f \label{2.7}.
\end{eqnarray}
если все указанные операторы определены.
\Endproc

Эта теорема позволяет доопределить операторы Бушмана--Эрдейи и на значения $Re \mu \geq 1$, переопределив их для натуральных $\mu$.

\textsc{Определение 9.} Пусть дано число $\sigma$, $Re \, \sigma \geq 1.$
Обозначим через $m$ наименьшее натуральное число, такое, что $\sigma= \mu +m$,
$Re \, \mu <1$. Тогда операторы Бушмана--Эрдейи доопределим по формулам

\begin{eqnarray}
& & B_{0+}^{\nu, \, \sigma}=B_{0+}^{\nu, \, \mu + m}=\lr{\frac{d}{dx}}^m \,
B_{0+}^{\nu, \, \mu}, \nonumber\\
& & E_{0+}^{\nu, \, \sigma}=E_{0+}^{\nu, \, \mu + m}=
E_{0+}^{\nu, \, \mu} \, \lr{\frac{d}{dx}}^m, \label{2.9}\\
& & B_{-}^{\nu, \, \sigma}=B_{-}^{\nu, \, \mu + m}=
B_{-}^{\nu, \, \mu} \lr{-\frac{d}{dx}}^m, \nonumber\\
& & E_{-}^{\nu, \, \sigma}=E_{-}^{\nu, \, \mu + m}=\lr{-\frac{d}{dx}}^m
E_{-}^{\nu, \, \mu}. \nonumber
\end{eqnarray}

Отметим, что при натуральных $\mu$ операторы Бушмана--Эрдейи определены и по формулам
\eqref{71}--\eqref{72}. Мы переопределяем их для этих значений $\mu$ по формуле
\eqref{2.9}. Таким образом, символами $B_{0+}^{\nu, \, \mu}$, $E_{0+}^{\nu, \, \mu}$,
$B_{-}^{\nu, \, \mu}$, $E_{-}^{\nu, \, \mu}$ далее мы будем обозначать операторы, определяемые
по формулам \eqref{71}--\eqref{72} при $Re \, \mu < 1$,  и по формулам \eqref{2.9}
при $Re \, \mu \geq 1$.

Будем рассматривать наряду с оператором Бесселя также тесно связанный с ним дифференциальный оператор
\begin{equation}
\label{75}
L_{\nu}=D^2-\frac{\nu(\nu+1)}{x^2}=\left(\frac{d}{dx}-\frac{\nu}{x}\right)
\left(\frac{d}{dx}+\frac{\nu}{x}\right),
\end{equation}
который при $\nu \in \mathbb{N}$ является оператром углового момента из квантовой механики.
Их взаимосвязь устанавливает
\Theorem{ 12} Пусть пара ОП $X_\nu, Y_\nu$ сплетают $L_{\nu}$ и вторую  производную:
\begin{equation}
\label{76}
X_\nu L_{\nu}=D^2 X_\nu , Y_\nu D^2 = L_{\nu} Y_\nu.
\end{equation}
Введём новую пару ОП по формулам
\begin{equation}
\label{77}
S_\nu=X_{\nu-1/2} x^{\nu+1/2}, P_\nu=x^{-(\nu+1/2)} Y_{\nu-1/2}.
\end{equation}
Тогда пара новых ОП $S_\nu, P_\nu$ сплетают оператор Бесселя и вторую производную:
\begin{equation}
\label{78}
S_\nu B_\nu = D^2 S_\nu, P_\nu D^2 = B_\nu P_\nu.
\end{equation}
\Endproc

Разумеется, по указанным формулам можно перейти и наоборот от ОП для оператора Бесселя к ОП для оператора углового момента. Мы сохраним за ОП, действующим по формулам (\ref{76}), названия ОП типа Сонина и Пуассона соответственно.

Определим класс $\Phi(B_{0+}^{\nu, \, \mu})$ как множество функций таких, что

1) $f(x) \in \mathfrak{D}(B_{0+}^{\nu, \, \mu})\bigcap C^2(0, \infty),$

2) $\lim\limits_{y \to 0}\left| \frac{\ln y}{\sq y} f(y)+\sq y \ln y \cdot f'(y) \right|=0,$

если $\nu=-\frac{1}{2}$, $ \mu \neq \frac{1}{2}$;
$$
\lim\limits_{y \to 0}[ (\nu+1) y^{\nu} f(y)- y^{\nu+1} f'(y) ]=0,
$$
если $\mu= \nu +1$, $ Re \, \nu \neq - \frac{1}{2}$; и, наконец,
$$
\lim\limits_{y \to 0}\left( \nu \frac{f(y)}{y^{\nu+1}}+\frac{f'(y)}{y^{\nu}} \right)=0
$$
во всех остальных случаях.

\Theorem{ 13} Пусть$f(x) \in \Phi(B_{0+}^{\nu, \, \mu})$,
$Re \, \mu \leq -1$. Тогда оператор $B_{0+}^{\nu, \, \mu}$ является оператором преобразования типа Сонина и удовлетворяет соотношению \eqref{76} на функциях $f(x)$.
\Endproc

Аналогичный результат справедлив и для других операторов Бушмана--Эрдейи. При этом
$E_{-}^{\nu, \, \mu}$ также является оператором типа Сонина, а $E_{0+}^{\nu, \, \mu}$
и $B_{-}^{\nu, \, \mu}$ -- операторами типа Пуассона.

Можно рассматривать случай, когда нижний предел в соответствующих интегралах \eqref{71}--\eqref{72}
равен произвольному числу $a>0$, или верхний предел в интегралах
равен произвольному конечному числу $b>0$. При этом все результаты этого пункта сохраняются, а их формулировки значительно упрощаются. В частности, все операторы Бушмана--Эрдейи в этом
случае определены при единственном условии $Re \mu < 1$ в форме \eqref{71}--\eqref{72}
и являются операторами преобразования на функциях $f(x)$ таких, что $f(a)=f'(a)=0$ или
($f(b)=f'(b)=0$).

Теперь сделаем важное замечание. Из полученной теоремы следует, что ОП Бушмана--Эрдейи связывают собственные функции операторов Бесселя и второй производной. Таким образом, половина ОП Бушмана--Эрдейи переводят тригонометрические или экспоненциальные функции в приведённые функции Бесселя, а другая половина наоборот. Эти формулы здесь не приводятся, их нетрудно выписать явно. Все они являются обобщениями исходных формул Сонина и Пуассона  (\ref{54}--\ref{55}) и представляют существенный интерес. Ещё раз отметим, что подобные формулы являютс
  непосредственными следствиями доказанных сплетающих свойств ОП Бушмана--Эрдейи, и могут быть непосредственно проверены при помощи таблиц интегралов от специальных функций.

Перейдём к вопросу о различных факторизациях операторов Бушмана--Эрдейи через операторы
Эрдейи--Кобера и дробные интегралы Римана--Лиувилля. Напомним, что операторы Эрдейи--Кобера определяются при $\alpha > 0$ по формулам

\begin{eqnarray*}
& & I_{0+;\, 2,\, y}^{\alpha} f = \frac{2}{\Gamma(\alpha)}x^{-2(\alpha+y)}
\int\limits_0^x (x^2-t^2)^{\alpha-1}t^{2y+1}f(t)\,dt, \\
& & I_{-;\, 2,\, y}^{\alpha} f = \frac{2}{\Gamma(\alpha)}x^{2y}
\int\limits_x^{\infty} (t^2-x^2)^{\alpha-1}t^{2(1-\alpha-y)-1}f(t)\,dt,
\end{eqnarray*}
а при значениях $\alpha > -n$, $n \in \mathbb{N}$ по формулам

\begin{eqnarray}
& & I_{0+;\, 2, y}^{\alpha} f =x^{-2(\alpha+y)} {\lr{\frac{d}{d x^2}}}^n  x^{2(\alpha
+y+n)}I^{\alpha+n}_{0+; \, 2,\, y}f \label{2.15} \\
& & I_{-;\, 2, y}^{\alpha} f = x^{2y} {\lr{-\frac{d}{d x^2}}}^n  x^{2(\alpha
-y)}I^{\alpha+n}_{-; \, 2,\, y-n}f \label{2.16}.
\end{eqnarray}

\Theorem{ 14} Справедливы следующие формулы факторизации операторов
Бушмана--Эрдейи через операторы дробного интегродифференцирования и Эрдейи--Кобера:
\begin{eqnarray}
& & B_{0+}^{\nu, \, \mu}=I_{0+}^{\nu+1-\mu} I_{0+; \, 2, \, \nu+ \frac{1}{2}}^{-(\nu+1)} {\lr{\frac{2}{x}}}^{\nu+1}\label{2.17}, \\
& & E_{0+}^{\nu, \, \mu}= {\lr{\frac{x}{2}}}^{\nu+1} I_{0+; \, 2, \, - \frac{1}{2}}^{\nu+1} I_{0+}^{-(\nu+\mu)}  \label{2.18}, \\
& & B_{-}^{\nu, \, \mu}= {\lr{\frac{2}{x}}}^{\nu+1}I_{-; \, 2, \, \nu+ 1}^{-(\nu+1)} I_{-}^{\nu - \mu+2}  \label{2.19}, \\
& & E_{-}^{\nu, \, \mu}= I_{-}^{-(\nu+\mu)} I_{-; \, 2, \, 0} ^{\nu+1} {\lr{\frac{x}{2}}}^{\nu+1}  \label{2.20}.
\end{eqnarray}
\Endproc

Все основные свойства операторов Бушмана--Эрдейи могут быть выведены из теоремы 14.
Так можно получить, что формально обратным к оператору Бушмана--Эрдейи с параметрами ($\nu,~\mu$)
является тот же оператор с параметрами ($\nu, ~2 - \mu$). При этом из двух операторов --- прямого и обратного --- всегда один будет иметь интегральное представление \eqref{71}--\eqref{72}, а другой определяется формулами \eqref{2.9}; один будет обязательно иметь положительный порядок гладкости, а другой отрицательный (кроме операторов нулевого порядка гладкости). Кроме того, учитывая \eqref{2.15}--\eqref{2.16}, теорема 14 позволяет доопределить операторы Бушмана--Эрдейи
на всю область значений параметров. Такое доопределение согласуется с \eqref{2.9}. Отметим, что факторизации \eqref{2.17}--\eqref{2.20} являются новыми по сравнению с факторизациями, приведёнными в [34].

Рассмотрим связь между операторами Бушмана--Эрдейи и сплетающими операторами Сонина--Пуассона--Дельсарта (СПД). Мы предпочтём дать новые определения для них, чтобы сохранить однообразие обозначений в этом пункте.

\textsc{Определение 10.} Введём операторы преобразования Сонина-- \\ Пуассона--Дельсарта по формулам

\begin{eqnarray}
& & {_0S_{0+}^{\nu}}=B^{\nu, \, \nu+2}_{0+}= I_{0+; \, 2, \, \nu+ \frac{1}{2}}^{-(\nu+1)} {\lr{\frac{2}{x}}}^{\nu+1} \label{2.21} \\
& & {_0P_{0+}^{\nu}}=E^{\nu, \, - \nu}_{0+}=  {\lr{\frac{x}{2}}}^{\nu+1} I_{0+; \, 2, \, - \frac{1}{2}}^{\nu+1} \label{2.22} \\
& & {_0P_{-}^{\nu}}=B^{\nu, \, \nu+2}_{-}= {\lr{\frac{2}{x}}}^{\nu+1}I_{-; \, 2, \, \nu+ 1}^{-(\nu+1)} \label{2.23} \\
& & {_0S_{-}^{\nu}}=E^{\nu, \, - \nu}_{-}= I_{-; \, 2, \, 0} ^{\nu+1} {\lr{\frac{x}{2}}}^{\nu+1} \label{2.24}
\end{eqnarray}
Это определение в сочетании с \eqref{2.15}--\eqref{2.16} приводит к следующим интегральным представлениям:

$$
{_0S_{0+}^{\nu}} f = \frac{2^{\nu+2}}{\Gamma(-\nu-1)}x \int\limits_0^x (x^2-t^2)^{-\nu-2}t^{\nu+1}f(t)\,dt, ~ Re \, \nu < -1,
$$

$$
= \frac{2^{\nu+1}}{\Gamma(-\nu)} \frac{d}{dx} \int\limits_0^x (x^2-t^2)^{-\nu-1}t^{\nu+1}f(t)\,dt, ~ Re \, \nu < 0;
$$

$$
{_0P_{0+}^{\nu}} f = \frac{1}{2^{\nu} \Gamma(\nu+1)}x^{-\nu} \int\limits_0^x (x^2-t^2)^{\nu}f(t)\,dt, ~ Re \, \nu > -1,
$$

$$
= \frac{1}{2^{\nu} \Gamma(\nu+2)}\frac{1}{x^{\nu+1}} \frac{d}{dx} \int\limits_0^x (x^2-t^2)^{\nu+1}f(t)\,dt, ~ Re \, \nu > -2;
$$

$$
{_0P_{-}^{\nu}} f = \frac{2^{\nu+2}}{\Gamma(-\nu-1)}x^{\nu +1} \int\limits_x^{\infty} (t^2-x^2)^{-\nu-2}tf(t)\,dt, ~ Re \, \nu < -1,
$$

$$
= \frac{2^{\nu+1}}{\Gamma(-\nu)} x^{\nu} \lr{-\frac{d}{dx}} \int\limits_x^{\infty} (t^2-x^2)^{-\nu-1}tf(t)\,dt, ~ Re \, \nu < 0;
$$

$$
{_0S_{-}^{\nu}} f = \frac{1}{2^{\nu} \Gamma(\nu+1)}\int\limits_x^{\infty} (t^2-x^2)^{\nu}t^{-\nu}f(t)\,dt, ~ Re \, \nu > -1,
$$

$$
= \frac{1}{2^{\nu+1} \Gamma(\nu+2)} \lr{-\frac{1}{x} \frac{d}{dx}} \int\limits_x^{\infty} (t^2-x^2)^{\nu+1}t^{-\nu}f(t)\,dt, ~ Re \, \nu > -2.
$$

Эти операторы являются сплетающими типа Сонина или Пуассона. Если построить новые ОП для оператора углового момента (см. выше), то получим операторы типа Сонина
$$
X_{\nu} f= {_0S_{0+}^{\nu- \frac{1}{2}}} x^{\nu} f= \frac{2^{\nu+\frac{3}{2}}}{\Gamma(-\nu-\frac{1}{2})}x \int\limits_0^x (x^2-t^2)^{-\nu-\frac{3}{2}}t^{2 \nu+1}f(t)\,dt
$$

если $Re \, \nu < -1/2$, а если $Re \, \nu < 1/2$, то
\be{2.25}{X_{\nu} f= S_{\nu} f= \frac{2^{\nu+\frac{1}{2}}}{\Gamma(\frac{1}{2}-\nu)}\frac{d}{dx} \int\limits_0^x (x^2-t^2)^{-\nu-\frac{1}{2}}t^{2 \nu+1}f(t)\,dt}

Аналогично получим оператор типа Пуассона вида

\be{2.26}{Y_{\nu} f= P_{\nu} f= \frac{1}{2^{\nu}\Gamma(\nu+1)}\frac{1}{x^{2 \nu}} \int\limits_0^x (x^2-t^2)^{\nu-\frac{1}{2}}f(t)\,dt}
при условии $Re \, \nu > -1/2$.

Оператор Бесселя является радиальной частью лапласиана в $\mathbb{R}^n$. Операторы преобразования, сплетающие $L_{\nu}$ и вторую производную, сводятся к операторам Эрдейи--Кобера, то есть к дробным производным по $x^2$. В пространствах более сложной  чем для
$\mathbb{R}^n$  геометрии также существуют подобные операторы преобразования. Они являются, например, отрицательными степенями оператора $d/{d\ch x}$  или оператора $d/{ d\cos x}$. Подобные операторы изучались в работах В.Я.~Ярославцевой [34], В.В.~Катрахова и его учеников.

Отметим, что из теоремы 14 выводятся и формулы \eqref{1.9}--\eqref{1.10}.

Перейдём теперь к изучению операторов \eqref{73}--\eqref{733}. Отметим, что если функция $f(x)$
допускает дифференцирование под знаком интеграла или интегрирование по частям, то операторы \eqref{73}--\eqref{733} принимают вид

\begin{eqnarray}
& & {_1S_{0+}^{\nu}}f=f(x)+\int\limits_0^x \frac{\partial}{\partial x}P_{\nu}(\frac{x}{y})f(y)dy, \label{3.1} \\
& & {_1P_{0+}^{\nu}}f=f(x)-\int\limits_0^x \frac{\partial}{\partial y}P_{\nu}(\frac{y}{x})f(y)dy, \label{3.2} \\
& & {_1P_{-}^{\nu}}f=f(x)+\int\limits_x^{\infty} \frac{\partial}{\partial y}P_{\nu}(\frac{y}{x})f(y)dy, \label{3.3} \\
& & {_1S_{-}^{\nu}}f=f(x)-\int\limits_x^{\infty} \frac{\partial}{\partial x}P_{\nu}(\frac{x}{y})f(y)dy. \label{3.4}
\end{eqnarray}
При этом для справедливости \eqref{3.2} и \eqref{3.3} соответственно дополнительно необходимы условия

$$
\lim\limits_{x \to 0} P_{\nu}(0) f(x) = 0,~\lim\limits_{x \to \infty} P_{\nu}(x) f(x) = 0.
$$

В дополнение к теореме 13 можно доказать, что при определённых условиях на функции операторы \eqref{73}--\eqref{733} являются операторами преобразования. Они сплетают оператор углового момента и вторую производную.
Из теоремы 14 вытекают следующие факторизации для операторов  Бушмана--Эрдейи нулевого порядка гладкости

\begin{eqnarray}
& & {_1S_{0+}^{\nu}}= I_{0+}^{\nu+1} I_{0+; \, 2, \, \nu+ \frac{1}{2}}^{-(\nu+1)} {\lr{\frac{2}{x}}}^{\nu+1} \label{3.5} \\
& & {_1P_{0+}^{\nu}}= {\lr{\frac{x}{2}}}^{\nu+1} I_{0+; \, 2, \, - \frac{1}{2}}^{\nu+1} I_{0+}^{-(\nu+1)} \label{3.6} \\
& & {_1P_{-}^{\nu}}= {\lr{\frac{2}{x}}}^{\nu+1}I_{-; \, 2, \, \nu+ 1}^{-(\nu+1)} I_{-}^{\nu+1} \label{3.7} \\
& & {_1S_{-}^{\nu}}= I_{-}^{-(\nu+1)} I_{-; \, 2, \, 0} ^{\nu+1} {\lr{\frac{x}{2}}}^{\nu+1} \label{3.8}
\end{eqnarray}

Теперь рассмотрим более подробно свойства ОП Бушмана--Эрдейи нулевого порядка гладкости, введённых по формулам (\ref{73}). Подобный оператор был построен В.~В.~Катраховым [233--234] путём домножения стандартного ОП Сонина на обычный дробный интеграл с целью взаимно компенсировать гладкость этих двух операторов и получить новый, который бы действовал в одном пространстве типа $L_2(0,\infty)$. Как впоследствии оказалось, это можно сделать известными средствами, так как ОП Сонина---это частный случай операторов Эрдейи--Кобера. Существует замечательная
  теорема А.~Эрдейи, позволяющая выделить стандартный дробный интеграл Римана--Лиувилля из дробного интеграла по любой функции [34]. В результате получается

\Theorem{ 15} Рассмотрим оператор дробного интегродифференцирования Эрдейи--Кобера по функции $g(x)=x^2$

$$
I_{0+; \, x^2}^{\alpha} f = \frac{1}{\Gamma(\alpha)} \int\limits_0^x (x^2-t^2)^{\alpha-1} 2t \cdot f(t)\,dt
$$
при значениях $Re\, \alpha > 0 $. Тогда при $0< Re\, \alpha  <1$ справедливо представление

\be{3.9}{I_{0+; \, x^2}^{\alpha} f = I_{0+}^{\alpha} \left[ (2x)^{\alpha} f(x) +\int_0^x \frac{\partial}{\partial x} P_{- \alpha}(\frac{x}{s})(2s)^{\alpha} f(s)\,ds \right],}
где $I_{0+}^{\alpha}$ --- обычный дробный интеграл Римана--Лиувилля.
\Endproc

Операторы нулевого порядка гладкости выделяются тем, что только для них можно доказать оценки в \textit{одном} пространстве типа $L_p(0,\infty)$. При этом, учитывая структуру этих операторов, удобно пользоваться техникой преобразования Меллина.

Напомним [235], что преобразованием Меллина функции $f(x)$ называется функция $g(s)$, которая определяется по формуле
\begin{equation}
\label{710}
g(s)=M[f](s)=\int_0^\infty x^{s-1} f(x)\,dx.
\end{equation}
Определим также свёртку Меллина
\begin{equation}
\label{711}
(f_1*f_2)(x)=\int_0^\infty  f_1\left(\frac{x}{y}\right) f_2(y)\,\frac{dy}{y},
\end{equation}
при этом оператор свёртки с ядром $K$ действует в образах преобразования Меллина как умножение на мультипликатор
\begin{eqnarray}
\label{712}
M[Af](s)=\int_0^\infty  K\left(\frac{x}{y}\right) f(y)\,\frac{dy}{y}=M[K*f](s)
=m_A(s)M[f](s),\\\nonumber m_A(s)=M[K](s).\phantom{1111111111111111111}
\end{eqnarray}

Заметим, что преобразование Меллина является обобщённым преобразованием Фурье на полуоси по мере Хаара  $\frac{dy}{y}$ [191]. Его роль велика в теории специальных функций, например, гамма--функция является преобразованием Меллина экспоненты. С преобразованием Меллина связан важный прорыв в 1970--х годах, когда в основном усилиями О.~И.~Маричева была полностью доказана и приспособлена для нужд вычисления интегралов известная теорема Джоан Люси Слейтер, позволяющая для большинства образов преобразований Меллина восстановить оригинал в явном ви
 е по простому алгоритму через гипергеометрические функции [235--236]. Этот результат, который несложно получить формально по общей формуле обращения Меллина--Барнса через вычеты, для своего строгого обоснования потребовал достаточно сложных и тщательных выкладок, связанных с обработкой асимптотик гипергеометрических функций вблизи полюсов и на бесконечности, а такие асимптотики весьма разнообразны и разнородны. Эта работа была только начата Люси Джоан Слейтер, а в основном проведена до конца Олегом Игоревичем Маричевым, данное обст
 ятельство часто недооценивается. Теорема Слейтер--Маричева позволила создать универсальный мощный метод вычисления интегралов, который впоследствии позволил решить многие задачи в теории уравнений с частными производными, а также воплотился в передовые технологии символьного интегрирования пакета MATHEMATICA (О.~И.~Маричев работает в Wolfram Research).

Меллин использовал введённое им преобразование для получения  формулы в виде простого явного ряда для решений  алгебраического уравнения произвольной степени. Эти формулы Меллина до сих пор остаются практически неизвестными большинству математиков, запуганных в этом вопросе специалистами по схоластической алгебре и теории групп. Вместе с тем методы решения уравнений произвольной степени на основе формулы Лагранжа для разложения обратной функции в ряд и теории гипергеометрических функций были известны ещё в 19 веке. Формулы  обра
 ения преобразования Меллина для основных специальных функций, которые мы теперь называем интегралами Меллина--Барнса, были получены в 1888 году итальянским математиком Сальваторе Пинкерле. Существует и также практически неизвестная теория вещественного обращения преобразований Лапласа и Меллина без выхода в комплексную плоскость и без использования формул с производными всё более высоких порядков, но мы не будем на этом останавливаться.

Добавим, что первый из финских математиков Хъялмар Меллин был интересной личностью. По своим убеждениям он был ярым анти--шведским националистом. Для такой позиции было достаточно оснований в то время, например, в конце 19 века все профессорские кафедры в Финляндии занимали шведы, и  Меллин был первым, как мы теперь говорим, лицом финской национальности, ставшим на своей родине профессором.

Для изучения операторов вида (\ref{712}) автором в [226--227] был предложен удобный алгебраический подход, который не содержит ничего нового, но в удобной форме позволяет быстро получать нужные оценки. Полезные факты будут собраны вместе как \\

\Theorem{ 16}  Пусть оператор $A$ действует по формуле (\ref{712}).
Тогда\\
а) Для того, чтобы он допускал расширение до ограниченного оператора в $L_2(0,\infty)$ необходимо и достаточно, чтобы
\begin{equation}
\label{713}
\sup_{\xi\in\mathbb{R}} |m_A(i\xi+\frac{1}{2})|=M_2<\infty,
\end{equation}
при этом $\|A\|_{L_2}=M_2$.\\
б) Для того, чтобы он допускал расширение до ограниченного оператора в $L_p(0,\infty), p>1$ при дополнительном условии неотрицательности ядра $K$ необходимо и достаточно, чтобы
\begin{equation}
\label{714}
\sup_{\xi\in\mathbb{R}} |m_A(i\xi+\frac{1}{p})|=M_p<\infty,
\end{equation}
при этом $\|A\|_{L_p}=M_p$.\\
в) Обратный оператор $A^{-1}$ действует также по формуле (\ref{712}) с мультипликатором $\frac{1}{m_A}$, для того, чтобы он допускал расширение до ограниченного оператора в $L_2(0,\infty)$ необходимо и достаточно, чтобы
\begin{equation}
\label{715}
\inf_{\xi\in\mathbb{R}} |m_A(i\xi+\frac{1}{2})|=m_2<\infty,
\end{equation}
при этом $\|A^{-1}\|_{L_2}=\frac{1}{m_2}$.\\
г) Пусть операторы $A,A^{-1}$ определены и ограничены в $L_2(0,\infty)$. Они унитарны тогда и только тогда, когда выполнено равенство
\begin{equation}
\label{716}
|m_A(i\xi+\frac{1}{2})|=1
\end{equation}
для почти всех $\xi$.
\Endproc

Последняя теорема суммирует результаты многих математиков: Шура, Харди, Литтвульда, Пойа, Кобера, Михлина и Хермандера.
Мне неизвестно, можно ли в формулировке пункта б) опустить требование неотрицательности ядра, хотя в письме С.~Г.~Самко сообщил мне, что можно. При этом он утверждал, что это следует из результатов одной статьи М.~Г.~Крейна, чего я проверить не смог. Неразобранным остаётся и диапазон $0<p<1$, хотя в общем случае тут оценок нет, на что было мне указано В.~И.~Буренковым на примере операторов Харди.
Первым математиком, использовавшим технику преобразования Меллина для оценки норм операторов Римана--Лиувилля для случая чисто мнимых степеней, был, насколько мне известно, Кобер [237]. Поэтому иногда часть б) приведённой теоремы называется леммой Кобера, что не совсем точно, так как он на самом деле доказал формулу для нормы  из части а) для случая знакопеременных функций.

Отметим, что из теоремы 16 можно сразу вывести симпатичную переформулировку знаменитой гипотезы Римана, которая оказывается эквивалентной обратимости во всех $L_p$ кроме $p=2$ интегральных операторов свёртки с дзета--функцией Римана в качестве мультипликатора. Было бы интересным построить этот оператор свертки в оригиналах (прообразах) в явном виде. По--видимому, этот факт ранее не отмечался, он приятен, но как и тысячи подобных переформулировок совершенно бесполезен в плане хоть какого--то приближения к доказательству знаменитой гипоте
 ы.

Более слабым, чем гипотеза Римана, но также до сих пор недоказанным утверждением является гипотеза Линделёфа, заключающаяся в том, что на критической прямой дзета--функция растёт медленнее любой степени. При исследовании гипотезы Линделёфа существенную роль играют так называемые формулы следа Сельберга--Кузнецова [377], а также дифференциальные уравнения с оператором Бесселя и разложения в ряды особого вида по функциям Бесселя. Поэтому не исключено, что ОП могут найти свои приложения в похожих задачах, когда--то давно автор обсуждал та
 кую возможность с Н.В.~Кузнецовым (чл.--корр. РАН Николай Васильевич Кузнецов скончался в 2010 г., вечная ему память).

\Theorem{ 17} Операторы Бушмана--Эрдейи нулевого порядка гладкости действуют по правилу \eqref{712}. Для их мультипликаторов справедливы формулы:

\begin{eqnarray}
& & m_{{_1S_{0+}^{\nu}}}(s)=\frac{\Gamma(-\frac{s}{2}+\frac{\nu}{2}+1) \Gamma(-\frac{s}{2}-\frac{\nu}{2}+\frac{1}{2})}{\Gamma(\frac{1}{2}-\frac{s}{2})\Gamma(1-\frac{s}{2})}= \label{3.11}; \\
& & =\frac{2^{-s}}{\sq{ \pi}} \frac{\Gamma(-\frac{s}{2}-\frac{\nu}{2}+\frac{1}{2}) \Gamma(-\frac{s}{2}+\frac{\nu}{2}+1)}{\Gamma(1-s)} , Re\, s < \min (2 + Re \, \nu, 1- Re\, \nu) \nonumber ; \\
& & m_{{_1P_{0+}^{\nu}}}(s)=\frac{\Gamma(\frac{1}{2}-\frac{s}{2})\Gamma(1-\frac{s}{2})}{\Gamma(-\frac{s}{2}+\frac{\nu}{2}+1) \Gamma(-\frac{s}{2}-\frac{\nu}{2}+\frac{1}{2})},~ Re\, s < 1; \label{3.12} \\
& & m_{{_1P_{-}^{\nu}}}(s)=\frac{\Gamma(\frac{s}{2}+\frac{\nu}{2}+1) \Gamma(\frac{s}{2}-\frac{\nu}{2})}{\Gamma(\frac{s}{2})\Gamma(\frac{s}{2}+\frac{1}{2})}, Re \, s > \max(Re \, \nu, -1-Re\, \nu); \label{3.13} \\
& & m_{{_1S_{-}^{\nu}}}(s)=\frac{\Gamma(\frac{s}{2})\Gamma(\frac{s}{2}+\frac{1}{2})}{\Gamma(\frac{s}{2}+\frac{\nu}{2}+\frac{1}{2}) \Gamma(\frac{s}{2}-\frac{\nu}{2})}, Re \, s >0 \label{3.14}
\end{eqnarray}

Кроме того, выполняются следующие соотношения для мультипликаторов:
\begin{eqnarray}
& & m_{{_1P_{0+}^{\nu}}}(s)=1/m_{{_1S_{0+}^{\nu}}}(s) ,~ m_{{_1P_{-}^{\nu}}}(s)=1/m_{{_1S_{-}^{\nu}}}(s) \label{3.15} \\
& & m_{{_1P_{-}^{\nu}}}(s)=m_{{_1S_{0+}^{\nu}}}(1-s) ,~ m_{{_1P_{0+}^{\nu}}}(s)=m_{{_1S_{-}^{\nu}}}(1-s) \label{3.16}
\end{eqnarray}
\Endproc

\Theorem{ 18} Справедливы следующие формулы для норм операторов Бушмана--Эрдейи нулевого порядка гладкости в $L_2$:

\begin{eqnarray}
& & \| _1{S_{0+}^{\nu}} \| = \| _1{P_{-}^{\nu}}\|= 1/ \min(1, \sq{1- \sin \pi \nu}), \label{3.22} \\
& & \| _1{P_{0+}^{\nu}}\| = \| _1{S_{-}^{\nu}}\|= \max(1, \sq{1- \sin \pi \nu}). \label{3.23}
\end{eqnarray}
\Endproc

Следствие 1. Нормы операторов \eqref{73} -- \eqref{733} периодичны по $\nu$ с периодом 2, то есть $\|x^{\nu}\|=\|x^{\nu+2}\|$, где $x^{\nu}$ --- любой из операторов \eqref{73} -- \eqref{733}.

Следствие 2. Нормы операторов ${_1 S_{0+}^{\nu}}$, ${_1 P_{-}^{\nu}}$ не ограничены в совокупности по $\nu$, каждая из этих норм не меньше 1. Если $\sin \pi \nu \leq 0$, то эти нормы равны 1. Указанные операторы неограничены в $L_2$ тогда и только тогда, когда $\sin \pi \nu = 1$ (или $\nu=2k + 1/2,~k \in \mathbb{Z}$).

Следствие 3. Нормы операторов ${_1 P_{0+}^{\nu}}$, ${_1 S_{-}^{\nu}}$
ограничены в совокупности по $\nu$, каждая из этих норм не больше $\sq{2}$. Все эти операторы ограничены в $L_2$ при всех $\nu$. Если $\sin \pi \nu \geq 0$, то их $L_2$ -- норма равна 1. Максимальное значение нормы, равное $\sq 2 $, достигается тогда и только тогда, когда $\sin \pi \nu = -1$ (или $\nu=- 1/2+2k ,~k \in \mathbb{Z}$).

Важнейшим свойством операторов Бушмана--Эрдейи нулевого порядка гладкости является их унитарность при целых $\nu$. Отметим, что при интерпретации $L_{\nu}$ как оператора углового момента в квантовой механике, параметр $\nu$ как раз и принимает целые неотрицательные значения.

\Theorem{ 19} Для унитарности в $L_2$ операторов \eqref{73} -- \eqref{733} необходимо и достаточно, чтобы число $\nu$ было целым. В этом случае пары операторов
$({_1 S_{0+}^{\nu}}$, ${_1 P_{-}^{\nu}})$ и  $({_1 S_{-}^{\nu}}$, ${_1 P_{0+}^{\nu}})$
взаимно обратны.
\Endproc

Перед формулировкой частного случая как  следствия предположим, что операторы \eqref{73} -- \eqref{733} заданы на таких функциях $f(x)$, что справедливы представления \eqref{3.1}--\eqref{3.4} (для этого достаточно предположить, что $x f(x) \to 0$ при $x \to 0$). Тогда при $\nu=1$
\be{3.25}{_1{P_{0+}^{1}}f=(I-H_1)f,~_1{S_{-}^{1}}f=(I-H_2)f,}
где $H_1,~ H_2$ -- операторы Харди (см., например, [261--262])
\be{3.26}{H_1 f = \frac{1}{x} \int\limits_0^x f(y) dy,~H_2 f = \int\limits_x^{\infty}  \frac{f(y)}{y} dy,}
$I$ --- единичный оператор.

Следствие 4. Операторы \eqref{3.25} являются унитарными взаимно обратными в $L_2$ операторами. Они сплетают дифференциальные выражения $d^2 / d x^2$ и $d^2 / d x^2 - 2/ x^2$.

Кроме того, можно показать, что операторы \eqref{3.25} являются преобразованиями Кэли  от симметричных операторов $\pm 2 i (x f(x))$ при соответствующем выборе областей определения.

В унитарном случае операторы Бушмана--Эрдейи нулевого порядка  гладкости образуют пару биортогональных преобразований Ватсона, а их ядра образуют пары несимметричных ядер Фурье [207]. Ограниченность операторов с подобными мультипликаторами изучалась ещё Лесли Фоксом.

Отметим важность изучения унитарности для теории интегральных  уравнений. В этом случае обратный оператор необходимо искать в виде интеграла с другими, чем у исходного, пределами интегрирования.

Рассмотрим случай $\nu = i \alpha - \frac{1}{2}$, $\alpha \in \mathbb{R}$, связанный с преобразованием Мелера--Фока.

\Theorem{ 20} Пусть $\nu = i \alpha - \frac{1}{2}$, $\alpha \in \mathbb{R}$. Тогда операторы \eqref{73} -- \eqref{733} ограничены в $L_2$ при всех таких $\nu$. Для их норм справедливы формулы

$$
\| _1{S_{0+}^{i \alpha - \frac{1}{2}}}\|=\| _1{P_{-}^{i \alpha - \frac{1}{2}}}\|=1
$$
\Endproc

Теперь рассмотрим весовые  пространства Лебега с нормой
$$
\|f\|_{L_{2,\, k}}^2=\int\limits_0^{\infty} |f(x)|^2 x^{2k+1} dx,~ k\in \mathbb{R}.
$$
При $k=-1/2$ получаем обычное невесовое пространство $L_2$. В силу того, что ограниченность оператора $A$, действующего из $L_{2, \, k }$ в $L_{2, \, k }$, эквивалентна ограниченности оператора $B=x^{k+\frac{1}{2}}Ax^{-(k+\frac{1}{2})}$ в невесовых пространствах $L_2$, все утверждения леммы 3.3 сохраняются и для весового случая с заменой величины $i \xi + \frac{1}{2}$ на $i \xi + k + 1$.

\Theorem{ 21} a) Оператор ${_1S_{-}^{\nu}}$ неограничен в $L_{2, \, k }$
при условии $k=-n$, $n \in \mathbb{N}$.

б) Оператор ${_1P_{0+}^{\nu}}$ неограничен в $L_{2, \, k }$
при условии $k=n$, $n \in \mathbb{N}$.

в) Оператор ${_1P_{-}^{\nu}}$ неограничен в $L_{2, \, k }$
при условии $k+\nu=-2 m_1-2$, $m_1 \in \mathbb{N}$ или $k-\nu=-2 m_2-1$, $m_2 \in \mathbb{N}_0$.

г) Оператор ${_1S_{0+}^{\nu}}$ неограничен в $L_{2, \, k }$
при условии $k+\nu=m$, $m \in \mathbb{N}_0$ или $k-\nu=n$, $n \in \mathbb{N}$.

д) При всех остальных значениях $\nu$, $k \in \mathbb{R}$ эти операторы ограничены
в $L_{2, \, k }$.
\Endproc

Покажем, что свойство унитарности операторов Бушмана--Эрдейи нулевого порядка гладкости присуще лишь случаю $k=-1/2$ невесовых пространств $L_2$.

\Theorem{ 22}Рассмотрим операторы ${_1 S_-^1}$ и ${_1 S_-^{-1/2}}$. Если оператор
${_1 S_-^1}$ унитарен в $L_{2, \, k}$, то $k=-1/2$. Если при $k > -1$ оператор ${_1 S_-^{-1/2}}$ унитарен в $L_{2, \, k}$, то также $k=-1/2$.
\Endproc

Оказывается, что отсутствие унитарности в весовом случае можно доказывать для всех $\nu$ из интервала $(-1, 0)$. На концах этого интервала все рассматриваемые операторы унитарны, так как превращаются в тождественные.

\Theorem{ 23} Оператор ${_1 S_{-}^{\nu}}$ не является изометрией в $L_{2, \, k}$, если выполнено хотя бы одно из следующих условий:

а) $-1<\nu<0$, $k>-1$;

б) $\nu$ --- комплексное число, $k>-1$, число $z=\nu(\nu+1)/(k+1)(k+2)$ лежит в круге с центром в точке $z_0=1$ единичного радиуса;

в) $\nu > 0$, $2m <\nu \leq 2m +(\sq7-1)/2$, $m \in \mathbb{N}_0$, $k > - \frac{1}{2}$.
\Endproc

Предыдущая теорема оставляет "зазор" \ в значениях параметра $\nu$, для которых отсутствие изометричности пока не установлено.

\Theorem{ 24} Пусть $k>-1$, $\nu \neq 0$, $\nu \neq -1$. Тогда из унитарности в $L_{2, k}$ любого из операторов Бушмана--Эрдейи нулевого порядка гладкости следует, что
$k=-1/2$, $\nu$ --- целое число.
\Endproc

Перейдём к вопросу о спектре рассматриваемых нами операторов. Ясно, что этот спектр может быть только непрерывным. Из теоремы 17 следует, что в $L_2$ спектр совпадает с замыканием образа прямой $Re\,s=1/2$ при отображении $m(s)$. Кривую, описывающую спектр в плоскости спектрального параметра, будем задавать функцией $W(z)$.
Рассмотрим операторы \eqref{3.25} при $\nu=1$, связанные с операторами Харди \eqref{3.26}.

\Theorem{ 25} Рассмотрим оператор
$$
_1{S_-^1} f=(I-H_2)f=f(x)-\int\limits_x^{\infty} \frac{f(y)}{y} dy.
$$
Тогда если, $k \neq 0$, то
а) оператор $_1{S_-^1}$ ограничен в $L_{2, \, k}$, и его норма равна
\be{4.3}{\|{_1S_-^1}\|_{L_{2,\,k}}= \left\{ \begin{matrix} \left|\frac{k+1}{k} \right|, & k \geq - \frac{1}{2},\\ 1, & k \leq - \frac{1}{2}, \end{matrix} \right.}
б) спектр оператора $_1{S_-^1}$ в $L_{2, \, k}$ есть окружность, проходящая через точку $z=1$. Её центр лежит на вещественной оси в точке $w_0$,  радиус равен $\rho$, где
\be{4.4}{w_0=1+\frac{1}{2k},~ \rho = \left|\frac{1}{2k}\right|}
в) в случае $k=0$ оператор неограничен в $L_{2, \, k}$, а его спектр есть прямая $Re \, w=1$.
\Endproc

\Theorem{ 26} Рассмотрим оператор
$$
_1{P_{0+}^1} f=(I-H_1)f=f(x)-  \frac{1}{x} \int\limits_0^x f(y)\,dy.
$$
Тогда если, $k \neq -1$, то

а) оператор $_{1 P_{0+}^1}$ ограничен в $L_{2, \, k}$, и его норма равна
\be{4.5}{\|{_1P_{0+}^1}\|_{L_{2,\,k}}= \left\{ \begin{matrix} 1, & k \leq - \frac{1}{2},\\  \left|\frac{k+1}{k} \right| & k \geq - \frac{1}{2}; \end{matrix} \right.}
б) спектр оператора $_1{P_{0+}^1}$ в $L_{2, \, k}$ есть окружность, проходящая через точку $z=1$. Её центр лежит на вещественной оси в точке $w_0$,  радиус равен $\rho$, где
\be{4.6}{w_0=1-\frac{1}{2k+2},~ \rho = |\frac{1}{2k+2}|\ ;}
в) в случае $k=-1$ оператор неограничен в $L_{2, \, k}$, а его спектр есть прямая $Re \, w=1$.
\Endproc

Далее перечислим некоторые общие свойства операторов, которые  действуют по правилу \eqref{712} как умножение на некоторый мультипликатор в образах преобразования  Меллина и одновременно являются сплетающими для второй производной и оператора углового момента.

\Theorem{ 27} Пусть оператор $S_{\nu}$ действует по формулам \eqref{712}
и \eqref{76}. Тогда

а) его мультипликатор удовлетворяет функциональному уравнению

\be{5.1}{m(s)=m(s-2)\frac{(s-1)(s-2)}{(s-1)(s-2)-\nu(\nu+1)};}

б) если функция $p(s)$ периодична с периодом 2 (то есть $p(s)=p(s-2)$),  то функция $p(s)m(s)$ является мультипликатором нового оператора  преобразования $S_2^{\nu}$, опять же сплетающего $L_{\nu}$ и вторую производную по правилу
\eqref{76}.
\Endproc

Последняя теорема ещё раз показывает, насколько удобно изучение ОП в терминах мультипликаторов преобразования Меллина.

Определим преобразование Стилтьеса (см., например, [34]) по формуле
$$
(S f)(x)= \int\limits_0^{\infty} \frac{f(t)}{x+t} dt.
$$
Этот оператор имеет вид \eqref{712} с мультипликатором  $p(s)= \pi /sin (\pi s)$ и ограничен в $L_2$. Очевидно, что $p(s)=p(s-2)$. Поэтому из теоремы 27 следует, что композиция преобразования Стилтьеса с ограниченными сплетающими операторами \eqref{73}--\eqref{733} снова является оператором преобразованием
того же типа, ограниченным в $L_2$.

Отметим, что из предыдущего изложения  следует, что
$$
\|S\|_{L_p}=|\pi / \sin \frac{\pi}{p}|,~p>1.
$$
С другой стороны,
$$
\|S\|_{L_{2,\, k}}=|\pi / \sin \pi k |,~k \notin  \mathbb{Z}.
$$
Аналогично получаются оценки в весовых пространствах $L_{p, \, k} ~ p>0 $.

Рассмотрим теперь оператор $H^{\nu}$ вида \eqref{712} с мультипликатором
\begin{equation}
\label{009}
m(s)=\sq{\frac{\sin \pi s - \sin \pi \nu}{\sin \pi s}}
\end{equation}
Из теоремы 17 получаем, что на прямой $Re \, s = \frac{1}{2}$ величина $m(s)$ обратна величине \eqref{009}. Тогда из теоремы 18 следует, что
\be{010}{\|H^{\nu}\|_{L_2}=\|_1{P_{0+}^{\nu}}\|_{L_2}=\|_1{S_{-}^{\nu}}\|_{L_2}.}
Поэтому для оператора $H^{\nu}$ справедливо следствие 3 теоремы 18.  В частности $H^{\nu}$ ограничен в $L_2$ при всех $\nu$.

Отметим, что формально этот оператор связан с преобразованием Стилтьеса формулой
\be{5.3}{H^{\nu}=(1-\frac{\sin \pi \nu}{\pi}S)^{\frac{1}{2}}.}
Одновременно с $H^{\nu}$ введём оператор $\mathfrak{D}^{\nu}$ с мультипликатором
$$
m_{\mathfrak{D}^{\nu}}(s)=\sq{\frac{\sin \pi s}{\sin \pi s - \sin \pi \nu}}
$$
Отсюда получаем, что
\be{5.4}{\|\mathfrak{D}^{\nu}\|_{L_2}=\|_1{S_{0+}^{\nu}}\|_{L_2}=\|_1{P_{-}^{\nu}}\|_{L_2}.}
и кроме того, оператор $\mathfrak{D}^{\nu}$ ограничен при $\sin  \pi \nu \neq 1$.

\Theorem{ 28} Рассмотрим композиции операторов
\begin{eqnarray}
& &  _3{S^{\nu}_{0+}}={_1S^{\nu}_{0+}}H^{\nu},~ _3{S^{\nu}_{-}}=\mathfrak{D}^{\nu} {_1S^{\nu}_{-}}, \label{5.5} \\
& &  _3{P^{\nu}_{0+}}=\mathfrak{D}^{\nu} {_1P^{\nu}_{0+}},~ _3{P^{\nu}_{-}}={_1P^{\nu}_{-}}H^{\nu}. \label{5.6}
\end{eqnarray}

Тогда операторы ${_3 S_{0+}^{\nu}}$, ${_3 S_{-}^{\nu}}$ являются новыми операторами  преобразования типа Сонина, а ${_3 P_{0+}^{\nu}}$, ${_3 P_{-}^{\nu}}$ --- типа Пуассона. Все эти операторы унитарны в $L_2$.
Кроме того, если $\sin  \pi \nu \neq 1$, то композиции \eqref{5.5} -- \eqref{5.6} можно вычислять в любом порядке.
\Endproc

Ниже будет получено явное интегральное представление операторов  преобразования, сплетающих $L_{\nu}$ и
$d^2/dx^2$, которые \textit{являются унитарными при всех $\nu \in \mathbb{R}$}.

Изучим вопрос о взаимосвязи разносторонних операторов Бушмана--Эрдейи. Полученные формулы аналогичны тем, которые связывают лево-- и правосторонние дробные интегралы Римана--Лиувилля (см. [34], с.~163--171). С этой целью введём оператор
\be{5.7}{C^{\nu} f= f(x)-\frac{sin \pi \nu}{\pi} S f,}
где $S$ --- преобразование Стильтьеса.
Приведём без доказательства свойства $C^{\nu}$:

1) $\|C^{\nu}\|_{L_2}=\min (1, 1 - \sin \pi \nu) \leq 1, ~ \nu \in \mathbb{R}$;

2) $\|C^{\nu}\|_{L_2}=1+\ch \pi \alpha, ~  \nu=i \alpha - \frac{1}{2},~ \alpha \in \mathbb{R}$;

3) при всех $\nu,~k$ можно вычислить весовую норму $C^{\nu}$. Например, если $\sin  \pi \nu \cdot \sin  \pi \nu > 0$, то
$$
\|C^{\nu}\|_{L_{2, \, k}}=1+\frac{\sin \pi \nu}{\sin \pi k}.
$$

\Theorem{ 29} При $\nu \in \mathbb{R}$ справедливы тождества для композиций
\begin{eqnarray}
& &  C^{\nu}={_1S_{-}^{\nu}} \ {_1P_{0+}^{\nu}}={_1P_{0+}^{\nu}} \  {_1S_{-}^{\nu}} , \label{5.8} \\
& &  {_1S_{-}^{\nu}}={_1S_{0+}^{\nu}} \ C^{\nu}, ~   {_1P_{0+}^{\nu}}= {_1P_{-}^{\nu}} \ C^{\nu} , \label{5.9} \\
& &  {_1S_{-}^{\nu}}=C^{\nu} \  {_1S_{0+}^{\nu}}, ~   {_1P_{0+}^{\nu}}=C^{\nu} \  {_1P_{-}^{\nu}} , ~ \sin \pi \nu \neq 1. \label{5.10}
\end{eqnarray}
\Endproc

Теперь определим и изучим операторы Бушмана--Эрдейи второго рода.

Введём новую пару операторов Бушмана--Эрдейи с функциями Лежандра второго рода [102] в ядре
\be{6.1}{{_2S^{\nu}}f=\frac{2}{\pi} \left( - \int\limits_0^x (x^2-y^2)^{-\frac{1}{2}}Q_{\nu}^1 (\frac{x}{y}) f(y) dy  +
\int\limits_x^{\infty} (y^2-x^2)^{-\frac{1}{2}}\mathbb{Q}_{\nu}^1 (\frac{x}{y}) f(y) dy\right),}
\be{6.2}{{_2P^{\nu}}f=\frac{2}{\pi} \left( - \int\limits_0^x (x^2-y^2)^{-\frac{1}{2}}\mathbb{Q}_{\nu}^1 (\frac{y}{x}) f(y) dy  -
\int\limits_x^{\infty} (y^2-x^2)^{-\frac{1}{2}}Q_{\nu}^1 (\frac{y}{x}) f(y) dy\right).}

При $y \to x \pm 0$ интегралы понимаются в смысле главного значения. Отметим без доказательства, что эти операторы определены и являются сплетающими при некоторых условиях на функции $f(x)$  (при этом оператор \eqref{6.1} будет типа Сонина, \eqref{6.2}  --- типа Пуассона).

\Theorem{ 30}  Операторы  \eqref{6.1} -- \eqref{6.2} представимы в виде \eqref{712}
с мультипликаторами
\begin{eqnarray}
& & m_{_2S^{\nu}}(s)=p(s) \ m_{_1S_{-}^{\nu}}(s), \label{6.3}\\
& & m_{_2P^{\nu}}(s)=\frac{1}{p(s)} \ m_{_1P_{-}^{\nu}}(s), \label{6.4}
\end{eqnarray}
где мультипликаторы операторов ${_1S_-^{\nu}}$, ${_1P_-^{\nu}}$ определены формулами \eqref{3.13} -- \eqref{3.14}, а функция $p(s)$ ( с периодом 2) равна

\be{6.5}{p(s)=\frac{\sin \pi \nu+ \cos \pi s}{\sin \pi \nu - \sin \pi s}.}
\Endproc

Лемма 1. Рассмотрим более общий чем \eqref{6.1} интегральный оператор при значениях $Re\, \nu < 1$:
\begin{eqnarray}
& & {_1S^{\nu}}f=\frac{2}{\pi} \left( \int\limits_0^x (x^2+y^2)^{-\frac{\mu}{2}} e^{-\mu \pi i} Q_{\nu}^{\mu}( \frac{x}{y}) f(y)\, dy\right. + \label{6.6} \\
& & +\left. \int\limits_x^{\infty} (y^2+x^2)^{-\frac{\mu}{2}}\mathbb{Q}_{\nu}^{\mu} (\frac{x}{y}) f(y)\, dy\right) \nonumber,
\end{eqnarray}
где $Q_{\nu}^{\mu}(z)$ --- функция Лежандра второго рода, $\mathbb{Q}_{\nu}^{\mu}(z)$ --- значение этой функции на разрезе.

Тогда на функциях из $C_0^{\infty}(0, \infty)$ оператор \eqref{6.6} определён и действует по формуле
$$
M[{_3S^{\nu}}](s)=m(s)\cdot M[x^{1-\mu} f](s), \label{6.7}
$$
\begin{eqnarray}
 m(s)=2^{\mu-1} \left( \frac{ \cos \pi(\mu-s) - \cos \pi \nu}{ \sin \pi(\mu-s) - \sin \pi \nu}  \right) \cdot\\
\cdot \left( \frac{\Gamma(\frac{s}{2})\Gamma(\frac{s}{2}+\frac{1}{2}))}{\Gamma(\frac{s}{2}+\frac{1-\nu-\mu}{2}) \Gamma(\frac{s}{2}+1+\frac{\nu-\mu}{2})} \right). \nonumber
\end{eqnarray}

\Theorem{ 31} Справедливы формулы для норм
\begin{eqnarray}
& &  \| {_2S^{\nu}} \|_{L_2}= \max (1, \sq{1+\sin \pi \nu}) , \label{6.9} \\
& &  \| {_2P^{\nu}} \|_{L_2}= 1 / {\min (1, \sq{1+\sin \pi \nu})} . \label{6.10}
\end{eqnarray}
\Endproc

Следствие.  Оператор ${_2S^{\nu}}$ ограничен при всех $\nu$. Оператор ${_2P^{\nu}}$ не является непрерывным тогда и только тогда, когда
$\sin \pi \nu=-1$.

\Theorem{ 32} Для унитарности в $L_2$ операторов ${_2S^{\nu}}$ и ${_2P^{\nu}}$ необходимо и достаточно, чтобы параметр $\nu$ был целым
числом.
\Endproc

\Theorem{ 33} Пусть $\nu=i \beta+1/2,~\beta \in \mathbb{R}$. Тогда
\be{6.11}{\| {_2S^{\nu}} \|_{L_2}=\sq{1+\ch \pi \beta},~\| {_2P^{\nu}} \|_{L_2}=1.}
\Endproc

\Theorem{ 34} Справедливы представления
\begin{eqnarray}
& & {_2S^0} f = \frac{2}{\pi} \int\limits_0^{\infty} \frac{y}{x^2-y^2}f(y)\,dy, \label{6.12} \\
& & {_2S^{-1}} f = \frac{2}{\pi} \int\limits_0^{\infty} \frac{x}{x^2-y^2}f(y)\,dy. \label{6.13}
\end{eqnarray}
\Endproc

Таким образом, в этом случае оператор ${_2S^{\nu}}$ сводится к паре известных преобразований Гильберта на полуоси [34].

Перейдём к построению операторов преобразования, унитарных при  всех $\nu$. Такие операторы определяются по формулам:
\begin{eqnarray}
& & S_U^{\nu} f = - \sin \frac{\pi \nu}{2} {_2S^{\nu}}f+ \cos \frac{\pi \nu}{2} {_1S_-^{\nu}}f, \label{6.14} \\
& & P_U^{\nu} f = - \sin \frac{\pi \nu}{2} {_2P^{\nu}}f+ \cos \frac{\pi \nu}{2} {_1P_-^{\nu}}f. \label{6.15}
\end{eqnarray}
Для любых значений $\nu \in \mathbb{R}$ они являются линейными комбинациями операторов преобразования Бушмана--Эрдейи 1 и 2 рода нулевого порядка гладкости. Их можно назвать операторами Бушмана--Эрдейи третьего рода. В интегральной форме эти операторы имеют вид:

\begin{eqnarray}
& & S_U^{\nu} f = \cos \frac{\pi \nu}{2} \left(- \frac{d}{dx} \right) \int\limits_x^{\infty} P_{\nu}(\frac{x}{y}) f(y)\,dy +  \label{6.16}\\
& & + \frac{2}{\pi} \sin \frac{\pi \nu}{2} \left(  \int\limits_0^x (x^2-y^2)^{-\frac{1}{2}}Q_{\nu}^1 (\frac{x}{y}) f(y)\,dy  -
 \int\limits_x^{\infty} (y^2-x^2)^{-\frac{1}{2}}\mathbb{Q}_{\nu}^1 (\frac{x}{y}) f(y)\,dy\right), \nonumber \\
& & P_U^{\nu} f = \cos \frac{\pi \nu}{2}  \int\limits_0^{x} P_{\nu}(\frac{y}{x}) \left( \frac{d}{dy} \right) f(y)\,dy - \label{6.17} \\
& &  -\frac{2}{\pi} \sin \frac{\pi \nu}{2} \left( - \int\limits_0^x (x^2-y^2)^{-\frac{1}{2}}\mathbb{Q}_{\nu}^1 (\frac{y}{x}) f(y)\,dy  -
\int\limits_x^{\infty} (y^2-x^2)^{-\frac{1}{2}} Q_{\nu}^1 (\frac{y}{x}) f(y)\,dy \right). \nonumber
\end{eqnarray}

\Theorem{ 35} Операторы \eqref{6.14}--\eqref{6.15} при всех $\nu$ являются унитарными, взаимно сопряжёнными и обратными в $L_2$. Они являются
сплетающими и действуют по формулам \eqref{75}. При этом $S_U^{\nu}$ является оператором типа Сонина, а $P_U^{\nu}$ --- типа Пуассона.
\Endproc

\newpage
Далее  приведены некоторые приложения и обобщения полученных результатов.

А. Оценки полунорм в функциональных пространствах.

Рассмотрим множество функций $\mathfrak{D}(0, \infty)$. Если $f(x) \in \mathfrak{D}(0, \infty)$, то $f(x) \in C^{\infty}(0, \infty),~ f(x)$ ---финитна на бесконечности. На этом множестве функций введём полунормы
\begin{eqnarray}
& & \|f\|_{h_2^{\alpha}}=\|\mathfrak{D}_-^{\alpha}f\|_{L_2(0, \infty)} \label{7.1} \\
& &  \|f\|_{\widehat{h}_2^{\alpha}}=\|x^{\alpha} (-\frac{1}{x}\frac{d}{dx})^{\alpha}f\|_{L_2(0, \infty)} \label{7.2}
\end{eqnarray}
где $\mathfrak{D}_-^{\alpha}$ --- дробная производная Римана--Лиувилля, оператор в \eqref{7.2} определяется по формуле
\be{7.3}{(-\frac{1}{x}\frac{d}{dx})^{\beta}=2^{\beta}I_{-; \, 2, \,0}^{-\beta}x^{-2 \beta},}
$I_{-; 2, \, 0}^{-\beta}$ --- оператор Эрдейи--Кобера, $\alpha$ --- произвольное действительное число. При $\beta = n \in \mathbb{N}_0$ выражение \eqref{7.3} понимается в обычном смысле, что согласуется с определением (7.2).

Лемма 1.  Пусть $f(x) \in \mathfrak{D}(0, \infty)$. Тогда справедливы тождества:
\begin{eqnarray}
& & \mathfrak{D}_-^{\alpha}f={_1S_-^{\alpha-1}} [x^{\alpha} (-\frac{1}{x}\frac{d}{dx})^{\alpha}] f, \label{7.4} \\
& & x^{\alpha} (-\frac{1}{x}\frac{d}{dx})^{\alpha}f={_1P_-^{\alpha-1}} \mathfrak{D}_-^{\alpha}f. \label{7.5}
\end{eqnarray}

Таким образом, операторы Бушмана--Эрдейи нулевого порядка гладкости первого рода осуществляют связь между дифференциальными операторами (при $\alpha \in \mathbb{N}$) из определений полунорм \eqref{7.1} и \eqref{7.2}.

Лемма 2.  Пусть $f(x) \in \mathfrak{D}(0, \infty)$. Тогда справедливы неравенства между полунормами
\begin{eqnarray}
& &  \|f\|_{h_2^{\alpha}} \leq \max (1, \sq{1+\sin \pi \alpha}) \|f\|_{\widehat{h}_2^{\alpha}}, \label{7.7}\\
& & \|f\|_{\widehat{h}_2^{\alpha}} \leq \frac{1}{\min (1, \sq{1+\sin \pi \alpha})} \|f\|_{h_2^{\alpha}}, \label{7.8}
\end{eqnarray}
где $\alpha$ --- любое действительное число, $\alpha \neq -\frac{1}{2}+2k,~k \in \mathbb{Z}$.

Постоянные в неравенствах \eqref{7.7}--\eqref{7.8} не меньше единицы, что будет далее использовано. В случае
$\sin \pi \alpha = -1 $ или $\alpha = -\frac{1}{2}+2k,~k \in \mathbb{Z}$, оценка \eqref{7.8} не имеет места.

Введём на $\mathfrak{D} (0, \infty )$ соболевскую норму
\be{7.9}{\|f\|_{W_2^{\alpha}}=\|f\|_{L_2 (0, \infty)}+\|f\|_{h_2^{\alpha}}.}
Введём также другую норму
\be{7.10}{\|f\|_{\widehat{W}_2^{\alpha}}=\|f\|_{L_2 (0, \infty)}+\|f\|_{\widehat{h}_2^{\alpha}}}
Пространства $W_2^{\alpha},~ \widehat{W}_2^{\alpha}$ определим как замыкания $\mathfrak{D}(0,
\infty)$ по нормам \eqref{7.9} и \eqref{7.10} соответственно.

\Theorem{ 36} а) при всех $\alpha \in \mathbb{R}$ пространство $\widehat{W}_2^{\alpha}$ непрерывно вложено в $W_2^{\alpha}$, причём
\be{7.11}{\|f\|_{W_2^{\alpha}}\leq A_1 \|f\|_{\widehat{W}_2^{\alpha}},}
где $A_1=\max (1, \sq{1+\sin \pi \alpha})$.

б) Пусть $\sin \pi \alpha \neq -1$ или $\alpha \neq -\frac{1}{2} + 2k, ~ k \in \mathbb{Z}.  $ Тогда справедливо обратное вложение $W_2^{\alpha}$  в $\widehat{W}_2^{\alpha}$, причём
\be{7.12}{\|f\|_{\widehat{W}_2^{\alpha}}\leq A_2 \|f\|_{W_2^{\alpha}},}
где $A_2 =1/  \min (1, \sq{1+\sin \pi \alpha})$.

в) Пусть $\sin \pi \alpha \neq -1$, тогда пространства $W_2^{\alpha}$  и $\widehat{W}_2^{\alpha}$ изоморфны, а их нормы эквивалентны.

г) Константы в неравенствах вложений \eqref{7.11}--\eqref{7.12} точные.
\Endproc

Эта теорема фактически является следствием результатов по ограниченности операторов Бушмана--Эрдейи нулевого порядка гладкости в $L_2$. В свою очередь, из теоремы об унитарности этих операторов вытекает

\Theorem{ 37} Нормы
\begin{eqnarray}
& & \|f\|_{W_2^{\alpha}} = \sum\limits_{j=0}^s \| \mathfrak{D}_-^j f\|_{L_2}, \label{7.13} \\
& & \|f\|_{\widehat{W}_2^{\alpha}}=\sum\limits_{j=0}^s \| x^j(-\frac{1}{x}\frac{d}{dx})^j f \|_{L_2} \label{7.14}
\end{eqnarray}
задают эквивалетные нормировки в пространстве Соболева при целых $s \in \mathbb{Z}$. Кроме того, каждое слагаемое в \eqref{7.13} тождественно равно соответствующему слагаемому в \eqref{7.14} с тем же индексом $j$.
\Endproc

И. А. Киприянов ввёл в [175] шкалу пространств, которые оказали существенное влияние на теорию уравнений в частных производных с оператором Бесселя по одной или нескольким переменным. Эти пространства можно определить следующим образом. Рассмотрим подмножество чётных функций на $\mathfrak{D}(0, \infty)$, у которых все производные нечётного порядка равны нулю при $x=0$. Обозначим это множество $\mathfrak{D}_c (0, \infty)$ и введём на нём норму
\be{7.15}{\|f\|_{\widetilde{W}_{2, k}^s} = \|f\|_{L_{2, k}}+\|B_k^{\frac{s}{2}}\|_{L_{2, k}}}
где $s$ --- чётное натуральное число,  $B^{s/2}_k$ --- итерация оператора Бесселя. Пространство И.А.~Киприянова при чётных $s$ определяется как замыкание $\mathfrak{D}_c (0, \infty)$ по норме \eqref{7.15}. Известно, что эквивалентная \eqref{7.15} норма может быть задана по формуле [175]
\be{7.16}{\|f\|_{\widetilde{W}_{2, k}^s} = \|f\|_{L_{2, k}}+\|x^s(-\frac{1}{x}\frac{d}{dx})^s f\|_{L_{2, k}}}
Это позволяет доопределить норму в $\widetilde{W}_{2, \, k}^s$ для всех $s$. Отметим, что по существу этот подход совпадает с одним из принятых в [175], другой подход основан на использовании преобразования Фурье--Бесселя. Далее будем считать, что $\widetilde{W}_{2, k}^s$ нормируется по формуле \eqref{7.16}.

Введём весовую соболевскую норму
\be{7.17}{\|f\|_{W_{2, k}^s} = \|f\|_{L_{2, k}}+\|\mathfrak{D}_-^s f\|_{L_{2, k}}}
и пространство $W_{2, \, k}^s$ как замыкание $\mathfrak{D}_c (0, \infty)$ по этой норме.

\Theorem{ 38} а) Пусть $k \neq -n, ~ n \in \mathbb{N}$. Тогда пространство  $\widetilde{W}_{2, \, k}^s$ непрерывно вложено в $W_{2, \, k}^s$, причём существует постоянная $A_3>0$ такая, что
\be{7.18}{\|f\|_{W_{2, k}^s}\leq A_3 \|f\|_{\widetilde{W}_{2, k}^s},}
б) Пусть $k+s \neq -2m_1-1, ~ k-s \neq -2m_2-2, ~ m_1 \in \mathbb{N}_0, ~ m_2 \in \mathbb{N}_0$. Тогда справедливо обратное вложение $W_{2, \, k}^s$ в $\widetilde{W}_{2, \, k}^s$, причём существует постоянная $A_4>0$, такая, что
\be{7.19}{\|f\|_{\widetilde{W}_{2, k}^s}\leq A_4 \|f\|_{W_{2, k}^s}.}
в) Если указанные условия не выполняются, то соответствующие вложения не имеют места.
\Endproc

Следствие 1.  Пусть выполнены условия: $k \neq -n, ~ n \in \mathbb{N}$; $k+s \neq -2m_1-1,  ~ m_1 \in \mathbb{N}_0; ~ k-s \neq -2m_2-2, ~ m_2 \in \mathbb{N}_0$. Тогда пространства И.А. Киприянова можно определить как замыкание $\mathfrak{D}_c (0, \infty)$ по весовой соболевской норме \eqref{7.17}.

Следствие 2.  Точные значения постоянных в неравенствах вложения \eqref{7.18}--\eqref{7.19} есть
$$
A_3 = \max (1, \|{_1S_-^{s-1}} \| _ {L_{2, k}}), ~ A_4=\max(1, \|{_1P_-^{s-1}}\|_{L_{2, k}}).
$$

Очевидно, что приведённая теорема и следствия из неё вытекают из приведённых выше результатов для операторов Бушмана--Эрдейи. Отметим, что  нормы операторов Бушмана--Эрдейи нулевого порядка гладкости в $L_{2, \, k}$  дают значения точных постоянных в неравенствах вложения \eqref{7.18}--\eqref{7.19}. Оценки норм операторов Бушмана--Эрдейи в банаховых пространтсвах $L_{p, \alpha}$ позволяют рассматривать вложения для соответствующих функциональных пространств.

Неравенство для полунорм, приводящее к вложению \eqref{7.18} ($s$ --- целое число), получено в [320]. Вложения типа полученных в теореме 38 ранее изучались в [321--322]. В последней работе рассматривался случай $k>-1/2, ~ s \in \mathbb{N}$, пространства $W_{p,\,k}^s$. Мы рассматриваем гильбертовы пространства $W_{2,\,k}^s$, $k$ и $s$ --- любые действительные числа. Кроме того, в теореме 38 уточнены условия отсутствия вложений из [321--322], которые содержали ошибки (см. ниже).   Отметим, что в теореме 38 фактически установлены более точные неравенства между соответствующими полунормами, чем
 в предшествующих работах [320--322]. Это стало возможным благодаря применению подробно изученных выше ОП Бушмана--Эрдейи.

Перейдём к рассмотрению правосторонних сплетающих операторов \eqref{73}--\eqref{733}. Мы покажем, что в общем случае они осуществлют изоморфизм пространства С.Л. Соболева и И.А. Киприянова. Определим пространства $H^{2s}$, $H_{\alpha}^{2s}$ и $K_{\alpha}^{2s}$ как замыкания множества функций $\mathfrak{D} (0, \infty)$ по нормам
\begin{eqnarray}
& & \|f\|_{H^{2s}} = \|f\|_{L_2}+\|\mathfrak{D}_-^{2s} f\|_{L_2}, \label{7.20} \\
& & \|f\|_{H^{2s}_{\alpha}} = \|f\|_{L_{2, {\alpha}}}+\|\mathfrak{D}_-^{2s} f\|_{L_{2, {\alpha}}}, \label{7.21} \\
& & \|f\|_{K^{2s}_{\alpha}} = \|f\|_{L_{2, {\alpha}}}+\|B_{\alpha}^s f\|_{L_{2, {\alpha}}}, \label{7.22}
\end{eqnarray}
$s$ --- натуральное число, $\alpha \in \mathbb{R}$. Определим также пару операторов типа \eqref{76}
\be{7.23}{{_1X_-^{\alpha}}={_1S_-^{\alpha-\frac{1}{2}}} x^{\alpha+\frac{1}{2}}, ~ {_1Y_-^{\alpha}}= x^{-(\alpha+\frac{1}{2})} {_1P_-^{\alpha-\frac{1}{2}}}. }

\Theorem{ 39} Пусть $\alpha \in \mathbb{R},~ s \in \mathbb{N}$. Тогда оператор ${_1X^2_-}$ действует непрерывно из $H^{2s}$ в $K^{2s}$, причём
\be{7.24}{\|{_1X_-^{\alpha}}f\|_{H^{2s}} \leq A_5 \|f\|_{K^{2s}} ,}
где $A_5=\|{_1X_-^{\alpha}}\|_{H^{2s} \to K^{2s}}=\|{_1S_-^{\alpha-\frac{1}{2}}}\|_{L_2}= \max (1, \, \sq{1+\cos \pi \alpha})$

Пусть $s \in \mathbb{N}, \alpha \neq 2k+1, ~ k \in \mathbb{Z}$ (или $\cos \pi \alpha \neq -1$). Тогда оператор ${_1Y_-^{\alpha}}$ действует непрерывно из в $K^{2s}$ в $H^{2s}$, причём справедливо неравенство
\be{7.24}{\|{_1Y_-^{\alpha}}f\|_{K^{2s}_{\alpha}} \leq A_6 \|f\|_{H^{2s}} ,}
с постоянной
$$
A_6=\|{_1Y_-^{\alpha}}\|_{K^{2s} \to H^{2s}}=\|{_1P_-^{\alpha-\frac{1}{2}}}\|_{L_2}=1/ \max (1, \, \sq{1+\cos \pi \alpha})
$$
\Endproc

Оператор Бесселя является радиальной частью лапласиана в $\mathbb{R}^n$. При такой интерпретации этого оператора в случае нечётномерных пространств  будет выполнено условие теоремы 39.

\Theorem{ 40} Пусть выполнены условия $\alpha \neq 2k+1, ~ k \in \mathbb{Z}$, $\alpha \neq -n, ~ n \in \mathbb{N}$; $\alpha+2s \neq -2m_1-1,  ~ m_1 \in \mathbb{N}_0; ~ \alpha-2s \neq -2m_2-2, ~ m_2 \in \mathbb{N}_0$. Тогда операторы \eqref{7.23} осуществляют топологический изоморфизм пространств Соболева $H^{2s}$ и весового пространства Соболева $H^{2s}_{\alpha}$.
\Endproc

Очевидно, что все условия теоремы 40 выполнены при полуцелых $\alpha \in \mathbb{R}$. Поэтому справедлива

\Theorem{ 41} Пусть $s \in \mathbb{N},~\alpha - \frac{1}{2} \in \mathbb{Z}$. Тогда операторы \eqref{7.23} осуществляют топологический изоморфизм пространств Соболева $H^{2s}$ и весовых пространств Соболева $H^{2s}_{\alpha}$.
\Endproc

Аналогично можно ввести по формулам типа \eqref{7.23} операторы ${_1X^{\alpha}}_{0+}$ и ${_1Y^{\alpha}}_{0+}$. В качестве приложения приведённых выше результатов можно также рассмотреть действие операторов \eqref{7.23} в пространствах с нормами \eqref{7.20}--\eqref{7.22} при произвольных весах, не согласованных с постоянной $\alpha$ в операторе Бесселя $B_{\alpha}$.

Принципиальная важность пространств И.А.\,Киприянова для теории уравнений в частных производных различных типов с операторами Бесселя отражает общий методологический подход, который автор услышал в виде красивого афоризма на пленарной лекции чл.--корр. РАН Л.Д.\,Кудрявцева: "\textit{КАЖДОЕ УРАВНЕНИЕ ДОЛЖНО ИЗУЧАТЬСЯ В СВОЁМ СОБСТВЕННОМ ПРОСТРАНСТВЕ!}"

\vskip 0.5cm

Б. Задача Коши для уравнения Эйлера--Пуассона--Дарбу (ЭПД).

Рассмотрим уравнение ЭПД в полупространстве
$$
B_{\alpha,\, t} u(t,x)= \frac{ \pr^2 \, u }{\pr t^2} + \frac{2 \alpha+1}{t} \frac{\pr u}{\pr t}=\Delta_x u+F(t, x),
$$
где $t>0,~x \in \mathbb{R}^n$. Дадим нестрогое описание процедуры, позволяющей получать различные постановки начальных условий при $t=0$ единым методом. Образуем по формулам \eqref{75} операторы преобразования $X_{\alpha, \, t}$ и $Y_{\alpha, \, t}$. Предположим, что существуют выражения $X_{\alpha, \, t} u=v(t,x)$, $X_{\alpha, \, t} F=G(t,x)$. Пусть обычная (несингулярная) задача Коши
\be{7.28}{\frac{ \pr^2 \, v }{\pr t^2} =\Delta_x v+G,~ v|_{t=0}=\varphi (x),~ v'_t|_{t=0}=\psi (x)}
корректно разрешима в полупространстве. Тогда в предположении, что $Y_{\alpha, \, t}=X^{-1}_{\alpha, \, t}$ получаем следующие начальные условия для уравнения ЭПД:
\be{7.29}{X_{\alpha} u|_{t=0}=a(x),~(X_{\alpha} u)'|_{t=0}=b(x).}
При этом различному выбору операторов преобразования  $X_{\alpha, t}$ (операторы Сонина, Бушмана--Эрдейи, Бушмана--Эрдейи нулевого порядка гладкости I или 2 рода, унитарные операторы третьего рода \eqref{6.16}, обобщенные операторы Бушмана--Эрдейи) будут соответствовать различные начальные условия. Следуя изложенной методике в каждом конкретном случае их можно привести к более простым аналитическим формулам.

Применение операторов преобразований позволяет сводить сингулярные (или иначе вырождающиеся) уравнения с операторами Бесселя по одной или нескольким переменным (уравнения ЭПД, сингулярное уравнение теплопроводности, $B$ --- эллиптические уравнения по определению И.А. Киприянова, уравнения обощённой осесимметрической теории потенциала --- теории $G A S P T$ --- А. Вайнстейна и другие) к несингулярным. При этом априорные оценки для сингулярного случая получаются как следствия соответствующих априорных оценок для регулярных уравнений, если т
 олько удалось оценить сами операторы преобразования  в нужных функциональных пространствах. Значительное число подобных оценок было приведено выше.
\vskip 03cm
В. Расмотрим оператор ${_1S_{0+}^{\nu}}$. Он имеет вид
\be{7.30}{{_1S_{0+}^{\nu}}=\frac{d}{dx}\int\limits_0^x K(\frac{x}{y}) f(y) \, dy,}
где ядро $K$ выражается по формуле $K(z)=P_{\nu}(z)$. Простейшие свойства специальных функций позволяют показать, что ${_1S_{0+}^{\nu}}$ можно рассматривать как частный случай оператора вида \eqref{7.30} с функцией Гегенбауэра в ядре
\be{7.31}{K(z)=\frac{\Gamma(\alpha+1)\Gamma(2 \beta)}{2^{p-\frac{1}{2}}\Gamma(\alpha+2\beta)\Gamma(\beta+ \frac{1}{2})}(z^{\alpha}-1)^{\beta-\frac{1}{2}}C^{\beta}_{\alpha}(z)}
при значениях параметров $\alpha=\nu, ~ \beta = \frac{1}{2}$ или с функцией Якоби в ядре
\be{7.32}{K(z)=\frac{\Gamma(\alpha+1)}{2^{\rho}\Gamma(\alpha+\rho+1)}(z-1)^{\rho}(z+1)^{\sigma}P^{(\rho, \sigma)}_{\alpha}(z)}
при значениях параметров $\alpha=\nu, ~ \rho=\sigma = 0$. Более общим являются операторы с гипергеометрической функцией Гаусса ${_2F_1}$ или с $G$ --- функцией  Майера в ядрах. Перейдём к их определению.

Операторы с функцией Гаусса ${_2F_1}(a, b; c; z)$ изучались в большом числе работ. Отметим здесь работы R.K. Saxena, S.L. Kalla (случай $z=y/x, ~ x/y$), A.C. Mabride, T.P. Higgins (случай $z=1-y^m/x^m$), E.R. Love (случай $z=1-y/x, ~ 1- x/y$), M. Saigo (для случая конечного отрезка $[a,b]$, $z=(x-y)/(x-a), ~ (y-x)/(b-x)$). Библиографию и другие примеры см. в [34, 353--358]. Подобные обобщения другого рода, в которых рассматриваются операторы с интегрированием по всей полуоси и ядра  выражаются через обобщённые функции Лежандра, изучались в [238--240]. Для исследования таких операторов могут оказаться полезными различные неравен
 тва для гипергеометрических функций, например, полученные в [241--251].

Введём ещё один класс таких операторов, обощающих операторы Бушмана--Эрдейи \eqref{71}--\eqref{72}.

Определение 17. Определим операторы Гаусса--Бушмана--Эрдейи по  следующим формулам:

\be{7.33}{{_1F_{0+}}(a, b, c)[f]=\frac{1}{2^{c-1}\Gamma(c)}.}
$$
\int\limits_0^x\lr{\frac{x}{y}-1}^{c-1}\lr{\frac{x}{y}+1}^{a+b-c} {_2F_1}\lr{^{a,b}_c| \frac{1}{2}-\frac{1}{2}\frac{x}{y}} f(y) \, dy,
$$
\be{7.34}{{_2F_{0+}}(a, b, c)[f]=\frac{1}{2^{c-1}\Gamma(c)}.}
$$
\int\limits_0^x\lr{\frac{y}{x}-1}^{c-1}\lr{\frac{y}{x}+1}^{a+b-c} {_2F_1}\lr{^{a,b}_c| \frac{1}{2}-\frac{1}{2}\frac{y}{x}} f(y) \, dy,
$$
\be{7.35}{{_1F_{-}}(a, b, c)[f]=\frac{1}{2^{c-1}\Gamma(c)}.}
$$
\int\limits_0^x\lr{\frac{y}{x}-1}^{c-1}\lr{\frac{y}{x}+1}^{a+b-c} {_2F_1}\lr{^{a,b}_c| \frac{1}{2}-\frac{1}{2}\frac{y}{x}} f(y) \, dy,
$$
\be{7.36}{{_2F_{-}}(a, b, c)[f]=\frac{1}{2^{c-1}\Gamma(c)}.}
$$
\int\limits_0^x\lr{\frac{x}{y}-1}^{c-1}\lr{\frac{x}{y}+1}^{a+b-c} {_2F_1}\lr{^{a,b}_c| \frac{1}{2}-\frac{1}{2}\frac{x}{y}} f(y) \, dy,
$$
\begin{eqnarray}
& {_3F_{0+}}[f]=\frac{d}{dx} {_1F_{0+}}[f], & {_4F_{0+}}[f]= {_2F_{0+}} \frac{d}{dx} [f], \label{7.37} \\
& {_3F_{-}}[f]={_1F_{-}}(-\frac{d}{dx} ) [f], & {_4F_{-}}[f]= (- \frac{d}{dx} ) {_2F_{-}}[f]. \label{7.38}
\end{eqnarray}

Символ ${_2F_1}$  в определениях \eqref{7.34} и \eqref{7.36} означает гипергеометрическую функцию Гаусса, а в определениях \eqref{7.33} и \eqref{7.35}
обозначает главную ветвь аналитического продолжения этой функции.

Операторы \eqref{7.33} -- \eqref{7.36} обобщают операторы Бушмана--Эрдейи \eqref{71} -- \eqref{72} соответственно. Они сводятся к последним при выборе параметров $a=-(\nu+\mu),~ b= 1+ \nu - \mu,~ c=1-\mu$. На операторы \eqref{7.33} -- \eqref{7.36} с соответствующими изменениями переносятся  все полученные выше результаты. В частности они факторизуются через более простые операторы \eqref{7.37} -- \eqref{7.38} при специальном выборе параметров.

Операторы \eqref{7.37} -- \eqref{7.38} обобщают операторы \eqref{73} -- \eqref{733}. Для них справедлива

\Theorem{ 42} Операторы \eqref{7.37} -- \eqref{7.38} могут быть расширены до изометричных в $L_2 (0, \infty)$ тогда и только тогда, когда они совпадают с операторами Бушмана--Эрдейи нулевого порядка гладкости I рода \eqref{73} -- \eqref{733} соответственно при целых значениях параметра $\nu=\frac{1}{2}(b-a-1)$.
\Endproc

Эта теорема выделяет операторы Бушмана--Эрдейи нулевого порядка гладкости среди их возможных обобщений \eqref{7.33} -- \eqref{7.38}. Сами операторы
\eqref{7.33} -- \eqref{7.36} интересны как новый класс преобразований, обобщающих операторы дробного интегродифференцирования. Аналогичные обобщения можно проделать и для операторов \eqref{6.1} -- \eqref{6.2}, \eqref{6.6}, \eqref{6.16} -- \eqref{6.17}.

Более общими являются операторы с $G$ -- функцией Майера в ядре. Например, один из таких операторов имеет вид

\begin{eqnarray}
& & {_1G_{0+}}(\alpha, \beta, \delta, \gamma)[f]= \frac{2^{\delta}}{\Gamma(1-\alpha)\Gamma(1-\beta)}\cdot \label{7.39} \\
& & \int\limits_0^x (\frac{x}{y}-1)^{-\delta}(\frac{x}{y}+1)^{1+\delta-\alpha-\beta} G_{2~2}^{1~2} \lr{\frac{x}{2y}-\frac{1}{2}|_{\gamma, \, \delta}^{\alpha, \, \beta}}  f(y)\, dy. \nonumber
\end{eqnarray}

Остальные получаются при изменении промежутка интегрирования и значений аргумента $G$ -- функции. При значениях параметров $\alpha=1-a,~ \beta=1-b, \delta=1-c, \gamma=0$ \eqref{7.39} сводится к \eqref{7.33}, а при значениях $\alpha=1+\nu,~ \beta=-\nu, \delta= \gamma=0$  \eqref{7.39} сводится к оператору Бушмана--Эрдейи I рода нулевого порядка гладкости ${_1S_{0+}^{\nu}}$.

Дальнейшие обобщения возможны в терминах функций Райта или Фокса.

\newpage

\section{8. В поисках вольтерровых унитарных операторов преобразования.}

Название этого пункта умышленно почти копирует название параграфа из известной работы Л.~Д.~Фаддеева [33]. Вопрос об унитарной эквивалентности возникает уже в простейшем случае двух матриц. Если общих содержательных критериев подобия не существует (кроме совпадения канонических форм), то для унитарного подобия получены интересные критерии Шпехта и Пирси [1]. Тем более интересным становится поиск унитарных ОП в бесконечномерном случае.

Ещё в работах Дельсарта и Лионса было замечено, что в классических пространствах  операторы преобразования Сонина--Пуассона--Дельсарта неограничены, так как изменяют гладкость функций. Тогда и возникла распространившаяся  в фольклоре задача, восходящая ещё к классикам теории ОП: обобщить операторы Сонина--Пуассона--Дельсартва с сохранением сплетающего свойства  так, чтобы обобщённые операторы действовали хотя бы ограниченно, а в идеальном случае были унитарными в пространстве $L_2(0,\infty)$.

Мечта о построении в явном виде унитарных ОП имеет и  корни в математической физике, а именно в задачах квантовой теории рассеяния. Там давно известны так называемые волновые операторы, которые сплетают возмущённый и невозмущённый операторы Шрёдингера и являются унитарными. Беда в том, что существование волновых операторов доказано и они имеют важные приложения, но их пока никто не построил в явном виде даже для простейших случаев за единичными исключениями. Это не удивительно, ведь и уравнение Шрёдингера решается в явном виде тольк
  для исключительных модельных случаев, с этим все уже смирились. А в теории ОП в определённом смысле наоборот---многие операторы построены в явном виде, но не для одного из них не изучалось свойство унитарности. Таким образом, решение задачи о построении унитарных обобщений ОП СПД позволяет соединить методы волновых операторов и ОП.

Решение поставленной задачи имеет некоторую историю. Первую модификацию ОП СПД предложил В.~В.~Катрахов в 1980 г. путём их домножения на оператор дробного интегрирования [233--234], но утверждение об унитарности построенных операторов для всех значений параметра $\nu$ содержало неточности. Оказалось, что построенные ОП унитарны лишь для целых значений параметра, что было показано автором в [252--253], см. также [254--260]. Следующее продвижение было связано с изучением ОП Бушмана--Эрдейи [226--229]. Окончательное решение задачи о построении унитарных обобщ
 ний ОП Сонина и Пуассона было получено автором [226--229] в рамках разработанного им \textit{композиционного метода} построения различных классов ОП [289--290, 373].

 Напомним, что для унитарности в $L_2(0,\infty)$ операторов Бушмана--Эрдейи нулевого порядка гладкости (\ref{73}--\ref{733}) необходимо и достаточно, чтобы число \  $\nu$ было целым.
Задача о построении унитарных ОП типа Сонина и Пуассона для оператора Бесселя (или углового момента) в общем случае для произвольного $\nu$ была окончательно решена в [289--290] в рамках композиционного метода (см. ниже). Это потребовало введения операторов достаточно сложной структуры.

\Theorem{ 43} Следующие операторы являются ОП типа Сонина и Пуассона, взаимно обратными и унитарными при всех $\nu$:
\begin{eqnarray}
\label{81}
S_U^\nu f =\cos \frac{\pi \nu}{2}\left(-\frac{d}{dx}\right) \int_x^\infty P_\nu\left(\frac{x}{y}\right)f(y)\,dy+\\\nonumber
+\frac{2}{\pi}\sin \frac{\pi \nu}{2}\Bigl(
\int_0^x \left(x^2-y^2\right)^{-\frac{1}{2}} Q_\nu^1 \left(\frac{x}{y}\right) f(y)\,dy - \\\nonumber
-\int_x^\infty \left(y^2-x^2\right)^{-\frac{1}{2}} \mathbb{Q}_\nu^1 \left(\frac{x}{y}\right) f(y)\,dy \Bigr),\\\nonumber
P_U^{\nu} f =\cos \frac{\pi \nu}{2} \Bigl( \int_0^{x} P_\nu \left(\frac{y}{x}\right)\left(\frac{d}{dy}\right)f(y)\,dy-\\\nonumber
-\frac{2}{\pi}\sin \frac{\pi \nu}{2}\bigl(
\int_0^x \left(x^2-y^2\right)^{-\frac{1}{2}} \mathbb{Q}_\nu^1 \left(\frac{y}{x}\right) f(y)\,dy - \\\nonumber
-\int_x^\infty \left(y^2-x^2\right)^{-\frac{1}{2}} Q_\nu^1 \left(\frac{y}{x}\right) f(y)\,dy \bigr)\Bigr),
\end{eqnarray}
где $P_\nu$---функция Лежандра первого рода, $Q_\nu^1$---функция Лежандра второго рода, $\mathbb{Q}_\nu^1$---функция Лежандра второго рода на разрезе [102].
\Endproc

Этим завершается история построения унитарных ОП типа Сонина и Пуассона.

Интересно рассмотреть частный случай, который вытекает из теоремы 19 при $\nu=1$. Получаем пару очень простых операторов
\begin{equation}
\label{82}
E_{0+}^{1,1}f=f(x)-\frac{1}{x}\int_0^x f(y)\,dy,\  E_{-}^{1,1}f=f(x)-\int_x^\infty \frac{f(y)}{y}\,dy,
\end{equation}
связанных со знаменитыми операторами Харди
\begin{equation}
\label{83}
H_1f=\frac{1}{x}\int_0^x f(y)\,dy, H_2f=\int_x^\infty \frac{f(y)}{y}\,dy.
\end{equation}
По поводу теории неравенств Харди см. [261--262]. Из наших результатов следует

\Theorem{ 43} Операторы (\ref{82}) образуют пару взаимнообратных унитарных в
$L_2(0,\infty)$ операторов. Они сплетают как ОП $\frac{d^2}{dx^2}$ и $\frac{d^2}{dx^2}-\frac{2}{x^2}$.
\Endproc

Как следует из (\ref{83}), операторы Бушмана--Эрдейи могут рассматриваться как \textit{обобщения операторов Харди}, а неравенства для их норм являются \textit{определёнными обобщениями неравенств Харди}, что позволяет взглянуть на этот класс операторов под новым интересным углом зрения. Кроме того, можно показать, что операторы (\ref{82}) являются преобразованиями Кэли от симметричных операторов $\pm 2i (xf(x))$ при соответствующем выборе областей определения. Их спектром является единичная окружность. В [226--227] эти вопросы рассмотрены и для пространств со
  степенным весом.

Результат об унитарности в частном случае теоремы 43 был недавно переоткрыт Куфнером, Перссоном и Малиграндой, давшими его элементарное доказательство. Теорема 19  позволяет выписать ещё несколько пар унитарных в $L_2(0,\infty)$  операторов  очень простого вида.
\begin{eqnarray*}
\label{84}
U_3f= f+\int_0^x f(y)\,\frac{dy}{y},\  U_4f= f+\frac{1}{x}\int_x^\infty f(y)\,dy,\\\nonumber
U_5f= f+3x\int_0^x f(y)\,\frac{dy}{y^2},\  U_6f= f-\frac{3}{x^2}\int_0^x y f(y)\,dy,\\\nonumber
U_7f=f+\frac{3}{x^2}\int_x^\infty y f(y)\,dy,\  U_8f=f-3x \int_x^\infty f(y)\frac{dy}{y^2},\\\nonumber
U_9f=f+\frac{1}{2}\int_0^x \left(\frac{15x^2}{y^3}-\frac{3}{y}\right)f(y)\,dy,\\\nonumber
U_{10}f=f+\frac{1}{2}\int_x^\infty \left(\frac{15y^2}{x^3}-\frac{3}{x}\right)f(y)\,dy.\\\nonumber
\end{eqnarray*}

Этот перечень можно существенно расширить. На мой взгляд подобных простых явных примеров очень не хватает в курсах функционального анализа.
\newpage

\section{9. Операторы преобразования для некоторых сингулярных дифференциальных операторов с переменными коэффициентами.}

В этом пункте будут рассмотрены ОП вида
\begin{equation}
\label{91}
S_\nu \left(B_\nu-q(x)\right)=B_\nu S_\nu, B_\nu=\frac{d^2}{dx^2}+\frac{\nu(\nu+1)}{x}\frac{d}{dx},
\end{equation}
которые обобщают ранее рассмотренные ОП для операторов Штурма--Лиувилля и Бесселя.

Первые работы по данной задаче рассматривали ОП на решениях соответствующих обыкновенных дифференциальных уравнений. При этом использовался метод функций Римана в том виде, как он был разработан в фундаментальной статье Б.~М.~Левитана [171]. Нужные ОП были построены в работах [263--264], в последней работе имеются неточности в формулировках. Следует отметить, что метод Б.~М.~Левитана имел, на мой взгляд, и один недостаток: им использовалось выражение функции Римана через гипергеометрическую функцию Гаусса, а последняя оценивалась на основан
 и общих асимптотик. В результате окончательные оценки для ядер ОП содержали неконтролируемые константы общего вида.

Определённый прорыв потребовал новых идей, он был совершен в замечательных работах [265--266] харьковского математика  Валерии Васильевны Сташевской (1924--2007). При этом была впервые установлена неожиданная связь между теорией ОП и теоремами типа Пэли--Винера из теории функций. Кандидатская диссертация В.В.~Сташевской была напечатана полностью отдельным изданием Харьковским университетом --- это единственный известный мне подобный случай.

Изложим  кратко план построения ОП вида (\ref{91}) по  В.~В.~Сташевской.

1. Рассматривается задача о нахождении решения дифференциального уравнения
\begin{equation}
\label{92}
y''-\left( q(x)+\frac{n(n-1)}{x^2}\right)y +\lambda^2 y=0,
\end{equation}
на полуоси $0<x<\infty$, где $n$---натуральное число, $q(x)$---вещественная функция, удовлетворяющая при $a>0, \mu <1/2$ ограничению
\begin{equation}
\label{93}
\int_0^a t^\mu |q(t)|\,dt<\infty,
\end{equation}
и условию в нуле
\begin{equation}
\label{94}
\lim_{x\to 0}\frac{y(x)}{x^n}=\frac{1}{2^{n-1/2}\Gamma(n+1/2)}.
\end{equation}

2. Доказывается существование решения поставленной задачи (\ref{92})--(\ref{94}) методом последовательных приближений. Это решение имеет вид
\begin{equation}
\label{95}
y(x)=y(x,\lambda)=\frac{\sqrt{\lambda x} J_\nu(\lambda x)}{\lambda ^n}+\frac{g(x,\lambda)}{\lambda ^n},
\end{equation}
где $J_\nu(t)$ --- функция Бесселя, $g(x,\lambda)$ есть целая функция от $\lambda$ при каждом фиксированном $x$.

3. Доказывается теорема Пэли--Винера для соответствующего преобразования Ханкеля (Фурье--Бесселя) полуцелой степени.

4. Устанавливается, что целая при каждом фиксированном $x$ функция $g(x,\lambda)$ из (\ref{95}) удовлетворяет всем условиям полученной теоремы  Пэли--Винера, в частности, имеет степень $\leq x$, поэтому для неё справедливо интегральное представление с интегрированием по промежутку $(0,x)$ с определённым ядром. Равенство Парсеваля приводит к нужной оценке этого ядра в пространстве $L_2(0,\infty)$.
На этом построение ОП закончено, для него получено интегральное представление и оценка ядра.

Впоследствии Н.~И.~Ахиезер доказал теорему Пэли--Винера и для преобразования Ханкеля любого порядка [267], что распространяет  результаты Сташевской и на все нецелые значения параметра $n$ в уравнении (\ref{92}).

Метод Сташевской излагается во всех монографиях по ОП, а также во многих книгах по обратным задачам и теории рассеяния. Он был обобщён для построения более общих ОП $S$, сплетающих по формуле
\begin{equation}
\label{96}
S(A-q(x))=AS, A=\frac{1}{v(x)}\frac{d}{dx}v(x)\frac{d}{dx}.
\end{equation}
Это направление разрабатывалось первоначально в работах Шебли [4--6, 268--271], а затем подробно изучалось в работах тунисского математика Халифы Тримеша и соавторов [4--6, 272--281].

Отметим, что подробное изучение ОП вида  (\ref{92}) и связанных с ними вопросов было проделано в работах другого харьковского математика А.С.~Сохина [378--381]. Вообще, как мы уже отмечали, вклад математиков Харьковской школы в развитие теории операторов преобразования очень велик: В.А.~Марченко, Н.И.~Ахиезер, Б.М.~Левитан, В.В.~Сташевская, А.Я.~Повзнер, Я.И.~Житомирский, А.С.~Сохин, Л.А.~Сахнович, В.Я.~Волк и другие.

Далее мы подробнее рассмотрим задачу о построении оператора преобразования $S_\alpha$ вида
\be{020}{ S_\alpha u(r)=u(r)+\int\limits_r^{\infty}P(r,t)u(t)\,dt ,}
определенного на любых функциях $u \in C^2(0, \infty)$ и сплетающего операторы
$B_\alpha - q(r)$ и $B_\alpha$ по формуле
\be{021}{S_{\alpha}(B_{\alpha}-q(r))u=B_{\alpha}S_{\alpha}u.}
Здесь $B_\alpha$ --- оператор Бесселя, который далее в этом пункте определяется так:
\be{022}{B_{\alpha}u =u''(r)+\frac{2\alpha}{r}u'(r),~ \alpha>0.}

Сформулированная задача о построении оператора преобразования по существу эквивалентна задаче  о нахождении решений дифференциального уравнения, коэффициенты которого имеют особенность в начале координат,
\be{023}{B_{\alpha}v(r)-q(r)v(r)+\lambda v(r)=0}
через решения невозмущенного уравнения $u(r)$, причем ищутся решения, представимые в виде (\ref{020}) с некоторым ядром $P(r,t)$.
Выбор пределов интегрирования в представлении (\ref{020}) приводит к тому, что искомое решение и его производная имеют ту же асимптотику, что и невозмущенное решение на бесконечности при выполнении очевидных требований к ядру $P(r, t)$.

Представление решений однородного уравнения (\ref{023}) по формуле (\ref{020}) обычно называется  представлением Йоста. Возможность такого представления с достаточно "хорошим"\  ядром $P$ для широкого класса потенциалов $q(r)$ лежит в основе классических методов решения обратных задач квантовой теории рассеяния.

Данный тип задач имеет свою историю.
Существование представления Йоста, сохраняющего асимптотику при $r \to \infty$, впервые было доказано для уравнения Штурма--Лиувилля ($\alpha=0,~ q$ --- суммируемая функция) Б.~Я.~Левиным в [45]. Операторы преобразования для оператора Бесселя впервые на русском языке были подробно изучены в широко известной работе Б.~М.~Левитана [171], и далее в их изучение основной вклад вносили математики Харьковской школы. Случай непрерывной $q,~ \alpha>0$ подробно рассмотрен в различных аспектах в работах  А.~С.~Сохина [378--381], а также ряда других авторов (см. подробнее [282]).
Оригинальная методика для решения поставленной задачи была разработана В.~В.~Сташевской [265--266], что позволило ей включить в рассмотрение сингулярные потенциалы с оценкой в нуле $|q(r)| \leq c x^{- 3/2+\varepsilon},~ \varepsilon > 0$ при целых $\alpha$, это методика получила широкое развитие.
В работе автора [255] условия на $q(r)$ были ослаблены до оценки $|q(r)| \leq C/r^2$ методом Сташевской.

Вместе с тем во многих математических и физических задачах необходимо рассматривать сильно сингулярные потенциалы, например, допускающие произвольную степенную особенность в нуле. В настоящей работе сформулированы результаты по интегральному представлению решений уравнений с подобными сингулярными потенциалами.

Далее приведём результаты автора, полученные в данном направлении. Изложение следует публикациям [255, 282--283, 382]

\begin{itemize}
\item Найдено интегральное  представление решений  для дифференциального уравнения с сингулярным потенциалом, имеющим  достаточно произвольную особенность в начале координат. От потенциала требуется лишь мажорируемость определенной функцией, суммируемой на бесконечности. В частности, к классу допустимых в данной работе относятся  сингулярный потенциал $q=r^{-2}$, сильно сингулярный потенциал со степенной особенностью $q=r^{-2-\varepsilon},~ \varepsilon > 0$, потенциалы Юкавы типа $q=e^{-\alpha r}/r$, потенциалы Баргмана и Батмана--Шадана~[31] и ряд других. При этом
 на функцию $q(r)$ не накладывается никаких дополнительных условий типа быстрой осцилляции в начале координат или знакопостоянства, что позволяет изучать притягивающие и отталкивающие потенциалы единым методом.

\item Доказательство существования оператора преобразования или интегрального представления решений возмущенного уравнения, что одно и то же, проводится по существу по известной схеме из работы Б.~М.~Левитана [171]. Мы вносим небольшое усовершенствование в эту схему, так как используемую в доказательстве функцию Грина как оказалось можно выразить не только через общую гипергеометрическую функцию Гаусса, но и более конкретно через функцию Лежандра, что позволяет избавиться от неопределенных постоянных в оценках из [171].

\item Ранее рассматривались лишь случаи одинаковых пределов (оба вида $[0;a]$ или $[a;\infty]$) в основном интегральном уравнении для ядра оператора преобразования. В данной работе впервые показано, что можно рассматривать случай различных пределов в основном интегральном уравнении. Именно такая расстановка пределов и позволила охватить более широкий класс потенциалов с особенностями в нуле. Кроме того, указанный способ позволяет рассматривать решения с более общими начальными условиями.

\item Изложенная техника полностью переносится и на задачу о построении неклассических операторов обобщенного сдвига. Данная задача по существу эквивалентна выражению решений уравнения
\be{024}{B_{\alpha,x} u(x,y) - q(x) u(x,y)=B_{\beta,y} u(x,y)}
через решения невозмущенного волнового уравнения при наличии дополнительных условий, обеспечивающих корректность.
\end{itemize}

Такие представления получаются уже из факта существования операторов преобразования и изучались для несингулярного случая ($\alpha=\beta=0$) в [24--26] как следствия общей теории обобщенного сдвига. Интересная оригинальная методика для получения подобных представлений также в несингулярном случае разработана в работах А.~В.~Боровских [286--287].
Из результатов настоящей работы следуют интегральные представления некоторого подкласса решений уравнения (\ref{024}) в общем сингулярном случае для достаточно произвольных потенциалов с особенностями в начале координат, причем оценки для решений не содержат неопределенных постоянных.

Как уже отмечалось, три задачи о построении оператора преобразования, представлении решений возмущенного уравнения и нахождении оператора обобщенного сдвига по существу эквивалентны, поэтому далее результаты приводятся для задачи о построении оператора преобразования.

Введем новые переменные и функции по формулам:
$$
\xi=\frac{t+r}{2}, ~ \eta=\frac{t-r}{2}, ~ \xi \geq \eta > 0;
$$
\be{025}{K(r, t)= \left(\frac{r}{t}\right)^\alpha P(r, t), ~ u(\xi, \eta)= K(\xi-\eta, \xi+\eta). }
Обозначим $\nu=\alpha-1$. Таким образом, для обоснования представления (\ref{020}) для решения  уравнения (\ref{023}) достаточно определить функцию $u(\xi, \eta )$. Доказано [255, 282--283, 382], что если существует дважды непрерывно дифференцируемое решение $u(\xi, \eta)$ интегрального уравнения
$$
u(\xi, \eta)=-\frac{1}{2}\int\limits_\xi^{\infty}R_\nu(s, 0; \xi, \eta) q(s) \, ds -
$$
$$
-\int\limits_\xi^{\infty} ds \int\limits_0^{\eta} q(s+\tau) R_{\nu}(s, \tau; \xi, \eta) u(s, \tau) \, d \tau,
$$
при условиях $ 0< \tau < \eta < \xi < s$,
то искомая функция $P(r, t)$ определяется по формулам (\ref{025}) через это решение $u(\xi, \eta)$. Функция $R_{\nu}=R_{\alpha-1}$ является функцией Римана, возникающей при решении следующей задачи для сингулярного гиперболического уравнения
$$
\frac{\partial^2 u(\xi, \eta)}{\partial \xi \partial \eta}+ \frac{4 \alpha(\alpha-1) \xi \eta}{(\xi^2-\eta^2)^2} u(\xi, \eta)=q(\xi+\eta)u(\xi, \eta),
$$
$$
u(\xi, 0)= - \frac{1}{2} \int\limits_\xi^\infty q(s)\, ds.
$$
Эта функция известна в явном виде, она выражается через гипергеометрическую функцию Гаусса $_2{F_1}$
$$
R_\nu=\left(\frac{s^2-\eta^2}{s^2-\tau^2}\cdot \frac{\xi^2-\tau^2}{\xi^2-\eta^2}\right)^{\nu}{_2F_1} \left(-\nu, -\nu; 1; \frac{s^2-\xi^2}{s^2-\eta^2}\cdot \frac{\eta^2-\tau^2}{\xi^2-\tau^2}\right).
$$
Это выражение упрощено в [255, 282], где показано, что функция Римана в рассматриваемом случае выражается через функцию Лежандра  по формуле
\be{026}{R_\nu (s, \tau, \xi, \eta)=P_\nu \left(\frac{1+A}{1-A}\right),~A=\frac{\eta^2-\tau^2}{\xi^2-\tau^2}\cdot \frac{s^2-\xi^2}{s^2-\eta^2}\ .}

Основной результат содержит

\Theorem{ 44}
Пусть функция $q(r)\in C^1 (0,\infty)$  удовлетворяет условию
\be{027}{|q(s+\tau)|\leq |p(s)|, ~ \forall s, \forall \tau, ~ 0< \tau <s,\ \int\limits_\xi^\infty |p(t)| \, dt < \infty, \forall \xi>0.}
Тогда существует интегральное представление вида (\ref{020}), ядро которого удовлетворяет оценке
$$
 \gathered
|P(r, t)| \leq \left(\frac{t}{r}\right)^ \alpha \frac {1}{2} \int\limits_{\frac{t+r}{2}}^\infty P_{\alpha-1}\left(\frac{y^2(t^2+r^2)-(t^2-r^2)}{2try^2}\right)|p(y)|\, dy \cdot
\\
\cdot \exp \left[ \left(\frac{t-r}{2}\right) \frac{1}{2} \int\limits_{\frac{t+r}{2}}^\infty P_{\alpha-1}\left(\frac{y^2(t^2+r^2)-(t^2-r^2)}{2try^2}\right)|p(y)|\, dy  \right].
\endgathered
$$
При этом ядро оператора преобразования $P(r,t)$, а также  решение уравнения (\ref{024}) являются дважды непрерывно дифференцируемыми на  $(0,\infty)$ по своим аргументам.
 \Endproc

 Перечислим классы потенциалов, для которых выполнены условия (\ref{027}). Если $|q(s)|$ монотонно убывает,  то можно принять $p(s)=|q(s)|$. Для потенциалов с произвольной особенностью в начале координат и возрастающих при $0<r<M$ (например, кулоновских $q=-\frac{1}{r}$), которые обрезаны нулем на бесконечности, $q(r)=0,~ r>M$, можно принять $p(s)=|q(M)|$, $s<M$, $p(s)=0$, $s \geq M$. Условию (\ref{027})  будут также удовлетворять потенциалы с оценкой $q(s+\tau) \leq c|q(s)|$ (на возможность подобного усиления теоремы 44 внимание автора обратил В.~В.~Катрахов).

Замечание. Фактически при доказательстве приведенной теоремы  не нужен явный вид функции Римана (\ref{026}).  Используется только  существование функции Римана, ее положительность и некоторое специальное свойство монотонности. Эти факты являются довольно общими, поэтому полученные результаты можно перенести на достаточно широкий класс дифференциальных уравнений.

Приведём упрощённую оценку для ядра ОП в первоначальных переменных. При её выводе существенным оказывается тот факт, что ядро выражено не через общую гипергеометрическую функцию Гаусса, а через более простую функцию Лежандра, для которой можно установить и использовать специальные свойства монотонности.

\Theorem{ 45}
Пусть выполнены условия теоремы 44. Тогда для ядра оператора преобразования $P(r, t)$ справедлива оценка
$$
 \gathered
 |P(r, t)| \leq \frac{1}{2} \left(\frac{t}{r}\right)^{\alpha} P_{\alpha-1} \left(\frac{t^2+r^2}{2tr}\right) \int\limits_r^\infty |p(y)|\, dy \cdot  \\
\cdot \exp \left[ \frac{1}{2} \left(\frac{t-r}{2}\right) P_{\alpha-1} \left(\frac{t^2+r^2}{2tr}\right) \int\limits_r^\infty |[p(y)| \, dy \right].
\endgathered
$$
 \Endproc

Отметим, что при $r \to 0$ ядро интегрального представления может иметь экспоненциальную особенность.

Для класса потенциалов со степенной сингулярностью вида
\be{029}{q(r)=r^{(-2 \beta +1 )},~ \beta > 0}
полученные оценки можно упростить не снижая их точности.

\Theorem{ 46} Рассмотрим потенциал вида (\ref{029}). Тогда теорема 44 выполняется с оценкой
$$
 |P(r, t)|  \leq \left(\frac{t}{r}\right)^{\alpha}  \frac{\Gamma(\beta)4^{\beta-1}}{(t^2-r^2)^\beta} \cdot P_{\alpha-1}^{- \beta} \left(\frac{t^2+r^2}{2tr}\right)\cdot
 $$
 $$
\cdot \exp\left[\left(\frac{t-r}{r}\right)  \frac{\Gamma(\beta)4^{\beta-1}}{(t^2-r^2)^\beta} P_{\alpha-1}^{- \beta} \left(\frac{t^2+r^2}{2tr}\right)   \right].
$$
 \Endproc

Данная оценка получается после довольно длинных вычислений с использованием знаменитой теоремы Слейтер--Маричева [235], которая помогает вычислить в терминах гипергеометрических функций необходимые интегралы после их сведения к свертке Меллина.

 Последняя оценка была получена в работе [255] для  потенциала $q(r)=cr^{-2}$, для которого $\beta=\frac{1}{2}$. Как следует из [102], в этом случае функция Лежандра $P_{\nu}^{-1/2}(z)$ может быть выражена через элементарные функции. Поэтому и соответствующая оценка   может быть выражена через элементарные функции.

Другим потенциалом, для которого рассматриваемая оценка оценка   может быть выражена через элементарные функции, является потенциал вида $q(r)=r^{-2 \beta + 1}$, для которого $\beta=\alpha-1$.

Результат подобный теоремам 44--46 можно использовать ещё двумя способами. Во--первых, он эквивалентен доказательству существования ООС, определяемых на решениях уравнения
\begin{equation}
\label{030}
\frac{\partial^2}{\partial x^2} u(x,y)+\frac{2\nu+1}{x}\frac{\partial}{\partial x} u(x,y)-q(x) u(x,y)=\frac{\partial^2}{\partial y^2} u(x,y).
\end{equation}
Во--вторых, получаются практически явные формулы для выражения решений уравнения (\ref{030}) через решения волнового уравнения. Метод получения явных интегральных представлений решений смешанных задач для уравнения вида (\ref{030})---это стандартное применение ОП и ООС, для случая второй производной он подробно разработан в [24--26]. Интересный другой подход к выводу подобных формул  для случая, когда сингулярное слагаемое в (\ref{030}) отсутствует (волновое уравнение для неоднородной среды),  предложен в [286--287], см. также [288].

Указанные результаты (построение ООС, представление решений смешанных задач через решения волнового уравнения) получаются и для более общего случая
\begin{eqnarray}
\label{031}
\frac{\partial^2}{\partial x^2} u(x,y)+\frac{2\nu+1}{x}\frac{\partial}{\partial x} u(x,y)-q(x) u(x,y)=\\ \nonumber
=\frac{\partial^2}{\partial y^2} u(x,y)+\frac{2\mu+1}{y}\frac{\partial}{\partial y} u(x,y)-p(y) u(x,y)
\end{eqnarray}
тем же методом, разработанным в [282--284].
\newpage

\section{10. Композиционный метод  построения операторов преобразования различных классов.}

Композиционный метод построения ОП был первоначально опубликован В.В.Катраховым и автором  в сборнике [289], посвящённом памяти новосибирского математика Бориса Алексеевича Бубнова, см. также [290, 373]. Этот метод основан на факторизации ОП через классические интегральные преобразования, такие как Фурье, Лапласа, Ханкеля и подобные им. При этом удаётся навести достаточно аккуратный порядок во всём многообразии ОП. Композиционный метод даёт алгоритмы не только для построения множества новых ОП, но содержит как частные случаи ОП СПД,  Векуа
 --Эрдейи--Лаундеса, Бушмана--Эрдейи, унитарные ОП (\ref{81}), обобщённые операторы Эрдейи--Кобера, а также введённые Р.~Кэрролом [4--6] классы эллиптических, гиперболических и параболических ОП и их обобщения.

Общая схема композиционного метода следующая. На вход подаётся пара операторов произвольного вида $A,B$, а также связанные с ними обобщённые преобразования Фурье $F(A), F(B)$, которые обратимы и действуют по формулам
\begin{equation}
\label{101}
F(A) A =g(t) F(A), F(B) B =g(t) F(B),
\end{equation}
где $t$---двойственная переменная (в простейших случаях классических интегральных преобразований можно принять $g(t)= - t^2$). На выходе получаем пару ОП, сплетающих $A$ и $B$.

\Theorem{ 47} Определим пару взаимно обратных операторов по формулам
\begin{equation}
\label{102}
S=F^{-1}(B) \frac{1}{w(t)} F(A), P=F^{-1}(A) w(t) F(B)
\end{equation}
с произвольной весовой функцией $w(t)$. Тогда они являются ОП,  которые удовлетворяют сплетающим соотношениям
$$SA=BS, \ \ \ PB=AP.$$
\Endproc

Разумеется, теорема 47 --- это только формальная схема для построения ОП методом композиции интегральных преобразований. Формулировка точного результата должна содержать перечисление пространств, в которых корректно действуют все задействованные операторы. Но можно рассматривать композиционный метод как формальный рецепт для построения с его помощью конкретных ОП, а уже после того, как получены явные формулы для конкретных ОП,  подбор нужных  функциональных пространств как правило является несложной задачей.

Несмотря на свою простоту, приведённая выше схема построения ОП композиционным методом позволяет получить все известные ранее в явном виде ОП, а также построить значительное число новых ОП  с заранее заданными полезными свойствами.

Отметим, что ввиду наличия весовых функций в предлагаемой схеме композиционного метода значительную роль в описании функциональных пространств, в которых действуют построенные ОП,  играют весовые оценки для различных классических интегральных преобразований. Это весьма разработанный раздел теории функций, содержащий большое число красивых и полезных результатов с многочисленными приложениями.

На практике чаще всего используются интегральные преобразования Фурье, синус и косинус преобразования, а также преобразование Ханкеля, которые мы определим так:
\begin{eqnarray*}
(Ff)(t)=\frac{1}{\sqrt{2\pi}}\int_0^\infty \exp(-ity)f(y)\,dy,\\
(F_c f)(t)=\sqrt{\frac{2}{\pi}}\int_0^\infty \cos(ty)f(y)\,dy,\\
(F_s f)(t)=\sqrt{\frac{2}{\pi}}\int_0^\infty \sin(ty)f(y)\,dy,\\
(H_\nu f)(t)=\frac{1}{t^\nu}\int_0^\infty J_\nu(ty)f(y)\,dy.
\end{eqnarray*}
При введённых нормировках все эти преобразования унитарны в $L_2(0,\infty)$ и за исключением преобразования Фурье три последних из них совпадают  со своими обратными [34, 291]. В дальнейшем сделаем самый очевидный выбор $g(t)= - t^2$.

В частном случае, когда строятся ОП, сплетающие оператор Бесселя и вторую производную, то есть $A=B_\nu, B=\frac{d^2}{dx^2}$ в обозначениях теоремы 47, схему композиционного метода можно конкретизировать.

\Theorem{ 48} Пусть в условиях предыдущей теоремы 47 рассматривается задача о построении ОП для пары сплетаемых операторов
$$A=B_\nu, B=\frac{d^2}{dx^2},$$
$F_c$ --- косинус преобразование Фурье, $H_\nu$ --- преобразование Ханкеля. Тогда оператор, построенный композиционным методом с произвольной весовой функцией $\frac{1}{w(t)}$ по формуле
\be{481}{S_\nu=F_c \left( \frac{1}{w(t)} H_\nu\right),}
на подходящих функциях является ОП типа Сонина. Для него справедливо интегральное представление
\be{482}{\left(S_\nu f \right)(x) = \sqrt{\frac{2}{\pi}}  \int_0^\infty K(x,y) f(y)\,dy,}
ядро которого выражается по формуле
\be{483}{K(x,y) = y^{\nu+1} \int_0^\infty \frac{\cos(xt) }{t^\nu \ w(t)} \ J_\nu(yt)\,dt,}
где $J_\nu(\cdot)$ --- функция Бесселя.

Формально обратный оператор, построенный композиционным методом с  весовой функцией $w(t)$ по формуле
\be{484}{P_\nu=H_\nu \left( w(t) F_c \right) ,}
на подходящих функциях является ОП типа Пуассона. Для него справедливо интегральное представление
\be{485}{\left(P_\nu  f \right)(x) = \sqrt{\frac{2}{\pi}}  \int_0^\infty G(x,y) f(y)\,dy,}
ядро которого выражается по формуле
\be{486}{G(x,y) = \frac{1}{x^\nu} \int_0^\infty w(t) \ t^{\nu+1} \cos(yt)  J_\nu(xt)\,dt,}
где $J_\nu(\cdot)$ --- функция Бесселя.
\Endproc

Рассмотрим кратко несколько примеров применения композиционного метода построения ОП из [289--290, 373].

Пример 1. Пусть $A=B_\nu, B=\frac{d^2}{dx^2}, F(A)=H_\nu, F(B)=F_c, w(t)=t^\alpha$.
Тогда при некоторых ограничениях на параметры получаем новое семейство обобщённых ОП Бушмана--Эрдейи, зависящих от двух параметров $\nu, \alpha$, которые в частных случаях являются обобщениями операторов  Эрдейи--Кобера, Сонина--Пуассона, а также ОП Бушмана--Эрдейи первого, второго и третьего родов. Свойства этих операторов подробно изучены автором, их ядра при всех значениях параметров выражаются через достаточно сложные комбинации функций Лежандра обоих родов. Эти результаты будут приведены в готовящейся к изданию монографии автора.

Рассмотрим самый простой частный случай, который привёл к построению рассмотренных выше унитарных представителей семейства ОП Бушмана--Эрдейи. Для этого напомним интегральное представление унитарных ОП Бушмана--Эрдейи третьего рода:
\begin{eqnarray}
\label{81}
S_U^\nu f =\cos \frac{\pi \nu}{2}\left(-\frac{d}{dx}\right) \int_x^\infty P_\nu\left(\frac{x}{y}\right)f(y)\,dy+\\\nonumber
+\frac{2}{\pi}\sin \frac{\pi \nu}{2}\Bigl(
\int_0^x \left(x^2-y^2\right)^{-\frac{1}{2}} Q_\nu^1 \left(\frac{x}{y}\right) f(y)\,dy - \\\nonumber
-\int_x^\infty \left(y^2-x^2\right)^{-\frac{1}{2}} \mathbb{Q}_\nu^1 \left(\frac{x}{y}\right) f(y)\,dy \Bigr),\\\nonumber
P_U^{\nu} f =\cos \frac{\pi \nu}{2} \Bigl( \int_0^{x} P_\nu \left(\frac{y}{x}\right)\left(\frac{d}{dy}\right)f(y)\,dy-\\\nonumber
-\frac{2}{\pi}\sin \frac{\pi \nu}{2}\bigl(
\int_0^x \left(x^2-y^2\right)^{-\frac{1}{2}} \mathbb{Q}_\nu^1 \left(\frac{y}{x}\right) f(y)\,dy - \\\nonumber
-\int_x^\infty \left(y^2-x^2\right)^{-\frac{1}{2}} Q_\nu^1 \left(\frac{y}{x}\right) f(y)\,dy \bigr)\Bigr),
\end{eqnarray}
где $P_\nu$---функция Лежандра первого рода, $Q_\nu^1$---функция Лежандра второго рода, $\mathbb{Q}_\nu^1$---функция Лежандра второго рода на разрезе [102].

При помощи композиционного метода получаются следующие естественные факторизации этих операторов через классические интегральные преобразования. Для этого достаточно в теореме 48 выбрать постоянную весовую функцию $w(t)=1$.

\Theorem{ 49} Для операторов Бушмана--Эрдейи третьего рода справедливы формулы факторизации через интегральные преобразования вида
\be{491}{S_U^\nu=F_c \ H_\nu , ~~~~~P_U^\nu=H_\nu \ F_c .}
\Endproc

Полученные формулы, представляющие операторы Бушмана--Эрдейи третьего рода в виде композиции двух унитарных преобразований, делают очевидным унитарность и самих этих операторов. Именно таким способом они и были впервые построены автором.

Пример 2. Пусть $A=B_\nu, B=B_\mu, F(A)=H_\nu, F(B)=H_\mu, w(t)=t^\alpha$. Тогда получаются операторы сдвига по параметру $T_{\nu,\mu}^{(\alpha)}$, ядра которых выражаются через гипергеометрические функции. Самым интересным из них представляется ОП при выборе $\alpha=0$
\begin{equation*}
(T_{\nu,\mu}^{(0)}f)(x)=\frac{2^{1-(\nu - \mu)}}{\Gamma(\nu - \mu)}
\int_x^\infty  y\left(y^2-x^2\right)^{\nu - \mu -1}f(y)\,dy,
\end{equation*}
который не зависит от самих значений $\nu, \mu$, а только от величины 'спуска' по разности параметров $\nu - \mu$. Это оператор типа Эрдейи--Кобера, открытый А.~Эрдейи [34].  Операторы спуска по параметру позволяют устанавливать формулы связи между решениями уравнений с операторами Бесселя с различными параметрами, это хорошо известный подход в теории уравнений с частными производными, например при спуске по размерности пространства к одномерному случаю. Другой подход к построению подобных операторов сдвига по параметру, это искать их в виде про
 зведения операторов типа Сонина и Пуассона $T_{\nu,\mu}=P_\mu S_\nu$, последние можно опять строить композиционным методом, причём с разными весовыми функциями.
Этот результат мы сформулируем отдельно.

\Theorem{ 50} Пусть дана пара $S_\nu, P_\mu$ ОП типа Сонина и Пуассона, удовлетворяющая сплетающим соотношениям
\be{501}{S_\nu B_\nu = \frac{d^2}{dx^2} S_\nu, P_\mu \frac{d^2}{dx^2} = B_\mu P_\mu, B_\gamma = \frac{d^2}{dx^2} + \frac{2\gamma+1}{x} \frac{d}{dx}}
между второй производной и оператором Бесселя $B_\gamma.$
Тогда по данной паре ОП   типа Сонина и Пуассона и произвольного коммутирующего со второй производной оператора $K(x)$ можно определить оператор сдвига по параметру по формулам
\be{}{T_{\nu,\mu}=P_\mu K(x) S_\nu, ~~~ K(x) \frac{d^2}{dx^2} = \frac{d^2}{dx^2} K(x),}
который будет удовлетворять определяющему соотношению
\be{503}{T_{\nu,\mu} B_\nu = B_\mu T_{\nu,\mu}}
и осуществлять сдвиг по параметру дифференциального оператора Бесселя.
\Endproc

Эта теорема ещё раз подчёркивает важность для теории ОП тех операторов, которые коммутируют со второй производной или вообще с производными. Например, в качестве таких операторов можно выбрать дробные интегралы Римана--Лиувилля на подходящих функциях. В качестве частного случая при $\nu=\mu$ из теоремы 50 получается алгоритм для конструирования операторов, которые коммутируют на подходящих функциях с дифференциальным оператором Бесселя.

Пример 3. При значениях $\nu=\mu$ из предыдущего примера  2 получается семейство операторов, коммутирующих с оператором Бесселя, ядра которых выражаются через функции Лежандра первого и второго родов.

Пример 4. Пусть выбрано $A=B_\mu, B=B_\nu - \lambda^2, w(t)=t^\mu\left(t^2+\lambda^2\right)^{-\frac{\mu}{2}}$,
\begin{equation}
\label{103}
F(A)=F(B)=\frac{1}{t^\nu} \int_0^\infty y^{\nu+1}J_\nu \left(y\sqrt{t^2+\lambda^2}\right) f(y)\,dy.
\end{equation}
Тогда по формулам (\ref{102}) получаем ОП типа Векуа--Эрдейи--Лаундеса вида
\begin{equation}
\label{104}
Sf=\lambda^{1+\nu - \mu} \int_0^\infty y\left(y^2-x^2\right)^{\frac{\nu - \mu -1}{2}}J_{\nu - \mu -1}\left(\lambda\sqrt{y^2-x^2}\right) f(y)\,dy.
\end{equation}

Это в точности ОП, введённый Лаундесом [34], он также зависит только от величины сдвига по параметрам. Отметим, что на основе композиционного метода можно изложить всю теорию ОП  Векуа--Эрдейи--Лаундеса, при этом роль обобщённого преобразования Ханкеля играет выражение (\ref{103}).

Обобщая классификацию Р.~Кэррола, мы назвали ОП для оператора Бесселя из приведённых примеров \textit{$B$--гипер\-бо\-лическими}. В работах [289--290] также рассмотрено применение композиционного метода к построению \textit{$B$--эллиптических} ОП со свойством $TB_\nu=-D^2T$ и \textit{$B$--параболических} ОП со свойством $TB_\nu=D T$. Эта классификация естественна, так как основана на типе уравнений в частных производных, которому удовлетворяют ядра операторов соответствующих преобразования.

В указанном методе может быть использован и оператор квадратичного преобразования Фурье (КПФ, дробное преобразование Фурье, преобразование Фурье--Френеля), см. [292--294]. Это важное интегральное преобразование недостаточно широко известно (пока), оно возникло из предложения Френеля заменить стандартные плоские волны с линейными аргументами в экспонентах на более общие волны с квадратичными аргументами в экспонентах, что позволило полностью объяснить парадоксы со спектральными линиями при дифракции Фраунгофера. Математически операто
 ры КПФ являются дробными степенями $F^\alpha$ обычного преобразования Фурье, достраивая его до полугруппы по параметру $\alpha$, они были определены Н.~Винером и А.~Вейлем. В теории всплесков, в которой принято каждую формулу считать новой и называть по--новому давно известные вещи, КПФ называется преобразованием Габора. Изложенный выше композиционный метод позволяет с помощью этого преобразования строить ОП для одномерного оператора Шрёдингера из квантовой механики. При этом может быть использовано и более общее квадратичное преобразован
 е Ханкеля (КПХ) ~[295].

Кратко изложим основные факты, относящиеся к теории дробного или квадратичного преобразования Фурье--Френеля, используя материалы диссертации Д.Б.\,Карпа [383]. К этому тексту мы пока отсылаем читателя и по поводу всех ссылок по теории КПФ.

Целые степени  (орбита) классического преобразования Фурье образуют циклическую группу порядка 4, при этом четвёртая степень этого преобразования даёт тождественный оператор. Поэтому, в частности,  спектр классического преобразования Фурье в $L_2(-\infty,\infty)$ состоит из четырех
точек, расположенных на единичной окружности:  $1$, $i$, $-1$ и $-i$.  . Идея включить эту
дискретную группу в непрерывную со спектром, целиком заполняющим единичную окружность,
принадлежит Винеру  и была реализована им в работе 1929 года.  Несколько
позже Кондон  (в 1937)  а затем Кобер  (в 1939) независимо
переоткрыли эту группу, которая стала называться дробным преобразованием Фурье (ДПФ).
Баргманну  принадлежит обобщение на многомерный случай.

Впоследствии дробное преобразование Фурье неоднократно переоткрывалось целым рядом
авторов.  В книге Антосика, Микусинского и Сикорского  под названием
''преобразование Фурье-Мелера'' упоминается циклическая группа произвольного порядка, в
которую можно включить преобразование Фурье. ДПФ изучалось Гинандом  и
Вольфом. Вавржинчик в  приходит к ДПФ рассматривая
классическое преобразование Фурье в виде экспоненты от производящего оператора.
В.Ф. Осиповым  независимо от предыдущих авторов построена
теория ДПФ на группах и введены соответствующие почти--периодические функции Бора--Френеля, изучены асимптотические свойства этого преобразования, рассмотрены приложения в гармоническом анализе и теории чисел [384--385]. Намиас  переоткрыл дробное преобразования Фурье и
использовал его для решения некоторых задач для уравнения Шр\"{e}\-дин\-гера.  Керр
исследовала дробное преобразование Фурье в пространстве $L_2$  и
пространстве Шварца $S$.

Отдельное направление исследований связано с дробным преобразованием Фурье при чисто
мнимых значениях группового параметра. В этом случае чаще всего встречается название
''полугруппа Эрмита''. Начало этому направлению положила статья Хилле 1926 года, в которой оператор ДПФ при мнимых значениях параметра возникает в связи с
суммированием методом Абеля разложений по полиномам Эрмита. Позднее, ''полугруппа Эрмита''
была использована Бабенко, Бекнером и Вайслером
 для получения неравенств в теории классического преобразования Фурье.

Отметим, что  квадратичное преобразование Фурье является одной из двух основных составляющих (наряду с неравенствами для средних значений в комплексной плоскости специального вида), которые были использованы сначала К.И.\,Бабенко для частных случаев, а затем Бекнером для общего случая при доказательстве знаменитых условий ограниченности обычного   преобразования Фурье в пространствах $L_p$ с точными постоянными. Другим интересным применением  квадратичного преобразования Фурье является круг вопросов, связанных со знаменитой задаче
 й Паули по определению функции по некоторому набору спектральных данных. Подобные задачи, обычно неразрешимые для классического преобразования, находят элегантные решения в рамках КПФФ.

Применения дробного преобразования Фурье столь же разнообразны и обширны как и
классического. Перечислим лишь некоторые.  Приложениям в квантовой механике посвящена уже
упоминавшаяся работа Намиаса. Использование ДПФ в оптике и анализе сигналов
разрабатывается  группой исследователей под руководством  Оцактаса. Участниками этой
группы написана книга , целиком по\-свя\-щ\"{e}н\-ная теории и
приложениям дробного преобразования Фурье [386],  в которой процитировано несколько сотен работ по данной теме. Автор с удовольствием приводит эту ссылку, так как это единственная западная монография, в которой мне выражена благодарность во введении за полезные обсуждения. Следует заметить, что Халдун  Оцактас --- это один из крупнейших специалистов в области математической оптики и её приложений, работавший долгое время в Стэнфорде.
Мастардом найдены аналоги неравенства Гейзенберга, инвариантные
относительно дробного преобразования Фурье.  Было показано, что
многомерное преобразование Вигнера равно корню шестой степени из обратного преобразования
Фурье, и, следовательно, также является частным случаем дробного преобразования Фурье.
Бьюн  получил некоторые новые формулы обращения для дробного преобразования
Фурье.  Аналог дробного преобразования Фурье для $q$--полиномов Эрмита был введён
Аски, Н.М. Атакишиевым и С.К. Сусловым.  Дальнейшие обобщения на операторы с ядрами
в форме билинейных производящих функций полиномов Аски--Вилсона рассматривались в
работах Аски и Рахмана.

Дробное преобразование Ханкеля (ДПХ) гораздо менее изучено.  Оно было введено Кобером в
 и изучалось Гинандом.  Затем было несколько раз переоткрыто, например,  Намиасом.  ДПХ было рассмотрено
Керр при действительных значениях группового параметра в пространстве $L_2(0,\infty)$  и в пространствах Земаняна.

В диссертации Д.Б.\,Карпа [383] рассмотрен достаточно общий подход к построению подобных преобразований при помощи разложения в ряды по известным системам  ортогональных функций. В частности, при выборе системы функций Эрмита получается классическое преобразование Фурье, при выборе системы функций Лагерра --- квадратичное преобразование Фурье--Френеля, а при выборе систем функций Лежандра, Чебышёва или Гегенбауэра построены новые полугруппы интегральных преобразований.

Выпишем явно интегральные формулы для квадратичных преобразований Фурье--Френеля (КПФФ) и Ханкеля (КПХ)
\begin{equation}\label{FrF}
(F^{\alpha}f)(y)=\frac{1}{\sqrt{\pi(1-e^{2i\alpha})}}\int\limits_{-\infty}^{\infty}\!\!
e^{-\frac{1}{2}i(x^2+y^2)\ctg\alpha}e^{ixy\cosec\alpha}f(x)dx;
\end{equation}

\begin{equation}\label{FrH}
(H_\nu^{\alpha}f)(y)=\frac{2\left(-e^{i\alpha}\right)^{-\frac{\nu}{2}}}{1-e^{i\alpha}}
\int\limits_0^{\infty}\!\!e^{-\frac{1}{2}i\ctg\frac{\alpha}{2}\left(x^{2}+y^{2}\right)}
(xy)^{\frac{1}{2}}J_{\nu}\left(\frac{2xy\sqrt{-e^{i\alpha}}}{1-e^{i\alpha}}\right)f(x)dx;
\end{equation}

Рассмотрим кратко формулы преобразований, которые образуют операционное исчисление для КПХ.

Введ\"{e}м следующие дифференциальные операторы
$$
A^{-}_\nu=x^{\nu +\frac{1}{2}}e^{-\frac{x^2}{2}}\frac{d}{dx} x^{-\nu
-\frac{1}{2}}e^{\frac{x^2}{2}}=-\frac{\nu + \frac{1}{2}}{x}+x+\frac{d}{dx},
$$
$$
A^{+}_\nu=x^{-\nu -\frac{1}{2}}e^{\frac{x^2}{2}} \frac{d}{dx}x^{\nu
+\frac{1}{2}}e^{-\frac{x^2}{2}}=\frac{\nu+\frac{1}{2}}{x}-x+\frac{d}{dx},
$$
$$
N_\nu=x^{\nu+\frac{1}{2}}\frac{d}{dx}x^{-\nu-\frac{1}{2}}=-\frac{\nu+\frac{1}{2}}{x}+\frac{d}{dx},
$$
$$
M_\nu=x^{-\nu-\frac{1}{2}}\frac{d}{dx}x^{\nu+\frac{1}{2}}=\frac{\nu+\frac{1}{2}}{x}+\frac{d}{dx}.
$$
Эти операторы связаны соотношениями
\begin{eqnarray}\label{AMN}
L_{\nu}= -\frac{1}{4}D^2-\frac{\nu^2-1/4}{x^2} + \frac{1}{4}x^2 - \frac{\nu+1}{2}=-\frac{1}{4}A^{+}_{\nu}A^{-}_{\nu},\nonumber\\~~~~A^{-}_{\nu}=N_{\nu}+x,
~~~~A^{+}_{\nu}=M_{\nu}-x,
\end{eqnarray}
\begin{equation}\label{LMN}
L_{\nu}=M_{\nu} N_{\nu},~~~~
x\frac{d}{dx}+\frac{d}{dx}x=N_{\nu}x+xM_{\nu}=M_{\nu}x+xN_{\nu},
\end{equation}

Обозначим через $X$ оператор умножения на независимую переменную.  Теперь мы можем, следуя [383],  выписать набор формул преобразования операций для КПХ.

\begin{equation}\label{H+tD}
H_{\nu+1}^{\alpha}Xf=\frac{1}{2}\left[(e^{-i\alpha}-1)N_{\nu}+
(e^{-i\alpha}+1)X\right]H_{\nu}^{\alpha}f,
\end{equation}
\begin{equation}\label{H+ND}
H_{\nu+1}^{\alpha}N_{\nu}f=\frac{1}{2}\left[(e^{-i\alpha}+1)N_{\nu}+
(e^{-i\alpha}-1)X\right]H_{\nu}^{\alpha}f,
\end{equation}
\begin{equation}\label{NHD}
N_{\nu}H_{\nu}^{\alpha}f=
\frac{1}{2}H_{\nu+1}^{\alpha}\left[(e^{i\alpha}+1)N_{\nu}+(e^{i\alpha}-1)X\right]f,
\end{equation}
\begin{equation}\label{xHD}
XH_{\nu}^{\alpha}f=\frac{1}{2}H_{\nu+1}^{\alpha}\left[(e^{i\alpha}-1)N_{\nu}+
(e^{i\alpha}+1)X\right]f,
\end{equation}
\begin{equation}\label{HtD}
H_{\nu}^{\alpha}Xf=
\frac{1}{2}\left[(1-e^{i\alpha})M_{\nu}+(1+e^{i\alpha})X\right]H_{\nu+1}^{\alpha}f,
\end{equation}
\begin{equation}\label{HMD}
H_{\nu}^{\alpha}M_{\nu}f=
\frac{1}{2}\left[(1+e^{i\alpha})M_{\nu}+(1-e^{i\alpha})X\right]H_{\nu+1}^{\alpha}f,
\end{equation}
\begin{equation}\label{MH+D}
M_{\nu}H_{\nu+1}^{\alpha}f=
\frac{1}{2}H_{\nu}^{\alpha}\left[(1+e^{-i\alpha})M_{\nu}+(1-e^{-i\alpha})X\right]f,
\end{equation}
\begin{equation}\label{xH+D}
XH_{\nu+1}^{\alpha}f=
\frac{1}{2}H_{\nu}^{\alpha}\left[(1-e^{-i\alpha})M_{\nu}+(1+e^{-i\alpha})X\right]f,
\end{equation}
\begin{equation}\label{A-H}
H_{\nu+1}^{\alpha}A^{-}_{\nu}f=e^{-i\alpha}A^{-}_{\nu}H_{\nu}^{\alpha}f,
~~~~A^{-}_{\nu}H_{\nu}^{\alpha}f=H_{\nu+1}^{\alpha}e^{i\alpha}A^{-}_{\nu}f,
\end{equation}
\begin{equation}\label{A+H}
H_{\nu}^{\alpha}A^{+}_{\nu}f=e^{i\alpha}A^{+}_{\nu}H_{\nu+1}^{\alpha}f,~~~~
A^{+}_{\nu}H_{\nu+1}^{\alpha}f=H_{\nu}^{\alpha}e^{-i\alpha}A^{+}_{\nu}f,
\end{equation}
\begin{equation}\label{Ht^2D}
H_{\nu}^{\alpha}X^{2}f=
\left[X^{2}\cos^{2}{\frac{\alpha}{2}}-\frac{1}{2}i\sin{\alpha}\left[XD+DX\right]-
\sin^{2}{\frac{\alpha}{2}}L_{\nu}\right]H_{\nu}^{\alpha}f,
\end{equation}
\begin{equation}\label{x^2HD}
X^{2}H_{\nu}^{\alpha}f=
H_{\nu}^{\alpha}\left[X^{2}\cos^{2}{\frac{\alpha}{2}}+\frac{1}{2}i\sin{\alpha}\left[XD+DX\right]-
\sin^{2}{\frac{\alpha}{2}}L_{\nu}\right]f,
\end{equation}
\begin{equation}\label{HLD}
H_{\nu}^{\alpha}L_{\nu}f=
\left[-X^{2}\sin^{2}{\frac{\alpha}{2}}-\frac{1}{2}i\sin{\alpha}\left[XD+DX\right]+
\cos^{2}{\frac{\alpha}{2}}L_{\nu}\right]H_{\nu}^{\alpha}f,
\end{equation}
\begin{equation}\label{LHD}
L_{\nu}H_{\nu}^{\alpha}f=
H_{\nu}^{\alpha}\left[-X^{2}\sin^{2}{\frac{\alpha}{2}}+\frac{1}{2}i\sin{\alpha}\left[XD+DX\right]+
\cos^{2}{\frac{\alpha}{2}}L_{\nu}\right]f,
\end{equation}
\begin{equation}\label{HtD+DtD}
H^\alpha_\nu\left[XD+DX\right]f=
\left[-i\sin{\alpha}(X^2+L_\nu)+\cos\alpha\left[XD+DX\right]\right]H_{\nu}^{\alpha}f,
\end{equation}
\begin{equation}\label{xD+DxHD}
\left[XD+DX\right]H_{\nu}^{\alpha}f=H_{\nu}^{\alpha}\left[i\sin{\alpha}(X^2+L_\nu)+
\cos\alpha\left[XD+DX\right]\right]f
\end{equation}

Пора вернуться к теме композиционного метода построения ОП. Для этих целей наиболее интересны соотношения (\ref{Ht^2D}--\ref{x^2HD}). Согласно общей схемы из теоремы 47 мы можем сконструировать ОП, сплетающие дифференциальные операторы $D^2$ и
$$
\left(\sin^{2}{\frac{\alpha}{2}}L_{\nu} - \frac{1}{2}i\sin{\alpha}\left(XD+DX\right) - X^{2}\cos^{2}{\frac{\alpha}{2}}
\right),
$$
где
$$
L_{\nu}= -\frac{1}{4}D^2-\frac{\nu^2-1/4}{x^2} + \frac{1}{4}x^2 - \frac{\nu+1}{2}
$$
с произвольными параметрами $\alpha, \nu$. Для этого в обозначениях теоремы 47 надо положить
$$
A=\left(\sin^{2}{\frac{\alpha}{2}}L_{\nu} - \frac{1}{2}i\sin{\alpha}\left(XD+DX\right) - X^{2}\cos^{2}{\frac{\alpha}{2}}\right), B=D^2,
$$
$$
F(A)=H_\nu^\alpha, F(B)=F_c, g(t)=-t^2,
$$
где $H_\nu^\alpha$ --- квадратичное преобразование Фурье--Френеля, $F_c$ --- косинус преобразование Фурье.

О пользе ОП для теории дискретного преобразования Фурье см. [296].
\newpage

\section{11. Некоторые приложения метода операторов преобразования.}

В этом пункте очень кратко перечислены некоторые приложения ОП.

1. \textbf{Теория обратных задач}. В самом простом случае это задача о восстановлении уравнения Штурма--Лиувилля по так называемой спектральной функции. Единственность была доказана В.~А.~Марченко, общий метод решения заключается в выписывании по спектральной функции некоторого интегрального уравнения для ядра ОП. Получившееся интегральное уравнение называется уравнением Гельфанда--Левитана (можно посплетничать по поводу этого названия, но мы не будем этого делать), потенциал по ядру восстанавливается из соотношения (\ref{tag10}). Подробное
 зложения вопроса в [4--6, 22--23, 28--29, 297--298]. Отметим, что М.~Г.~Крейном был создан свой метод решения обратных задач, но он не использует технику ОП.

В настоящее время теория обратных задач --- это огромный бурно развивающийся раздел современной математики. Просто приведём список известных монографий  Марченко В.А., Левитана Б.М., Фаддеева Л.Д., Кэррола Р.,  Захарова, Манакова, Новикова, Питаевского, Абловица М. и Сигура Х., Рамма А.Г., Шадана К. и  Сабатье П., Абловица М. и Сигура М., Юрко В.А., Исакова В. Нелинейные обратные задачи для гиперболических уравнений в
различных постановках изучались в работах М.М. Лаврентьева,
В.Г.~Романова, Ю.Е. Аниконова, Б.А. Бубнова, С.И. Кабанихина, А.И.\,Кожанова, А.И.
Прилепко, А.Х. Амирова, Е.Г. Саватеева, Е.С. Глушковой, Д.И.
Глушковой, Т.Ж. Елдесбаева, A. Lorenzi, А.М. Денисова, M.
Grasseli, М. Клибанова, М. Ямамото.

2.  \textbf{Теория рассеяния}. В самом простом случае это задача о восстановлении уравнения Штурма--Лиувилля по данным рассеяния, например, по двум спектрам. Схема решения та же, что и для обратных задач, определяющее уравнение в этом случае называется уравнением Марченко, см. [4--6, 21--23, 31--33, 299--306]. Через ОП факторизуется сам оператор рассеяния. Метод ОП является также основным при исследовании решений Йоста и задач для системы Дирака [29]. Укажем также работы [390--393].

3.  \textbf{Операторы Штурма--Лиувилля и их обобщения}. Можно без преувеличения сказать, что методы ОП тривиализировали эту теорию, предоставив единый мощный метод для большинства задач. К ним относятся асимптотические формулы для решений и собственных значений, формулы следов, асимптотика спектральной функции, см. [22--23, 28--29]. В работах [76--77] строятся ОП для уравнений вида Штурма--Лиувилля, но с изменённым вхождением спектрального параметра в потенциал. В монографии А.~П.~Хромова [40] рассматриваются ОП для более общих операторов, содержащих ди
 фференциальные и интегральные части, а также вопросы подобия двух операторов Вольтерра. Эти результаты относятся к общей схеме применения ОП для того случая, когда сплетаемые операторы интегральные или интегро--дифференциальные, а сам ОП является интегральным.

Построены ОП и для так называемого оператора Эйри из квантовой механики вида $D^2 -x$. В этом случае ядра выражаются через функции Эйри. По--видимому, аккуратных доказательств существования и формул представления со всеми деталями для этого случая не было опубликовано, краткие выкладки для случая возмущённого оператора Штарка  $D^2 -x - q(x)$ приведены в [387]. В указанной работе приведены интересные приложения ОП данного вида к обратной задаче рассеяния и математическому описанию эффекта Штарка, построению решений Йоста и уравнения Левитана, о
 писанию спектров и решению уравнений Липмана--Швингера, связям с волновыми операторами и преобразованием Фурье--Эйри, нахождению условий унитарности и построению матрицы рассеяния.

В теории интегральных уравнений применяется метод операторных тождеств, разработанный В.А.\,Амбарцумяном, Л.А.\,Сахновичем и другими авторами. Его сутью является рассмотрение интегральных уравнений, ядра которых удовлетворяют гиперболическим уравнениям, аналогичным уравнениям для ядер ОП, возникающих при построении ОП для операторов типа Штурма--Лиувилля.
Решения таких интегральных уравнений  имеют специфическую природу и интересные приложения, см. работу Е.А.\,Аршавы [388].

4. \textbf{Нелинейные уравнения}. В теории нелинейных дифференциальных уравнений метод Лакса положил начало использованию  ОП для доказательства существования решений и построения солитонов. ОП используются также по  схеме Лакса в методе обратной задачи рассеяния (МОЗР), а также
так называемом методе одевания и преобразованиях Дарбу [303--306].
Эти результаты относятся к общей схеме применения ОП для того случая, когда и сплетаемые операторы,  и сам ОП является дифференциальным. К такому типу относятся многие дифференциальные подстановки, линеаризирующие соответствующие нелинейные дифференциальные уравнения, например, подстановки Миуры.

5. \textbf{Сингулярные и вырождающиеся краевые задачи для уравнений с частными производными}. В работе [234] В.~В.~Катраховым был предложен новый подход к постановке краевых задач для эллиптических уравнений с особенностями. Например, для уравнения Пуассона им рассматривалась задача в области, содержащей начало координат, в которой решения могут иметь особенности произвольного роста. В точке начала координат им было предложено новое нелокальное краевое условие типа свертки, которое мы назовём К---следом. В определении классов для решений, к
 оторые обобщают пространства С.~Л.~Соболева на случай функций с существенными особенностями, фундаментальную роль играют различные ОП. Основные результаты состоят в доказательстве корректной разрешимости поставленных задач во введённых пространствах. На более сложные уравнения и области в этом направлении  результаты обобщались в работах [307--311].

Отдельно выделим теорию ОП, построенную для дифференциальных операторов гипербесселева типа высоких порядков. Основные результаты по данной тематике получены чл.--корр. Болгарской национальной АН И.Х.\,Димовски и его учениками. Многочисленные приложения получили ОП Сонина--Димовски и Пуассона--Димовски.

6. \textbf{Сингулярные псевдодифференциальные операторы}. Подобные ПДО, связанные с оператором Бесселя, были определены и изучены В.~В.~Катраховым и И.~А.~Киприяновым [312]. Р.~Кэррол посвятил этим ПДО Киприянова--Катрахова отдельную главу в книге [5], существенно переработав изложение материала.

7. \textbf{Задачи об убывании решений дифференциальных уравнений}.

Е.~М.~Ландисом была поставлена такая задача [313]: доказать, что любое решение во всём пространстве уравнения $\Delta u-q(x)u=0$ при условиях $|q(x)|\leq\lambda^2, |u(x)|\leq  A \exp(-(\lambda+\varepsilon)|x|, \lambda\geq 0, \varepsilon>0$ равно нулю. В работах [314--316] автором было дано частичное решение этой задачи для случаев радиально-симметричного или зависящего от одной переменной потенциалов, а также для некоторых ультрагиперболических уравнений. Решение основывалось на использовании ОП со специальными оценками ядер прямых и обратных операторов для этого случая. Для общих потенциал
 ов утверждение задачи оказалось неверным, В.~З.~Мешковым был доказан удивительный факт [317--318], что существуют решения со скоростью убывания $\exp(|x|^{-4/3})$. Доказательства Мешкова потребовали обобщения известных неравенств Карлемана и использовали методы теории чисел. Отметим, что по данным математического портала www.mathnet.ru на статью Виктора Захаровича Мешкова ссылаются в трёх работах, проиндексированных на этом портале. Зато одна из ссылок --- в статье выдающегося математика, лауреата премии Филдса Жана Бургейна, другая --- в совместной ста
 тье Бургейна и Карлоса Кёнига, третья --- в знаменитой работе В.И.\,Юдовича "одиннадцать великих проблем математической гидродинамики".

7. \textbf{Применение ОП к доказательству вложений функциональных пространств}.

В [226--227] ОП Бушмана--Эрдейи применены автором к доказательству вложений и изометрий пространств И.~А.~Киприянова [175--177] в весовые пространства С.~Л.~Соболева. При этом основную роль играют описанные выше оценки норм и свойство унитарности. Подобные вопросы рассматривались также в [319--322]. Отметим, что есть противоречия в формулировках окончательных результатов между моими оценками и оценками из работы [321]. Это некоторое время расстраивало меня, но потом удалось разобраться, что в указанной работе контрпримеры строятся на функциях синуса,
  про которые сказано, что они, разумеется, принадлежат пространству пробных функций. Но пробные функции для уравнений с оператором Бесселя всегда выбираются чётными. Поэтому, если считать синус чётной функцией, то верны результаты из [321--322], а если думать иначе, то мои.

8. \textbf{ ОП для дифференциально--разностных операторов}.

В последнее время теория  ОП для дифференциально--разностных операторов в основном разрабатывалась для так называемого оператора Дункла [323--333]. Этот оператор является в простейшем случае весовой суммой обычной производной и симметричной разности, а в общем многомерном случае он связан  с симметриями и группой Кокстера. При определённых условиях оператор Дункла является одной из неожиданных явных реализаций квадратного корня для дифференциального оператора Бесселя. Построены явные формулы для ОП типа Сонина и Пуассона и их обобщен
 ий, ООС, соответствующих рядов и обобщённых почти--периодических функций. На мой взгляд, часть этих результатов может быть выведена из результатов монографии А.~П.~Хромова [40], если использовать хитрое строение оператора Дункла. ОП для разностного уравнения Хилла изучаются в [334--335]. Сейчас исследование различных вопросов для оператора Дункла --- это чрезвычайно интересная тема, которой посвящены многочисленные работы. Эти результаты относятся к общей схеме применения ОП для того случая, когда  сплетаются дифференциальный и дифференциа
 льно--разностный операторы, а сам ОП является интегральным или интегро--дифференциальным оператором.

9. \textbf{ ОП и свёртки}.

Операторы дробного интегродифференцирования, обобщённого сдвига (ООС), различные свёртки --- все эти конструкции прямо или опосредовано связаны с теорией ОП. Связи с дробным  интегродифференцированием рассмотрены выше, ООС напрямую выражаются через ОП по явным формулам и строятся по одной с ними схеме. Под свёрткой в различных разделах математики понимаются похожие, но несколько различающиеся конструкции. Классическая свёртка для различных интегральных преобразований изучалась во многих работах, наиболее общая теория построена в
 аботах В.А.\,Какичева и его ученицы Л.Е.\,Бритвиной.

Несколько другая свёртка вводится для дифференциальных операторов. Эта теория развита в работах И.\,Димовски и его учеников, в том числе Н.\,Божинова. Среди результатов получена так называемая свёртка Димовски для второй производной, которая служит основой для построения специального операционного исчисления. В этой теории используются ОП, которые позволяют, зная свёртку для одного дифференциального оператора и ОП, который сплетает его со вторым  дифференциальным оператором, сразу построить и свёртку для второго дифференциального
 ператора. Свёртки типа Димовски --- это интересный раздел теории, который и уже имеет, и, по моему мнению,  будет иметь в дальнейшем многочисленные важные приложения.

В заключение просто отметим, что в монографиях Р.~Кэррола [4--6] подробно изложены приложения теории ОП в задачах  теории вероятностей и случайных процессов, линейном стохастическом оценивании, фильтрации, стохастических случайных уравнениях, обратных задачах геофизики и трансзвуковой газодинамики.
\newpage

\section{12. Избранные задачи.}

В заключение приведём несколько задач, решение которых по мнению автора было бы существенным для теории операторов преобразования.

Задача 1. В теории ОП есть своя теорема Ферма. Это задача о построении многомерных ОП, сплетающих стационарные операторы Шрёдингера и Лапласа $T(\Delta-q(x))=\Delta T$.
Есть оптимисты, верящие, что эта задача решена, при этом они обычно ссылаются на работу Л.~Д.~Фаддева [33]. Но чтение этой работы лучше начать с конца, где честно написано, что '... исследование должно быть проведено для того,
чтобы сделать строгими формальные рассуждения этой главы.
Имеющиеся в нашем распоряжении его варианты слишком
громоздки для того, чтобы их можно было помещать в
настоящем обзоре. Мы надеемся, что приведенная здесь
формальная схема решения многомерной обратной задачи
рассеяния может стать стимулом для некоторых читателей,
которые разработают более адекватные аналитические модели
для её оправдания'. Насколько мне известно, такой одарённый читатель пока к сожалению так и не появился. Думаю, что и в принципе перспектива построения таких ОП в обозримом виде представляется туманной. Ведь их ядра должны удовлетворять уже ультрагиперболическим уравнениям, методы исследований которых пока практически не разработаны. Думаю, что и множество решений уравнения Шрёдингера уж слишком разнородно и многообразно, чтобы его можно было описать по существу единой формулой по аналогии с уравнениями Штурма--Лиувилля.

Но самое обидное, что существование таких ОП  доказано под именем волновых операторов для всех разумных потенциалов! Не удаётся только построить их в обозримом виде, например, как многомерные операторы Вольтерра.

Периодически появляются отдельные публикации, посвящённые решению этой задачи в общем или частном случаях, автор затрудняется дать им адекватную оценку.

Задача 2. Обобщить результаты М.~М.~Маламуда из пункта 4, заменив обычные производные степенями оператора Бесселя.

Задача 3. Исторически сложилось, что изучение ОП началось с операторов второго порядка. Построить теорию ОП для операторов первого порядка, которая оказалась пропущена. Некоторые построения сделаны в [41] для ОП вида $T(D-q)=DT$.

Задача 4.  Можно ли с помощью ОП вида $T(D \pm\lambda)=DT$ сконструировать ОП вида  $T(D^2-\lambda^2)=D^2 T$ ? Тогда удалось бы построить, в частности, теорию ОП Векуа--Эрдейи--Лаундеса из более простых 'кирпичиков'. (Эта задача сложнее, чем кажется с первого взгляда. Возможно, для её решения потребуется использование гиперкомплексных чисел).

Задача 5. Применить к оценкам ядер ОП обобщения неравенств Коши--Буняковского--Шварца и Янга, метод построения которых разработан автором в [336--350].

Задача 6. Обобщить построения теории ОП, заменив обычные производные их $q$--аналогами. Упомянутые выше операторы Дункла---это лишь первый небольшой шаг в данном направлении.

Задача 7. Изучить основные задачи для уравнения с дробным оператором Бесселя (\ref{611}). Для этого изучить асимптотику и распределение нулей ядра резольвентного оператора (\ref{900}).

Задача 8. Обобщить метод операторов преобразования на простейшие уравнения, заданные на графах. При этом в качестве исходных могут быть использованы результаты А.М.\,Боровских, В.Л.\,Прядиева, О.М.\,Пенкина.

Задача 9. Обобщить метод операторов преобразования на простейшие уравнения, заданные на стратифицированных множествах, составленных из подмножеств разной размерности.  При этом в качестве исходных могут быть использованы результаты  О.М.\,Пенкина.

Задача 10. Построить ОП для дифференциальных уравнений бесконечного порядка, хотя бы модельных с постоянными коэффициентами. (Эта задача возникла у автора после прекрасной лекции проф. Юрия Фёдоровича Коробейника летом 2010 года на конференции во Владикавказе. Отметим, что за несколько дней до этой лекции Юрий Фёдорович отпраздновал свой 80--летний юбилей).

\section{ Благодарности.}

Автор многим обязан в студенческие годы своему  научному руководителю, а впоследствии  и  другу Валерию Вячеславовичу Катрахову, благодаря которому я начал заниматься теорией операторов преобразования и вообще математикой. Я горжусь тем, что принадлежу к числу учеников Ивана Александровича Киприянова, ветерана Великой Отечественной войны, основателя Воронежской школы теории уравнений в частных производных с особенностями.
Значительную роль сыграло обсуждение различных вопросов, и личное, и письменное с Анатолием Александровичем Килбасом. К глубокому сожалению, этих  прекрасных людей и замечательных математиков уже нет в живых.

Автор искренне благодарит:\\

Анатолия Георгиевича Кусраева за предложение написать этот обзор в 2007 году и огромное терпение, проявленное за время его подготовки;\\

Леонида Аркадьевича Минина, подробные обсуждения с которым помогали и помогают мне глубже разбираться в различных математических задачах;\\

Дмитрия Борисовича Карпа, с которым, не видя друг друга уже много лет после моего отъезда из Владивостока, мы практически  ежедневно обсуждаем нашу работу и являемся соавторами статей по специальным функциям;\\

участников в настоящем и прошлом Воронежского семинара по Анализу, теории функций и дифференциальным уравнениям (МНРС), на заседаниях которого неоднократно докладывались и обсуждались результаты по теории операторов преобразования, а именно  А.~В.~Боровских, А.~В.~Лободу, В.~В.~Меньших, И.~Я.~Новикова, О.~М.~Пенкина, И.~П.~Половинкина, В.~Л.~Прядиева,  В.~А.~Родина, Н.~Н.~Удоденко;\\

автор искренне благодарит за различную помощь, дружеское участие, полезные обсуждения, а также предоставленные монографии, статьи и другие материалы:\\

Роберта Вэйна Кэррола (США),  Ивана Димовски и Виржинию Кирякову (Болгария), Алоиса Куфнера (Чехия),
Халдуна Оцактаса (Турция); Ишхана Гвидоновича Хачатряна (Армения), Армена Джрбашяна (Армения),  Николая Иосифовича Юрчука (Белоруссия), Нину Афанасьевну Вирченко (Украина), Елену Александровну Аршаву (Украина), Ирину Ивановну Мороз (Украина), Виктора Гайдея (Украина),  Ростислава Олеговича Гринив (Украина), Оксану Валерьевну Скоромник (Белоруссия), Вагифа Гулиева (Азербайджан). \\

Виктора Захаровича Мешкова, Юрия Фёдоровича Коробейника, Адама Маремовича Нахушева, Владимира Павловича Глушко, Августа Петровича Хромова, Александра Павловича Солдатова, Евгения Александровича Радкевича, Сергея Борисовича Климентова, Арсена Владимировича Псху, Сергея Сергеевича Платонова, Александра Васильевича Глушака, Владимира Павловича Лексина, Олега Александровича Репина, Вячеслава Анатольевича Юрко, Александра Назарова,    Любовь Бритвину,  Геннадия Ляховецкого,     Владимира Воловича.

\newpage
 \Lit

 \begin{enumerate}
 \itemsep=0pt\parskip=0pt

\bib{Хорн~Р., Джонсон~Ч.}{Матричный анализ.---М.:~Мир,~1989.---655~с.}
\bib{Гантмахер~Ф.~Р.}{Теория матриц.---М.:~Наука,~1988.---552~с.}
\bib{Тыртышников~Е.~Е.}{Матричный анализ и линейная алгебра.---М.:~Физматлит,~2007.---480~с.}

\bib{Carroll~R.}{Transmutation and Operator Differential
Equations.---North Holland,~1979.---245~p.}
\bib{Carroll~R.}{Transmutation, Scattering Theory and Special Functions.---North Holland,~1982.---457~p.}
\bib{Carroll~R.}{Transmutation Theory and
Applications.---North Holland,~1986.---351~p.}
\bib{Gilbert~R., Begehr~H.}{Transformations, Transmutations and Kernel Functions. Vol. 1--2.---Longman,
Pitman,~1992.}
\bib{Trimeche~Kh.}{Transmutation Operators and Mean-Periodic Functions Associated with Differential Operators (Mathematical Reports, Vol 4, Part 1).---Harwood Academic Publishers,~1988.---282~p.}

\bib{Фаге~Д.~К., Нагнибида~Н.~И.}{Проблема эквивалентности обыкновенных
дифференциальных операторов.---Новосибирск:~Наука,~1977.---280~с.}
\bib{Gilbert~R.}{Constructive Methods for Elliptic Equations.~Springer Lecture
Notes Math,~365,~1974.}
\bib{Carroll~R.~W., Showalter R.E.}{Singular and Degenerate Cauchy problems.---N.Y.:~Academic Press,~1976.---333~p.}
\bib{Colton~D.}{Solution of Boundary Value Problems by the Method of Integral Operators.---London:~Pitman Press,~1976.}
\bib{Colton~D.}{Analytic Theory of Partial Differential Equations.---London:~Pitman Press,~1980.}
\bib{Gilbert~R.}{Function Theoretic Methods in Partial Differential Equations.---N.Y.:~Academic Press,~1969.}
\bib{Carroll~R.}{Topics in Soliton Theory.---North Holland,~1991.---428~p.}
\bib{Antimirov~M.~Ya., Kolyshkin~A.~A., Vaillancourt~R.}{Applied Integral Transforms.---CRM Monograph Series, AMS,~1993.---265~p.}
\bib{Lions~J.~L.}{Equations differentielles operationnelles et probl\'emes aux limites.---Springer,~1961.}
\bib{Trimeche~Kh.}{Generalized Harmonic Analysis and Wavelet Packets (an Elementary Treatment of Theory and Applications).---Taylor \& Francis,~2001.---320~p.}
\bib{Kiryakova~V.}{Generalized Fractional Calculus and Applications.--- Pitman Research Notes in Math. Series No. 301.---Longman Sci. UK.---1994.---402~p.}
\bib{Edited by Connett~W.C., Gebuhrer~M.--O., Schwartz~F.L.}{Applications of Hypergroups and Related Measure Algebras.---1995.---AMS Contemporary Mathematics.---No.~83.}

\bib{Агранович~З.~С., Марченко~В.~А.}{Обратная задача теории рассеяния.---Харьков: изд. ХГУ,~1960.---268~с.}
\bib{Марченко~В.~А.}{Спектральная теория операторов
Штурма--Лиувилля.---Киев: Наукова Думка,~1972.---220~с.}
\bib{Марченко~В.~А.}{Операторы
Штурма -- Лиувилля и их приложения.---Киев: Наукова Думка,~1977.---331~с.}
\bib{Левитан~Б.~М.}{Разложение по собственным функциям дифференциальных уравнений второго порядка.---М.:~Гостехиздат,~1950.}
\bib{Левитан~Б.~М.}{Операторы обобщённого сдвига и некоторые их применения.---М.: ГИФМЛ,~1962.---324~с.}
\bib{Левитан~Б.~М.}{Теория операторов обобщённого сдвига.---М.: Наука, 1973.---312~с.}
\bib{Левитан~Б.~М.}{Почти--периодические функции.---М.:~ГИТТЛ,~1953.---398~с.}
\bib{Левитан~Б.~М.}{Обратные задачи Штурма--Лиувилля.---М.: Наука,~1984.---240~с.}
\bib{Левитан~Б.~М., Саргсян~И.~С.}{Операторы Штурма--Лиувилля и Дирака.---М.: Наука,~1988.---432~с.}
\bib{Киприянов~И.~А.}{Сингулярные эллиптические краевые задачи.---М.: Наука-Физматлит,~1997.---204~с.}
\bib{Шадан~К., Сабатье~П.}{Обратные задачи в квантовой теории рассеяния.---М.: Мир,~1980.---408~с.}
\bib{Фаддеев~Л.~Д.}{Обратная задача квантовой теории рассеяния-1~/\!/~УМН.---1959.---Т.~14.---№~4.---С.~57--119.}
\bib{Фаддеев~Л.~Д.}{Обратная задача квантовой теории рассеяния-2~/\!/~"Итоги науки и техники", "Современные проблемы математики. т.~3".---ВИНИТИ.---1974.---С.~93--180.}
\bib{Килбас~А.~А., Маричев~О.~И., Самко~С.~Г.}{Интегралы и производные дробного порядка и некоторые
их приложения.---Минск:~Наука и техника,~1987.---688~с.}
\bib{Нижник~Л.~П.}{Обратная нестационарная задача теории рассеяния.---Киев: Наукова Думка,~1973.---183~с.}
\bib{Кошляков~Н.~С., Глинер~Э.~Б., Смирнов~М.~М.}{Уравнения в частных производных математической физики.---М.:~Высшая школа,~1962 (второе издание 1970 г.).---712~с.}
\bib{Положий~Г.~Н.}{Уравнения математической физики.---М.:Высшая школа,~1964.---560~с.}
\bib{Наймарк~М.~А.}{Линейные дифференциальные операторы.---М.:~Наука,~1969.---528~с.}
\bib{Левин~Б.~Я.}{Распределение корней целых функций.---М.:~ГИТТЛ,~1956.---632~с.}
\bib{Хромов~А.~П.}{Конечномерные возмущения вольтерровых операторов~/\!/~Современная математика. Фундаментальные направления---2004.---№~10.---С.~3--163.}
\bib{Ситник~С.~М.}{Элементарная теория операторов преобразования.(Рукопись подготовлена к печати).---150~с.}

\bib{Левитан~Б.~М., Повзнер~А.~Я.}{Дифференциальные уравнения  Штурма--Лиувилля на полуоси и теорема Планшереля~/\!/~Докл. АН СССР.---1946.---Т.~52.---№~6.---С.~483--486.}
\bib{Повзнер~А.~Я.}{О дифференциальных уравнениях типа Штурма--Лиувилля на полуоси~/\!/~Матем. сборник.---1948.---Т.~23~(65).---№~1.---С.~3--52.}
\bib{Блох~А.~Ш.}{Оп определении дифференциального оператора по его спектральной матрице--функции~/\!/~Докл. АН СССР.---1953.---Т.~92.---№~2.---С.~209--212.}
\bib{Левин~Б.~Я.}{Преобразования типа Фурье и Лапласа при помощи решений дифференциального уравнения второго порядка~/\!/~Докл. АН СССР.---1956.---Т.~106.---№~2.---С.~187--190.}

\bib{Соболев~С.~Л.}{Уравнения математической физики.---М.:~Наука,~1992.---432~с.}
\bib{Трикоми~Ф.}{Лекции по уравнениям в частных производных.---М.:~ИЛ,~1957.---444~с.}
\bib{Комеч~А.~И.}{Практическое решение уравнений математической физики.---М.:~МГУ,~1993.---155~с.}
\bib{Copson~E.~T.}{On the Riemann--Green Function~/\!/~Archive for Rational Mechanics and Analysis.---1957/58.---V.1.---P.~324--348.}
\bib{Ибрагимов~Н.~X.}{Опыт группового анализа обыкновенных дифференциальных
уравнений.---Серия~"Математика. Кибернетика".---М.:~Знание,~1991.---47~с.}
\bib{Волкодавов~В.~Ф., Лернер~М.~Е., Николаев~Н.~Я., Носов~В.~А.}{Таблицы некоторых функций Римана, интегралов и рядов.---Куйбышев:~изд. Куйбышев. гос. пед. инст.,~1982.---56~с.}
\bib{Волкодавов~В.~Ф.,  Николаев~Н.~Я.}{Интегральные уравнения Вольтерра первого рода с некоторыми специальными функциями в ядрах и их приложения.---Самара:~изд.-во 'Самарский университет',~1992.---100~с.}
\bib{Волкодавов~В.~Ф., Захаров~В.~Н.}{Таблицы функций Римана и Римана--Адамара для некоторых дифференциальных уравнений в n--мерных евклидовых пространствах.----Самара,~1994.---31~с.}
\bib{Волкодавов~В.~Ф., Николаев~Н.~Я., Быстрова~О.~К., Захаров~В.~Н.}{Функции Римана для некоторых дифференциальных уравнений в n--мерных евклидовых пространствах и их приложения.---Самара:~изд.-во 'Самарский университет',~1995.---76~с.}
\bib{Волкодавов~В.~Ф., Захаров~В.~Н.}{Функции Римана для одного класса дифференциальных уравнений в трёхмерном эвклидовом пространстве и её применения.--- – Самара:~изд.-во 'Самарский университет', 1996.---53~с.}
\bib{Пулькин~С.~П.}{Некоторые краевые задачи...~/\!/~Уч.зап. Куйбышевского
пед. ин-та.---1958.---Вып.~21.--- С.~3--54.}
\bib{Лернер~М.~Е.}{Принципы максимума для уравнений гиперболического типа и новые свойства функции Римана.---Самара:~Самарский государственный технический университет,~2001.---113~с.}
\bib{Малышев~Ю.~В.}{Решение линейных дифференциальных уравнений с переменными коэффициентами и сложными функциями~/\!/~Вестник Сам. ГТУ., серия 'физико--математические науки'---2002.---Вып.~16.---С.~5--9.}
\bib{Малышев~Ю.~В.}{Факторизация дифференциальных уравнений в частных производных~/\!/~Дифференциальные уравнения и процессы управления.---2003.---№~2.---С.~78--86}
\bib{Жегалов~В.~И., Миронов~А.~Н.}{Дифференциальные уравнения со старшими частными производными.---Казань:~Казанское математическое общество,~2001.---226~с.}
\bib{Аксёнов~А.~В.}{Метод построения функции Римана гиперболического уравнения второго порядка~/\!/~В сб.: Тезисы докладов международной конференции 'Дифференциальные уравнения, теория функций и приложения', посвящённой 100--летию со дня рождения академика И.Н.Векуа.---Новосибирск:~НГУ,~2007.---С.~47--48.}
\bib{Аксёнов~А.~В.}{Симметрии линейных уравнений с частными производными и фундаментальные решения~/\!/~Доклады АН России.---1995.---Т.~342.---№~2.---С.~151--153.}
\bib{Кощеева~О.~А.}{О построении функции Римана для уравнения Бианки в четырёхмерном пространстве~/\!/~В сб.:~Труды участников международной школы--семинара по геометрии и анализу памяти Н.В.Ефимова, Абрау--Дюрсо, "Лиманчик". (отв. ред. С.~Б.~Климентов)---Ростов--на--Дону:~2006.---С.~236--238.}
\bib{Романовский~Р.~К.}{О матрицах Римана первого и второго рода~/\!/~Матем. сборник.---1985.---Т.~121.---№.~4.---С.~494--501.}
\bib{Воробьёва~Е.~В., Романовский~Р.~К.}{Метод характеристик для гиперболических краевых задач на плоскости~/\!/~СМЖ.---2000.---Т.~41.---№~3.---С.~531--540.}
\bib{Энбом~Е.~А., Волкодавов~В.~Ф.}{Неклассическая задача для вырождающегося гиперболического
уравнения третьего порядка~/\!/~Известия ву-
зов. Математика. -Казань, 2003. -Деп. в ВИНИТИ 23.07.2003. №~1445-В2003.}

\bib{Колмогоров~А.~Н., Фомин~С.~В.}{Элементы теории функций и функционального анализа.---М.:~Наука,~1981.---544~с.}
\bib{Siersma~J.}{Thesis.---Groningen,~1979.}
\bib{Чернятин~В.~А.}{Обоснование метода Фурье в смешанной задаче для уравнений в частных производных.---М.:~МГУ,~1991.---112~с.}
\bib{Нейман-заде~М.~И.,  Шкаликов~А.~А.}{Операторы Шрёдингера с сингулярными потенциалами из пространств мультипликаторов~/\!/~Матем. заметки.---1999.---
Т.~66.---№~5.---С.~723–-733.}
\bib{Савчук~А.~М., Шкаликов~А.~А.}{Операторы Штурма–-Лиувилля с сингулярными потенциалами~/\!/~Матем. заметки.---1999.---Т.~66.---№~6.---С.~897–-912.}
\bib{Савчук~А.~М., Шкаликов~А.~А.}{Формула следа для операторов Штурма–Лиувилля с сингулярными потенциалами~/\!/~Матем. заметки.---2001.---Т.~69.---№~3.---С.~427–-442.}
\bib{Савчук~А.~М., Шкаликов~А.~А.}{О собственных значениях оператора Штурма–Лиувилля с потенциалами из пространств Соболева~/\!/~Матем. заметки.---2006.---Т.~80.---№~6.---С.~864–-884.}
\bib{Ширяев~Е.~А., Шкаликов~А.~А.}{Регулярные и вполне регулярные дифференциальные операторы~/\!/~Матем. заметки.---2007.---Т.~81.---№~4.---С.~636–-640.}
\bib{Boumenir~A., Vu Kim Tuan.}{Existence and construction of the transmutation operator~/\!/~J. Math. Phys.---2004.---Vol.~45.---Issue~7.---P.~2833--2843.}
\bib{Гусейнов~И.~М.}{Об одном операторе преобразования~/\!/~Матем. заметки.---1997.---Т.~62.---№~2.---С.~206--215.}
\bib{Гусейнов~И.~М., Набиев~А.~А., Пашаев~Р.~Т.}{Операторы преобразования и асимптотические формулы для собственных значений полиноминального пучка операторов Штурма--Лиувилля~/\!/~СМЖ.---2000.---Т.~41.---№~3.---С.~554--566.}

\bib{Векуа~И.~Н.}{О решениях уравнения $\Delta u +
\lambda^2 u$~/\!/~Сообщения Акад. Наук Груз. ССР.---1942.---Т.~III.---No~4.---C.~307--314.}
\bib{Векуа~И.~Н.}{Обращение одного интегрального
преобразования и его некоторые применения~/\!/~Сообщения Акад. Наук
Груз. ССР.---1945.---Т.~VI.---No~3.---C.~177--183.}
\bib{Векуа~И.~Н.}{Новые методы решения эллиптических уравнений.---
М.-Л.:~ГИТТЛ,~1948.---296~c.}
\bib{Векуа~И.~Н.}{Обобщённые аналитические функции.---М.:~Наука,~1988.---512~с.}
\bib{Erdelyi~A.}{Some applications of fractional integration~/\!/~
Boeing Sci. Res. Labor. Docum. Math. Note D1--82--0286.---1963.---No~316.---23~ p.}
\bib{Erdelyi~A.}{An application of fractional integrals~/\!/~J. Analyse Math.---1965.---Vol.~14.---P.~113--126.}
\bib{Erdelyi~A.}{On the Euler--Poisson--Darboux equation~/\!/~J.
Analyse Math.---1970.---Vol.~23.---P.~89--102.}
\bib{Lowndes~J.S.}{An
application of some fractional integrals~/\!/~Glasgow
Math. J.---1979.---Vol.~20.---No~1.---P.~35--41.}
\bib{Lowndes~J.S.}{On some generalizations of Riemann--Liouville
and Weil fractional integrals and their applications~/\!/~Glasgow
Math. J.---1981.---Vol.~22.---No~2.---P.~73--80.}
\bib{Lowndes~J.S.}{Cauchy problems for second order hyperbolic differential
equations with constant coefficients~/\!/~Proc. Edinburgh Math. Soc.---1983.---Vol.~26.---No~3.---P.~97--105.}
\bib{Бергман~С.}{Интегральные операторы в теории уравнений с частными производными.---М.:~Мир,~1964.---309~с.}
\bib{Bers~L., Gelbart~A.}{On a class of differential equations in mechanics of continua~/\!/~Quart. of Appl. Math.---1943.---V.~5.---No.~1.---P.~168--188.}
\bib{Bers~L., Gelbart~A.}{On a class of functions defined by partial differential equations~/\!/~Trans. Amer. Math. Soc.---1944.---V.~56.---P.~67--93.}
\bib{Bers~L.}{A remark on an applications of pseudo--analytic functions~/\!/~Amer. J. Math.---1956.---V.~78.---No.~3.---P.~486--496. }
\bib{Weinstein~A.}{Discontinuos integrals and generalized theory of potential~/\!/~Trans. Amer. Math. Soc.---1948.---V.~63.---No.~2.---P.~342--354.}
\bib{Weinstein~A.}{Generalized axially symmetric potential theory~/\!/~Bull. Amer. Math. Soc.---1953.---V.~59.---P.~20---38.}
\bib{Положий~Г.~Н.}{Обобщение теории аналитических функций комплексного переменного.---Киев:~изд. КГУ,~1965.---444~с.}
\bib{Климентов~С.~Б.}{Классы Харди обобщённых аналитических функций~/\!/~Изв. ВУЗов, Сев.--Кавказ. регион, серия 'Естественные науки'.---2003.---№~3.---С.~6--10.}
\bib{Климентов~С.~Б.}{Классы Смирнова обобщённых аналитических функций~/\!/~Изв. ВУЗов, Сев.--Кавказ. регион, серия 'Естественные науки'.---2005.---№~1.---С.~13--17.}
\bib{Климентов~С.~Б.}{Классы ВМО обобщённых аналитических функций~/\!/~Владикавказский математический журнал.---2006.---Т.~8.---В.~1.---С.~27--39.}
\bib{Климентов~С.~Б.}{Теорема двойственности для классов Харди обобщённых аналитических функций~/\!/~ В сб.:
Исследования по комплексному анализу, теории операторов
 и математическому моделированию.---Владикавказ:~
 Изд-во ВНЦ РАН,~2006.---С.~63--73.}
\bib{Ляховецкий~Г.~В., Ситник~С.~М.}{Операторы преобразования
Векуа--Эрдейи--Лаундеса~/\!/~Препринт института автоматики и
процессов управления Дальневосточного отделения РАН.---Владивосток:~ДВО РАН,
1994.---24~с.}
\bib{Lyahovetskii~G.~V., Sitnik~S.~M.}{The Vekua --
Erdelyi -- Lowndes transmutations~/\!/~Препринт института автоматики и
процессов управления Дальневосточного отделения РАН.---Владивосток:~ДВО РАН,
1994.---14~с.}
\bib{Ситник~С.~М.}{Построение операторов преобразования Векуа--Эрдейи--Лаундеса~/\!/~В сб.: Тезисы докладов международной конференции 'Дифференциальные уравнения, теория функций и приложения', посвящённой 100--летию со дня рождения академика И.~Н.~Векуа.---Новосибирск:~НГУ,~2007.---С.~469--470.}

\bib{Бейтмен~Г., Эрдейи~А.}{Высшие трансцендентные функции, т.1.---М.:~Наука, Гл. ред. ФМЛ,~1973.---296~с.}
\bib{Бейтмен~Г., Эрдейи~А.}{Высшие трансцендентные функции, т.2.---М.:~Наука, Гл. ред. ФМЛ,~1966.---296~с.}
\bib{Бейтмен~Г., Эрдейи~А.}{Высшие трансцендентные функции, т.3.---М.:~Наука, Гл. ред. ФМЛ,~1967.---299~с.}
\bib{Andrews~G.~E.,  Askey~R., Roy~R.}{Special Functions.---Cambridge University Press,~1999.---665~p.}
\bib{Luke~Y.~L.}{Mathematical functions and their approximations.---Academic Press,~1975.---584~p.}
\bib{Luke~Y.~L.}{The special functions and their approximations. Volume 1.---Academic Press,~1969.---362~p.}
\bib{Люк~Ю.}{Специальные математические функции и их аппроксимации.---М.:~Мир,~1980.---610~с.}
\bib{Под ред. Абрамовиц~М., Стиган~И.}{Справочник по специальным функциям.---М.:~Наука,~1979.---832~с.}
\bib{Уиттекер~Э., Ватсон~Дж. }{Курс современного анализа. Часть вторая.  Трансцендентные функции.---М.:~ГИФМЛ,~1963.---516~с.}
\bib{http://functions.wolfram.com}{}
\bib{http://dlmf.nist.gov/}{Digital Library of Mathematical Functions}
\bib{Slater~L.~J.}{Generalized hypergeometric functions.---Cambridge University Press,~1966.---285~p.}
\bib{Dwork~B.}{Generalized hypergeometric functions.---Oxford,~1990.---196~p.}
\bib{Dwork~B., Gerotto, Sullivan.}{Introduction to G-functions.---Princeton,~1994.---335~p.}
\bib{Bailey~W.~N.}{Generalized hypergeometric series.---Cambridge,~1964.---117~p.}
\bib{Karlsson~P.~W., Srivastava~H.~M.}{Multiple Gaussian hypergeometric series.---Ellis Horwood
Series:~Mathematics and its Applications, New York, 1985.}
\bib{Appell~P., Kampe de Feriet~J.}{Fonctions Hypergeometriques et Hyperspheriques; Polynomes d'Hermite.---Paris:~Gauthier-Villars,~1926.}
\bib{Exton~H.}{Multiple Hypergeometric Functions and Applications.---New York:~John Wiley and Sons.---1976.}
\bib{Kilbas~A.~A., Saigo~M.}{H - transforms. Theory and applications. Chapman and Hall,~CRC.---2004.}
\bib{Килбас~А.~А., Сайго~М., Жук~В.~А.}{~/\!/~Дифференциальные уравнения.---1991.---Т.~27.---№~9.---С.~1640--1642.}

\bib{Delsarte~J.}{Sur certaines transformation fonctionnelles relative aux \'{e}quations lin\'eares aux d\`eriv\'ees partielles du seconde ordre~/\!/~C.~R. Acad. Sci. Paris.---1938.---206.---1780--1782.}
\bib{Lions~J.~L.}{Op\'erateurs de Delsarte et probl\`eme mixte~/\!/~Bull. Soc. Math. France.---1956.---No.~84.---9--95.}
\bib{Lions~J.~L.}{Quelques applications  d'op\'erateurs de transmutations~/\!/~Colloques Internat. Nancy.---1956.---125--142.}
\bib{Delsarte~J., Lions~J.~L.}{Transmutations d'op\'erateurs diff\'erentiels dans le domaine complexe~/\!/~Comm. Math. Helv.---1957.---No.~32.---113--128.}
\bib{Delsarte~J., Lions~J.~L.}{Moyennes g\'en\'eralis\'ees~/\!/~Comm. math. Helv.---1959.---No.~34.---59--69.}
\bib{Арнольд~В.~И.}{Математические эпидемии 20 века~/\!/~Независимая газета.---2001.}
\bib{Фаге~М.~К.}{Построение операторов преобразования и решение одной проблемы моментов для обыкновенных линейных дифференциальных уравнений произвольного порядка~/\!/~УМН.---1957.---Т.~12.---Вып.~1~(73).---С.~240--245.}
\bib{Фаге~М.~К.}{Операторно--аналитические функции одной независимой переменной~/\!/~Докл. АН СССР.---1957.---Т.~112.---№~5.---С.~1008--1011.}
\bib{Фаге~М.~К.}{Интегральные представления операторно--аналитических функций одной независимой переменной~/\!/~Докл. АН СССР.---1957.---Т.~115.---№~5.---С.~874--877.}
\bib{Фаге~М.~К.}{Операторно--аналитические функции одной независимой переменной~/\!/~Тр. Моск. матем. об--ва.---1958.---Т.~7.---С.~227--268.}
\bib{Фаге~М.~К.}{Интегральные представления операторно--аналитических функций одной независимой переменной~/\!/~Тр. Моск. матем. об--ва.---1958.---Т.~8.---С.~3--48.}
\bib{Сахнович~Л.~А.}{Спектральный анализ вольтерровских операторов и обратные задачи~/\!/~Докл. АН СССР.---1957.---Т.~115.---№~4.---С.~666--669.}
\bib{Сахнович~Л.~А.}{Обратная задача для дифференциальных операторов порядка $n>2$ с аналитическими коэффициентами~/\!/~Матем. сборник.---1958.---Т.~46.---№~1.---С.~61--76.}
\bib{Мацаев~В.~И.}{О существовании оператора преобразования для дифференциальных уравнений высших порядков~/\!/~Докл. АН СССР.---1960.---Т.~130.---№~3.---С.~499--502.}
\bib{Сахнович~Л.~А.}{Необходимые условия наличия операторов преобразования для уравнения четвёртого порядка~/\!/~УМН.---1961.---Т.~16.---Вып.~5.---С.~199--205.}
\bib{Хачатрян~И.~Г.}{Об операторах преобразования для дифференциальных уравнений высших порядков~/\!/~Изв. АН Арм.ССР. Сер. Математика.---1978.---Т.~13.---№~3.---С.~215--236.}
\bib{Хачатрян~И.~Г.}{Об операторах преобразования для дифференциальных уравнений высших порядков, сохраняющих асимптотику решений~/\!/~Изв. АН Арм.ССР. Сер. Математика.---1979.---Т.~14.---№~6.---С.~424--445.}
\bib{Поляцкий~В.~Т.}{О свойствах решений некоторого уравнения~/\!/~УМН.---1965.---Т.~17.---№.~4.---С.~119--124.}
\bib{Маламуд~М.~М.}{Об операторах преобразования для обыкновенных дифференциальных уравнений высших порядков~/\!/~Матем. анализ и теория вероятн. Киев:~Наукова думка,~1978.---С.~108--111.}
\bib{Маламуд~М.~М.}{Необходимые условия существования оператора преобразования для уравнений высших порядков~/\!/~Функцион. анализ и его прил.---1982.---Т.~16.---№~3.---С.~74--75.}
\bib{Маламуд~М.~М.}{К вопросу об операторах преобразования~/\!/~Препринт ИМ АН УССР.---Киев,~1984.---48~с.}
\bib{Маламуд~М.~М.}{Операторы преобразования для  уравнений высших порядков~/\!/~Матем. физика и нелин. механика. Киев:~Наукова думка, 1986.---№~6.---С.~108--111.}
\bib{Маламуд~М.~М.}{К вопросу об операторах преобразования для обыкновенных дифференциальных уравнений~/\!/~Тр. Моск. матем. об--ва.---1990.---Т.~53.---С.~69--97.}
\bib{Morrey~C.~B.}{Multiple integrals in the calculus of variations.---Springer,~1966.}
\bib{Коробейник~Ю.~Ф.}{Операторы сдвига на числовых семействах.---Ростов--на--Дону:~ИРУ,~1983.---155~с.}
\bib{Коробейник~Ю.~Ф.}{О разрешимости в комплексной области некоторых общих классов линейных интегральных уравнений.---Ростов--на--Дону:~2005.---245~с.}
\bib{Фетисов~В.~Г.}{Операторы и уравнения в локально ограниченных пространствах~/\!/В книге:
Операторы и уравнения в линейных топологических пространствах.---Владикавказ:
 Изд-во ВНЦ РАН,~2006.---С.~7--142.}
\bib{Кусраев~А.~Г.}{Мажорируемые операторы.---М.:~Наука,~2003.---619~с.}
\bib{Пасенчук~А.~Э.}{Абстрактные сингулярные операторы.---Новочеркасск,~1993. }
\bib{Фишман~М.~К.}{Матем. сборник.---Т.~68~(110).---№~1.---С.~73.}
\bib{Марченко~В.~А.}{Некоторые вопросы теории дифференциального оператора второго порядка~/\!/~Докл. АН СССР.---1950.---Т.~72.---№~3.---С.~457--460.}
\bib{Марченко~В.~А.}{Операторы преобразования~/\!/~Докл. АН СССР.---1950.---Т.~74.---№~2.---С.~185--188.}
\bib{Марченко~В.~А.}{О формулах обращения, порождаемых линейным дифференциальным оператором второго порядка~/\!/~Докл. АН СССР.---1950.---Т.~74.---№~4.---С.~657--660.}
\bib{Марченко~В.~А.}{Некоторые вопросы теории одномерных дифференциальных операторов второго порядка, I~/\!/~Тр. Моск. матем. об--ва.---1952.---I.---С.~327--420.}
\bib{Марченко~В.~А.}{Некоторые вопросы теории одномерных дифференциальных операторов второго порядка, II~/\!/~Тр. Моск. матем. об--ва.---1953.---II.---С.~3--82.}
\bib{Cheikh~B.}{Relations between harmonic analysis associated with two differential operators of different orders~/\!/~Journal of Computational and Applied Mathematics.---2003.---Vol.~153.---No.~1.---P.~61--71.}
\bib{Jafford~K.}{Inversion of the Lions transmutation operators using generalized wavelets~/\!/~Applied and Computational Harmonic Analysis.---1997.---Vol.~4.---No.~1.---P.~97--112.}
\bib{Samoilenko~A.M., Prykarpatsky~Ya.~A., Prykarpatsky~A.~K.}{The generalized de Rham--Hodge theory aspects of Delsarte-Darboux type transformations in multidimension~/\!/~Central European Journal of Mathematics.---2005.---Vol.~3.---No.~3.---P.~529--557.}
\bib{Golenia~J., Samoilenko~A.M., Prykarpatsky~Ya.~A., Prykarpatsky~A.~K.}{The general differential-geometric structure of
multidimensional Delsarte transmutation
operators in parametric functional spaces and
their applications in soliton theory~/\!/~arXiv:~math-ph/0404016.---2004.---10~p.}
\bib{Trimeche~Kh.}{Inversion of the Lions transmutation operators using generalized wavelets~/\!/~Applied and Computational Harmonic Analysis.---1997.---Vol.~4.---No.~1.---P.~97-–112, }
\bib{Kiryakova~V.~S.}{An explanation of Stokes phenomenon by the method of transmutations~/\!/~Proc. Conf. Diff. Equations and Appl., Rousse, 1982.---P.~349--353.}
\bib{Carroll~R., Boumenir~A.}{Toward a general theory of
transmutation~/\!/~arXiv:~funct-an/9501006.---1995.---19~p.}

\bib{Юшкевич~А.~П.}{Леонард Эйлер.---М.:~Знание,~1982.---64~с.}
\bib{Ватсон~Г.~Н.}{Теория бесселевых функций,~т.~1.---М.:~ИЛ,~1949.---728~с.}
\bib{Ситник~С.~М.}{Леонард Эйлер и теория специальных функций~/\!/~В сб.:~Материалы Международной научной конференции
'Леонард Эйлер и современная наука'.---Санкт-Петербург,~2007.---С.~ 192-200.}
\bib{Сонин~Н.~Я.}{Исследования о цилиндрических  функциях и специальных полиномах.---М.:~Гостехтеоретиздат,~1954.---241~с.}
\bib{Delsarte~J.}{Sur une extension de la formule de Taylor~/\!/~Journ. Math. pures et appl.---1938.---17.---217--230.}
\bib{Delsarte~J.}{Une extension nouvelle de la th\'eory de fonction presque p\'eriodiques de Bohr~/\!/~Acta Math.---1939.---69.---259--317.}
\bib{Delsarte~J.}{Hypergroupes et operateurs de permutation et de transmutation~/\!/~Colloques Internat. Nancy.---1956.---29--44.}
\bib{Левитан~Б.~М.}{Разложения по функциям Бесселя в ряды и интегралы Фурье~/\!/~УМН.---1951.---Т.~6.---Вып.~2.---С.~102--143.}
\bib{Платонов~С.~С.}{Гармонический анализ Бесселя и приближение функций на полупрямой~/\!/~Изв. РАН. Сер. матем.---2007.---Т.~71.---№~5.---C.~149–-196.}
\bib{Платонов~С.~С.}{Обобщённые сдвиги Бесселя и некоторые задачи теории приближения функций в метрике $L_2$--1~/\!/~Труды ПетрГУ.---Сер. Математика.---2000.---Вып.~7.---С.~70--82.}
\bib{Платонов~С.~С.}{Обобщённые сдвиги Бесселя и некоторые задачи теории приближения функций в метрике $L_2$--2~/\!/~Труды ПетрГУ.---Сер. Математика.---2001.---Вып.~8.---С.~20--36.}
\bib{Киприянов~И.~А.}{Преобразования Фурье--Бесселя и теоремы вложения для весовых классов~/\!/~Труды матем. ин--та АН СССР им. В.А.Стеклова.---1967.---Т.~89.---С.~130--213.}
\bib{Трибель~Х.}{Теория интерполяции, функциональные пространства, дифференциальные операторы.---М.:~Мир,~1980.---664~с. }
\bib{Кудрявцев~Л.~Д., Никольский~С.~М.}{Пространства дифференцируемых функций многих переменных и теоремы вложения~/\!/~ВИНИТИ. Итоги науки и техники. Современные проблемы математики. Фундаментальные направления.---М.,~1988.---С.~5--157.}
\bib{Келдыш~М.~В.}{О некоторых случаях вырождения уравнений эллиптического типа на границе области~/\!/~Докл. АН СССР.---1951.---Т.~77.---№~1.---С.~181--183.}
\bib{Ключанцев~М.~И.}{О построении $r$--чётных решений сингулярных дифференциальных уравнений~/\!/~Докл. АН СССР.---1975.---Т.~224.---№~5.---С.~1004--1007.}
\bib{Ключанцев~М.~И.}{Интегралы дробного порядка и сингулярные краевые задачи~/\!/~Дифференц. уравнения.---1976.---Т.~12.---№~6.---С.~983--990.}
\bib{Rodrigues~J.}{Operational calculus for the generalized Bessel operator~/\!/~SERDICA, Bulgaricae mathematicae publicationes.---1989.---Vol.~15.---P.~179--186.}
\bib{Dimovski~I.~H., Kiryakova~V.S.}{Generalized Poisson transmutations and corresponding representations of hyper--Bessel functions~/\!/~C. R. de l'Acad. bulgare des Sci.---1986.---T.~39.---№~10.---P.~29--32.}
\bib{Kamoun~L., Sifi~M.}{Bessel-Struve intertwining operator and generalized Taylor series on the real line~/\!/~Integral Transforms and Special Functions.---2005.---Vol.~16.---No.~1.---P.~39--55.}
\bib{Ali~I., Kiryakova~V., Kalla~S.L.}{Solutions of fractional multi-order integral and differential equations using a Poisson-type transform~/\!/~Journal of Mathematical Analysis and Applications.---2002.---Vol.~269.---No.~1.---P.~172--199.}
\bib{Rachdi~L.~T.}{Fractional powers of Bessel operator and inversion formulas for Riemann-Liouville and Weyl transforms~/\!/~Journal of Mathematical Sciences.---2001.---Vol.~12.---No.~1.}
\bib{Dimovski~I.~H., Kiryakova~V.S.}{Transmutations, convolutions and fractional powers of Bessel-type operators via Meijer G-functions.---Proc. Conf. Complex Anal. and Appl., Varna'~1983.---Sofia,~1985.---P.~45--66.}
\bib{Kiryakova~V.S.}{Applications of the generalized Poisson transformation for solving hyper-Bessel differential equations (In Bulgarian)~/\!/~Godishnik VUZ. Appl. Math.---1986.---Vol.~22.---No~4.---P.~129--140.}
\bib{Baccar~C., Hamadi~N.~B., Rachdi~L.~T.}{Inversion formulas for riemann-liouville transform and its dual associated with
singular partial differential operators~/\!/~International Journal of Mathematics and Mathematical Sciences.---2006.---P.~1–-26.}
\bib{Ярославцева~В.~Я.}{Об одном классе операторов преобразования и их приложении к дифференциальным уравнениям~/\!/~Докл. АН СССР.---1976.---Т.~227.---№~4.---С.~816--819.}
\bib{Ярославцева~В.~Я.}{Неоднородная граничная задача в полупространстве для одного класса сингулярных уравнений~/\!/~Деп. ред. "Дифференц. уравнения".---Минск,~1989.---11~с.}
\bib{Хелгасон~С.}{Группы и геометрический анализ.---М.:~Мир, 1987.---736~с.}

\bib{Kilbas~A.~A., Srivastava~H.~M., Truhillo~J.~J.}{Theory and Applications of Fractional Differential Equations.---North Holland Mathematical Studies.---Vol.~204.--- Elsevier,~2006.---523~p.}
\bib{Нахушев~А.~М.}{Уравнения математической биологии.---М.:~Высшая Школа, 1995.---301~с.}
\bib{Нахушев~А.~М.}{Элементы дробного исчисления и их применение.---Нальчик,~2000.---300~с.}
\bib{Нахушев~А.~М.}{Дробное исчисление и его применение.---М.:~Физматлит,~2003.---273~с.}
\bib{Псху~А.~В.}{Краевые задачи для дифференциальных уравнений с частными производными дробного и континуального порядка.---Нальчик,~2005.---186~с.}
\bib{Kiryakova~V.}{All the special functions are fractional differintegrals of elementary functions~/\!/~J. Physics A: Math. \& General.---1997.---Vol.~30.---No.~14.---P.~5085--5103.}
\bib{Ситник~С.~М.}{О некоторых обобщениях дробного  интегродифференцирования~/\!/ В сб.:~Материалы международного симпозиума "Уравнения смешанного типа и родственные проблемы анализа и информатики".---Нальчик-Эльбрус,~2003.--- С.~86-87.}
\bib{Ситник~С.~М.}{ Дробное интегродифференцирование для дифференциального оператора Бесселя~/\!/ В сб.:~Материалы международного Российско-Казахского симпозиума "Уравнения смешанного типа и родственные проблемы анализа и информатики".---Нальчик-Эльбрус,~2004.---С.~163-167.}
\bib{Ситник~С.~М.}{Об обобщении формулы Хилле-Тамаркина для резольвенты на случай операторов дробного интегрирования Бесселя~/\!/ В сб.:~
III Международная конференция:
" Нелокальные краевые задачи и родственные проблемы
математической биологии, информатики и физики".---Нальчик,~2006.---С.~269-270.}
\bib{Ситник~С.~М.}{Операторы дробного интегро-дифференцирования для дифференциального оператора Бесселя~/\!/ В сб.:~Труды четвёртой Всероссийской научной конференции
с международным участием:
"Математическое моделирование и краевые задачи". Часть 3.---Самара,~2007.---С.~158-160.}
\bib{Коновалова~Д.~С., Ситник~С.~М.}{Формула Тэйлора для операторов типа Бесселя~/\!/ В сб.:~Тезисы докладов Воронежской весенней математической школы 'Современные методы в теории краевых задач. Понтрягинские чтения --- VII'.---Воронеж,~1996.---С.~102.}
\bib{Катрахов~В.~В., Катрахова~А.~А.}{Формула Тэйлора с оператором
Бесселя для функций одной и двух переменных~/\!/~Деп. ВИНИТИ.---Воронеж,~
1982.---23~с.}
\bib{Sprinkhuizen-Kuyper~I.~G.}{A fractional integral operator corresponding to negative powers of a certain second-order differential operator~/\!/~J. Math.Analysis and Applications.---1979.---No.~72.---P.~674-702.}
\bib{Репин~О.~А.}{Краевые задачи со смещением для уравнений
гиперболического и смешанного типов.---Самара,~1992.}
\bib{Hille~E., Tamarkin~J.~D.}{On the theory of linear integral equations~/\!/~Ann. Math.---1930.---Vol.~31.---P.~479-528.}
\bib{Джрбашян~М.~М.}{Интегральные преобразования и представления функций в комплексной области.---М.:~Наука.---1966.---672~с.}
\bib{Кочубей~А.~Н.}{Задача Коши для эволюционных уравнений дробного порядка~/\!/~Дифференциальные уравнения.---1989.---Т.~25.---№~8.---С.~1359--1369.}
\bib{Кочубей~А.~Н.}{Диффузия дробного порядка~/\!/~Дифференциальные уравнения.---1990.---Т.~26.---№~4.---С.~660--770.}
\bib{Kilbas~A.~A., Srivastava~H.~M., Trujillo~J.~J.}{Theory and applications    of fractional differential equations.---North--Holland Mathematics Studies, Vol.~204, Elsever:~Amsterdam,~2006.}
\bib{Kilbas~A.~A., Trujillo~J.~J.}{Differential equations of fractional order: methods, results and problems. Part I~/\!/~Journal of Applicable Analysis.---2001.---Vol.~78.---Nos.~1--2.---P.~153--192.}
\bib{Kilbas~A.~A., Trujillo~J.~J.}{Differential equations of fractional order: methods, results and problems. Part II~/\!/~Journal of Applicable Analysis.---2001.---Vol.~81.---No.~2.---P.~435--493.}
\bib{Глушак~А.~В.}{Задача типа Коши для абстрактного дифференциального уравнения с дробными производными~/\!/~Матем. заметки.---2005.---Т.~77.---№~1.---С.~28--41.}
\bib{Килбас~А.~А., Репин~О.~А.}{ Аналог задачи Бицадзе--Самарского для уравнения смешанного типа с дробной производной~/\!/~Дифференциальные уравнения.---2003.---Т.~39.---№~5.---С.~638--644.}
\bib{Андреев~А.~А., Огородников~Е.~Н.}{Некоторые локальные и нелокальные аналоги задачи Коши-Гурса для системы уравнений типа Бицадзе-Лыкова с инволютивной матрицей~/\!/~Вестник Самарского государственного  технического  университета.  Сер. "Физико-математические науки".---2002.--- Вып.~16.---С.~19--35.}

\bib{Buschman~R.~G.}{An inversion integral for a general Legendre transformation~/\!/~SIAM Review.---1963.---Vol.~5.---No.~3.---P.~232--233.}
\bib{Buschman~R.~G.}{An inversion integral for a Legendre transformation~/\!/~Amer. Math. Mon.---1962.---Vol.~69.---No.~4.---P.~288--289.}
\bib{Erdelyi~A.}{An integral equation involving Legendre functions~/\!/~SIAM Review.---1964.---Vol.~12.---No.~1.---P.~15--30.}
\bib{Erdelyi~A.}{Some integral equations involving finite parts of divergent integrals~/\!/~Glasgow Math. J.---1967.---Vol.~8.---No.~1.---P.~50--54.}
\bib{Higgins~T.~P.}{A hypergeometric function transform~/\!/~SIAM J.---1964.---Vol.~12.---No.~3.---P.~601--612.}
\bib{Ta Li}{A new class of integral transform~/\!/~Proc. AMS.---1960.---Vol.~11.---No.~2.---P.~290--298.}
\bib{Ta Li}{A note on integral transform~/\!/~Proc. AMS.---1961.---Vol.~12.---No.~6.---P.~556.}
\bib{Love~E.~R.}{Some integral equations involving hypergeometric functions~/\!/~Proc. Edinburgh Math. Soc.---1967.---Vol.~15.---No.~3.---P.~169--198.}
\bib{Love~E.~R.}{Two more hypergeometric integral equations~/\!/~Proc. Cambridge Phil. Soc.---1967.---Vol.~63.---No.~4.---P.~1055--1076.}
\bib{Динь Хоанг Ань}{Интегральные уравнения с функцией Лежандра в ядрах в особых случаях~/\!/~ДАН БССР.---1989.---Т.33.---№7.---С.591--594.}
\bib{Ситник~С.~М.}{Унитарность и ограниченность операторов Бушмана--Эрдейи нулевого порядка гладкости~/\!/~Препринт ИАПУ ДВО РАН.---Владивосток,~1990.---45~С.}
\bib{Ситник~С.~М.}{Факторизация и оценки норм  в весовых лебеговых пространствах операторов Бушмана-Эрдейи~/\!/~ДАН СССР.---1991.---т.320.---№6.---С.1326--1330.}
\bib{Ситник~С.~М., Ляховецкий~Г.~В.}{Формулы композиций для операторов Бушмана-Эрдейи~/\!/~Препринт ИАПУ ДВО РАН.---Владивосток,~1991.---11~С.}
\bib{Ситник~С.~М.}{Об одной паре операторов преобразования~/\!/ В сб.:~Краевые задачи для неклассических уравнений  математической
физики.
(ответственный редактор В.Н. Врагов).---Новосибирск,~1987.---С.~168-173.}
\bib{Ludwig~D.}{The Radon transform on Euclidean space~/\!/~Math. Meth. in the Appl. Sci.---1980~(2).---P.~108--109.}
\bib{Наттерер~Ф.}{Математические аспекты компьютерной томографии
.---М.:~Мир,~1990.---280~с.}
\bib{Deans~S.~R.}{The Radon Transform and Some of Its Applications.---Dover,~1990.---304~p.}
\bib{Катрахов~В.~В.}{Изометрические операторы преобразования и спектральная функция для одного класса одномерных сингулярных псевдодифференциальных операторов~/\!/~ДАН СССР.---1980.---Т.~251.---№~5.---С.~1048--1051.}
\bib{Катрахов~В.~В.}{Об одной краевой задаче для уравнения Пуассона~/\!/~ДАН СССР.---1981.---Т.~259.---№~5.---С.~1041--1045.}
\bib{Маричев~О.~И.}{Метод вычисления интегралов от специальных функций.---Минск:~Наука и техника, 1978.---312~с.}
\bib{Прудников~А.~П., Брычков~Ю.~А., Маричев~О.~И.}{Вычисление интегралов и преобразование Меллина~/\!/~ВИНИТИ, Итоги науки и техники, Математический анализ.---М.,~1989.---Т.~27.---С.~3--146.}
\bib{Kober~H.}{On a Theorem of Schur and On Fractional Integrals of Purely Imaginary Order~/\!/~Transactions of the American Mathematical Society.---1941.---Vol.~50.---No.~1.---P.~160--174.}
\bib{Virchenko~N., Fedotova~I.}{Generalized Associated Legendre Functions and Their Applications.---World Scientific,~2001.---220~p.}
\bib{Virchenko~N.}{On Some Generalized Symmetric Integral
Operators of Buschman-Erdelyi’s Type~/\!/~Nonlinear Mathematical Physics.---1996.---Vol.3.---No.3–-4.---P.421-–425.}
\bib{Вирченко~Н.~А.}{О некоторых приложениях обобщённых ассоциированных функций Лежандра~/\!/~Укр. Матем. Журн.---1987.---Т.~39.---№~2.---С.~149--156.}
\bib{Karp~D., Savenkova~A., Sitnik~S.~M.}{Series expansions for the third incomplete elliptic integral via partial fraction decompositions~/\!/Journal of Computational and Applied Mathematics.~2007.---Vol.~207.---No.~2.---P.~331--337.}
\bib{Karp~D., Sitnik~S.~M.}{Asymptotic approximations for the first incomplete elliptic integral near logarithmic singularity~/\!/~Journal of Computational and Applied Mathematics.---2007.---Vol.~205.---P.~186--206.}
\bib{Karp~D., Sitnik~S.~M.}{Inequalities and monotonicity of ratios for generalized hypergeometric function~/\!/~Arxiv:~math.~CA/0703084.---2007.---14~P.}
\bib{Karp~D., Sitnik~S.~M.}{Two-sided inequalities for generalized hypergeometric function~/\!/~RGMIA Research Report Collection.---2007.--- 10(2).---Article~7.---14~P.}
\bib{Karp~D., Sitnik~S.~M.}{Asymptotic approximations for the first incomplete elliptic integral near logarithmic singularity~/\!/~2006.---Arxiv:~math.~CA/0604026.---20~P.}
\bib{Karp~D., Savenkova~A., Sitnik~S.~M.}{Series expansions and asymptotics for incomplete elliptic integrals
via partial fraction decompositions~/\!/~Proceedings of the fifth annual conference of the
Society for special functions and their applications (SSFA).---2004.---Lucknow (India).---P.~4--30.}
\bib{Ситник~С.~М.}{Неравенства для функций Бесселя~/\!/~ДАН СССР.---1995.---т.~340.---№1.---С.~29--32.}
\bib{Волович~В.~В., Минин~Л.~А., Мордвинов~В.~В., Ситник~С.~М.}{Формула Обрешкова и её применение  в  теории  аппроксимаций
Паде~/\!/~~Препринт ИАПУ ДВО РАН.---Владивосток,~1994.---20~С.}
\bib{Минин~Л.~А., Ситник~С.~М.}{Аппроксимации Паде элементарных и специальных функций~/\!/~~Препринт ИАПУ ДВО РАН.---Владивосток,~1991.---22~С.}
\bib{Ситник~С.~М.}{Неравенства для полных эллиптических интегралов Лежандра~/\!/~~Препринт ИАПУ ДВО РАН.---Владивосток,~1994.---17~С.}
\bib{Ситник~С.~М.}{Неравенства для остаточного члена ряда Тэйлора
экспоненциальной функции~/\!/~~Препринт ИАПУ ДВО РАН.---Владивосток,~1994.---30~С.}

\bib{Ситник~С.~М.}{Краевая задача с интегральными граничными условиями для
одного класса уравнений с сильным вырождением~/\!/~В сб.:~Материалы XXI Всесоюзной научной студенческой конференции
"Студент и научно-технический прогресс". Математика.---Новосибирск,~1983.---С.~55--58.}
\bib{Катрахов~В.~В., Ситник~С.~М.}{Краевая задача для стационарного уравнения Шрёдингера с сингулярным потенциалом~/\!/~ДАН СССР.---1984.---т.278.---№4.---С.797--799.}
\bib{Ситник~С.~М.}{Об унитарных операторах преобразования. Рукопись депонирована в ВИНИТИ 13.11.1986, N 7770--В86.---1986.---Воронежский университет, Воронеж.---9~с.}
\bib{Ситник~С.~М.}{Операторы преобразования  для  дифференциального  выражения Бесселя. Рукопись депонирована в ВИНИТИ 23.01.1987, N 535--В87.---1987.---Воронежский университет, Воронеж.---28~с.}
\bib{Ситник~С.~М.}{$L_2$--теория операторов преобразования и её приложения к  эллиптическим уравнениям~/\!/ В сб.:~10-е Чехословацко--Советское совещание: Применение фукциональных методов и методов теории функций к  задачам математической физики.---1988.---Стара Тура, Чехословакия, С.~46.}
\bib{Ситник~С.~М., Фадеев~С.~А.}{Об одной паре операторов преобразования.
Рукопись депонирована в ВИНИТИ  13.07.1988, N 5629--В88. Воронежский политехнический институт.---1988.---Воронеж.}
\bib{Ситник~С.~М.}{Операторы преобразования  Бушмана--Эрдейи  в  функциональных пространствах~/\!/ В сб.:~"Линейные операторы в функциональных пространствах". Тезисы докладов Северо--Кавказской региональной конференции.---1989.---Грозный.---С.~150.}
\bib{Ситник~С.~М.}{Операторы Бушмана--Эрдейи~/\!/В сб.:~ Дифференциальные и интегральные  уравнения.  Математическая
физика и специальные функции. Международная научная конференция.---1992.---Самара.---С.~233--234.}
\bib{Sitnik~S.~M.}{On unitary transmutations for the Bessel operator~/\!/ В сб.:~Международный семинар Day on diffraction 2004.---Санкт--Петербург.---2004.---С.~71.}
\bib{Opic~B., Kufner~A.}{Hardy--Type Inequalities.---Longman,1990.}
\bib{Kufner~A., Maligranda~L., Persson~L.-E.}{The Hardy Inequality.---Pilsen,~2007.---162~P.}

\bib{Волк~В.~Я.}{О формулах обращения для дифференциального уравнения с особенностью при $x=0$~/\!/~УМН.---1953.---Т.~111.-Вып.~4~(56).---С.~141--151.}
\bib{Сохин~А.~С.}{Об одном классе операторов преобразования~/\!/~Тр. физ.--тех. ин-та низких температурур АН УССР.---1969.---Вып.~1.---С.~117--125.}
\bib{Сташевская~В.~В.}{Метод операторов преобразования~/\!/~ДАН CCCР.---1953.---Т.~113.---№~3.}
\bib{Сташевская~В.~В.}{Об обратной задаче спектрального анализа для дифференциального оператора с особенностью в нуле~/\!/~Уч. зап. Харьковского матем. об--ва.---1957.---№~5.---С.~49--86.}
\bib{Ахиезер~Н.~И.}{К теории спаренных интегральных уравнений~/\!/~~Уч. зап. Харьковского гос. универ.---1957.---Т.~80.---С.~5--21.}
\bib{Chebli~H.}{Op\'{e}rateurs de translation g\'{e}n\'{e}ralises et semigroupes de convolution.---Springer Lect. Notes.---1974.---No.~404.---P.~35--59.}
\bib{Chebli~H.}{Sur un th\`eor\'eme de Paley--Winer associ\'e \`a la d\'ecomposition spectrale d'un op\'erateur de Sturm--Liouville sur $(0,\infty)$~/\!/~J. Funct. Anal.---1974.---No.~17.---P.~447--461.}
\bib{Chebli~H.}{Positivit\'e des op\'erateurs de "translation g\'{e}n\'{e}ralises" associ\'e \`a un op\'erateur de Sturm--Liouville et quelques applications a l'analyse harmonique.---These.---Strasbourg,~1974.}
\bib{Chebli~H.}{Th\`eor\'eme de Paley--Winer associ\'e \`a un op\'erateur diff\'erentiel singulier sur $(0,\infty)$~/\!/~Jour. Math. Pures Appl.---1979.---No.~58.---P.~1--19.}
\bib{Trimeche~K.}{Transformation int\'egrale de Weil et th\`eor\'eme de Paley--Winer associ\'es \`a un op\'erateur diff\'erentiel singulier sur $(0,\infty)$~/\!/~Jour. Math. Pures Appl.---1981.---No.~60.---P.~51--98.}
\bib{Trimeche~K.}{Transformation int\'egrale de Riemann--Liouville g\'{e}n\'{e}ralises et convergence des series de Taylor g\'{e}n\'{e}ralis\'ees au sens de Delsarte~/\!/~Rev. Fac. Sci. Tunis.---1981.---No.~1.---P.~7--14.}
\bib{Bloom~R., Xu~Z.}{Fourier multipliers for $L_p$ on Chebli--Trimeche hypergroups~/\!/~Proc. of the London Mathematical Society.---2000.---Vol.~80.---P.~643--664.}
\bib{Kamel~J., Yacoub~C.}{Huygens’ principle and equipartition of energy for the modified wave equation associated to a generalized radial Laplacian~/\!/~Annales mathematiques Blaise Pascal.---2005.---No.~12.---P. 147--160.}
\bib{Mourou~M., Trimeche~K.}{Op\'erateurs de transmutation et th\`eor\'eme de Paley–Wiener associ\'es a un op\'erateur aux derivees et differences sur $\mathbb{R}$.~/\!/~\\C. R. de l’Acad. des Sci. Series I -- Mathematics.---2001.---Vol.~332.---Issue~5.---P.~397--400.}
\bib{Boumenir~A., Nashed~M.~Z.}{Paley--Wiener type theorems by transmutations~/\!/~Journal of Fourier Analysis and Applications.---2001.---Vol.~7.---No.~4.---P.~395--417.}
\bib{Mahmoud~N.~H.}{Partial Differential Equations with Matricial Coefficients and Generalized Translation Operators~/\!/~Transactions of the American Mathematical Society.---2000.---Vol.~352.---No.~8.---P.~3687--3706.}
\bib{Chebli~H., Fitouhi~A., Hamza~M.~M.}{Expansion in series of Bessel functions for perturbed Bessel operators~/\!/~J. Math. Anal. Appl.---1994.---Vol.~181.---P.~789--802.}
\bib{Fahem~N.~H.}{Theoreme de Paley-Wiener associe a un operateur differentiel singulier a coefficients matriciels~/\!/~C. R. Acad. Sci. Paris.---1985.---Vol.~301.---P.~821--823.}
\bib{Mahmoud~N.~H.}{Differential operators with matrix coefficients and transmutations~/\!/~Contemporary Mathematics.---1995.---Vol.~183.---P.~261--268.}

\bib{Ситник~С.~М.}{Оператор преобразования  и  представление Йоста для уравнения с сингулярным потенциалом~/\!/~~Препринт ИАПУ ДВО РАН.---Владивосток,~1993.---21~С.}
\bib{Катрахов~В.~В., Ситник~С.~М.}{Оценки решений Йоста одномерного уравнения Шредингера с сингулярным потенциалом~/\!/~ДАН СССР.---1995.---Т.~340.---№~1.---С.~18--20.}
\bib{Sitnik~S.~M.}{Explicit solutions to a singular differential
equation with Bessel operator~/\!/~Days on diffraction 2007. International conference. Abstracts.---Saint Petersburg,~2007.---P.~78.}
\bib{Шадан~К., Сабатье~П.}{Обратные задачи в квантовой теории рассеяния.---М.:~Мир,~1980.---408~с.}
\bib{Боровских~А.~В.}{Формула распространяющихся волн для одномерной неоднородной среды~/\!/~Дифференц. уравнения.---2002.---Т.~38.---№~6.---С.~758--767.}
\bib{Боровских~А.~В.}{Метод распространяющихся волн~/\!/~Труды семинара имени И.~Г.~Петровского.---2004.---Вып.~24.---С.~3--43.}
\bib{Ким Хе Ен.}{Продолжение решений систем дифференциальных уравнений в частных производных~/\!/~Вест. МГУ, сер. Математика, Механика.---1991.---№~2.---С.~75--78.}

\bib{Катрахов~В.~В., Ситник~С.~М.}{Метод факторизации в теории операторов преобразования~/\!/~В сб.:~(Мемориальный сборник памяти Бориса Алексеевича Бубнова). Неклассические уравнения и уравнения смешанного типа.
(ответственный редактор В.~Н.~Врагов).---1990,~Новосибирск.---С.~104--122.}
\bib{Катрахов~В.~В., Ситник~С.~М.}{Композиционный метод
построения В--эллиптических, В--гиперболических и В--параболических операторов преобразования~/\!/~ДАН СССР.---1994.---Т.~337.---№~3.---С.~307--311.}
\bib{Титчмарш~Е.}{Введение в теорию интегралов Фурье.---М.--Л.:~ГИТТЛ,~1948.---480~с.}
\bib{Ozaktas~H., Zalevsky~Z., Kutay~M.}{The Fractional Fourier Transform: with Applications in Optics and Signal Processing.---Wiley,~2001.---513~p.}
\bib{Абжандадзе~З.~Л., Осипов~В.~Ф.}{Преобразование Фурье--Френеля и некоторые его приложения.---Санкт--Петербург,~изд.--во С.--П.У.---224~с.}
\bib{Осипов~В.~Ф.}{Почти периодические функции Бора--Френеля.---Санкт--Петербург,~изд.--во С.--П.У.---312~с.}
\bib{Карп~Д.~Б., Ситник~С.~М.}{Дробное преобразование Ханкеля и его приложения в математической физике.---Препринт. Институт автоматики и процессов управления ДВО РАН.---1994.---Владивосток.---24~с.}
\bib{Matveev~V.B.}{Intertwining relations between the Fourier transform
and discrete Fourier transform, the related functional
identities and beyond~/\!/~Inverse Problems.---2001.---Vol.~17.---P.~633-–657.}

\bib{Абловиц~А., Сигур~Х.}{Солитоны и метод обратной задачи.---М.:~Мир, 1979.---479~с.}
\bib{Захаров~В.~Е., Манаков~С.~В., Новиков~С.~П., Питаевский~Л.~П.}{Теория солитонов: метод обратной задачи.---М.:~Наука,~1980.---319~с.}

\bib{Колтон~Д., Кресс~Р.}{Методы интегральных уравнений в теории рассеяния.---М.:~Мир,~1987.---311~с.}
\bib{Лакс~П., Филлипс~Р.}{Теория рассеяния для автоморфных функций.---М.:~Мир,~1979.---330~с.}
\bib{Лакс~П., Филлипс~Р.}{Теория рассеяния.---М.:~Мир,~1971.---312~с.}
\bib{Рамм~А.~Г.}{Многомерные обратные задачи теории рассеяния.---М.:~Мир,~1994.---494~с.}

\bib{Carroll~R.~W.}{Topics in Soliton Theory.---North--Holland Mathematics Studies,~167.---1991.---428~P.}
\bib{Журавлёв~В.~М.}{Нелинейные волны. Точно решаемые задачи.---Ульяновск,~2001.---212~с.}
\bib{Ньюэлл~А.}{Солитоны в математике и физике.---М.:~Мир,~1989.---326~с.}
\bib{Марченко~В.~А.}{Нелинейные уравнения и операторные алгебры.---Киев:~Наукова Думка,~1986.---152~с.}

\bib{Катрахов~В.~В.}{Об одной сингулярной краевой задаче для уравнения Пуассона.---Матем. сб.---1991.---Т.~182.---№~6.---С.~849–-876.}
\bib{Киприянов~И.~А., Катрахов~В.~В.}{Об одной сингулярной эллиптической краевой задаче в областях на сфере~/\!/~Препринт ИПМ ДВО РАН.---1989.---10~с.}
\bib{Киприянов~И.~А., Катрахов~В.~В.}{Сингулярные краевые задачи для некоторых эллиптических уравнений высших порядков~/\!/~Препринт ИПМ ДВО РАН.---1989.---26~с.}
\bib{Катрахов~В.~В.}{Сингулярные краевые задачи для некоторых эллиптических уравнений в областях с угловыми точками~/\!/~ДАН СССР.---1991.---С.~1047--1050.}
\bib{Киприянов~И.~А., Катрахов~В.~В.}{Об одной краевой задаче для эллиптических уравнений второго порядка в областях на сфере~/\!/~ДАН СССР.---1990.---Т.313.---№3.---С.~545--548.}

\bib{Киприянов~И.~А., Катрахов~В.~В.}{Об одном классе многомерных сингулярных псевдодифференциальных операторов~/\!/~Матем. сб.---1977.---Т.~104.---№~1.---С.~49--68.}

\bib{Ландис~Е.~М.}{Задачи \ Е.~М.~Ландиса~/\!/~УМН.---1982.---Т.~37.---№~6.---С.~278--281.}
\bib{Ситник~С.~М.}{О скорости убывания решений некоторых эллиптических
и  ультраэллиптических уравнений~/\!/~Дифференциальные уравнения.---1988.---Т.~24.---№~3.---С.~538--539.}
\bib{Ситник~С.~М.}{О скорости убывания решений стационарного уравнения
Шрёдингера с потенциалом, зависящим от одной переменной~/\!/~В сб.:~
Неклассические уравнения математической физики. (ответственный редактор В.Н. Врагов).---Новосибирск,~1985.---С.~139--147.}
\bib{Ситник~С.~М.}{О скорости убывания решений некоторых эллиптических и  ультраэллиптических уравнений.---Рукопись депонирована в ВИНИТИ 13.11.1986, N  7771--В~86.---1986.---Воронежский университет,~Воронеж.---19~с.}
\bib{Мешков~В.~З.}{Весовые дифференциальные неравенства и их применение для оценок скорости убывания на бесконечности решений эллиптических уравнений второго порядка~/\!/~Тр. МИАН СССР им. В.~А.~Стеклова.---1989.---№~190.---С.~139--158.}
\bib{Мешков~В.~З.}{О возможной скорости убывания на бесконечности решений уравнений в частных производных второго порядка~/\!/~Матем. сб.---1991.---Т.~182.---№~3.---С.~364–-383}

\bib{Успенский~С.~В.}{О теоремах вложения для весовых классов~/\!/~Тр. МИАН СССР им. В.~А.~Стеклова.---1961.---№~60.---С.~282--303.}
\bib{Лизоркин~П.~И.}{Классы функций, построенные на основе усреднений по сферам. Случай пространств Соболева~/\!/~Тр. МИАН СССР им. В.~А.~Стеклова.---1990.---№~192.---С.~122--139.}
\bib{Лейзин~М.~А.}{К теоремам вложения для одного класса сингулярных дифференциальных операторов в полупространстве~/\!/~Дифференц. уравнения.---1976.---Т.~12.---№~6.---С.~1073--1083.}
\bib{Лейзин~М.~А.}{О вложении некоторых весовых классов~/\!/~В сб.:~Методы решений операторных уравнений.---Воронеж:~изд.--во ВГУ,~1978.---С.~96--103.}

\bib{Dunkl~Ch.}{Differential--difference operators associated to reflection groups~/\!/~.---Trans. Amer.
Math. Soc.---1989.---Vol.~311.---P.~167--183.}
\bib{Dunkl~Ch.}{Intertwining operators associated to the group S3.---Trans. Amer. Math. Soc.---1995.---Vol.~347.---P.~3347--3374.}
\bib{Dunkl~Ch.}{An Intertwining Operator for the Group B2.---arXiv:~math.~CA/~0607823.---2006.---27~p.}
\bib{Gallardo~L., Trimeche~Kh.}{Un analogue d'un theoreme de Hardy pour la transformation de Dunkl~/\!/~C. R. Acad. Sci. Paris, Ser. I.---2002.---Vol.~334.---P.~849–-854.}
\bib{Gallardo~L., Trimeche~Kh.}{A  version of Hardy's theorem for the Dunkl transform~/\!/~J. Aust. Math. Soc.---2004.---Vol.~77.---P.~371-385.}
\bib{Chouchene~F.}{Harmonic analysis associated with the Jacobi-Dunkl operator~/\!/~Journal of Computational and Applied Mathematics.---2005.---Vol.~178.---Issue~1--2.---P.~75--89.}
\bib{Mourou~M.}{Taylor series associated with a differential--difference operator on the real line~/\!/~Journal of Computational and Applied Mathematics.---2003.---Vol.~153.---Issues~1--2.---P.~343--354.}
\bib{Dimovski~I., Hristov~V., Sifi~M.}{Commutants of the Dunkl operators in C(R)~/\!/~Fractional Calculus and Applied Analysis.---2006.---Vol.~9.---No.~3.---P.~195--213.}
\bib{R\"osler~M.}{Positivity of Dunkl’s intertwining operator~/\!/~Duke Math. J.---1999.---Vol.~98.---P.~445--463.}
\bib{R\"osler~M.}{Dunkl operators: theory and applications. Orthogonal Polynomials and Special
Functions (Leuven 2002), (E. Koelink and W. Van Assche, eds.), LNIM 1817, Springer,
Berlin--Heidelberg--New York,~2003.---P.~93--135.}
\bib{Trimeche~Kh.}{Inversion of the Dunkl intertwining
operator and its dual using Dunkl wavelets~/\!/~Rocky Mountain
Journal Of Mathematics.---2002.---Vol.~32.---No.~2.---P.889--895.}
\bib{Лянце~В.~Э.}{Несамосопряжённый разностный оператор~/\!/~ДАН СССР.---1967.---Т.~173.---№.~6.---С.~1260--1263.}

\bib{Ханмамедов~Аг.~Х.}{Операторы преобразования для возмущённого разностного уравнения Хилла и их одно приложение~/\!/~СМЖ.---2003.---Т.~44.---№~4.---С.~926--937.}

\bib{Sitnik~S.~M.}{Refinements of the Bunyakovskii -- Schwartz inequalities with applications to special functions estimates~/\!/~ Conference in Mathematical Analysis and Applications in Honour of Lars Inge Hedberg's 60 th Birthday. Link\"{o}ping University, Sweden, 1996.---P.~97.}
\bib{Ситник~С.~М.}{Уточнение интегрального неравенства  Коши--Буняковского~/\!/~Вестник Самарского государственного  технического  университета.  Сер. "Физико-математические науки".---2000.---№~9.---С.~37--45.}
\bib{Ситник~С.~М.}{Некоторые приложения уточнений неравенства
Коши--Буняковского~/\!/~Вестник Самарской государственной экономической академии.---2002.---№1~(8).---С.~302--313.}
\bib{Ситник~С.~М.}{Обобщения неравенства Коши--Буняковского в пространствах
с индефинитной метрикой~/\!/~В сб.:~Материалы шестой Казанской международной летней школы-конференции. Теория функций, её приложения и смежные вопросы.---Труды математического центра имени Н.~И.~Лобачевского.---Том 19.---Казань,~2003.---C.~202--203.}
\bib{Ситник~С.~М.}{Об уточнениях интегрального неравенства Коши-Буняковского~/\!/~В сб.:~Современные проблемы теории функций и их приложения. Тезисы докладов 11--й Саратовской зимней школы, посвящённой памяти выдающихся профессоров МГУ Н.~К.~Бари  и  Д.~Е.~Меньшова.---Саратов,~2002.---C.~ 191--192.}

\bib{Ситник~С.~М.}{Обобщения неравенств Коши-Буняковского и их приложения~/\!/~В сб.:~Abstracts of Bunyakovsky International Conference.---Kyiv,~2004.---C.~119--120.}
\bib{Ситник~С.~М.}{Метод получения  последовательных уточнений неравенства Коши--Буняковского и его применения к оценкам специальных функций~/\!/~В сб.:~Современные методы теории краевых задач.
Материалы Воронежской  весенней математической школы 'Понтрягинские чтения-VII'.---Воронеж,~1996.---С.~164.}
\bib{Ситник~С.~М.}{Уточнения интегрального неравенства Коши-Буняковского и их
приложения к дифференциальным уравнениям~/\!/~В сб.:~ Международная конференция 'Дифференциальные уравнения и смежные вопросы', посвящённая памяти И.~Г.~Петровского. Тезисы докладов.---Москва,~МГУ.---2004.---С.~210.}
\bib{Ситник~С.~М.}{Обобщения неравенств Коши-Буняковского методом средних значений и их приложения~/\!/~Чернозёмный альманах научных исследований. Серия "Фундаментальная математика".---2005---№1~(1).---С.~3--42.}
\bib{Дикарева~Е.~В., Телкова~С.~А., Ситник~С.~М.}{Некоторые приложения уточнений неравенства Коши-Буняковского
в математическом анализе~/\!/~Вестник Воронежского Института МВД России.---Воронеж,~2005.---№5~(24).---С.~65--69.}
\bib{Дикарева~Е.~В., Телкова~С.~А., Ситник~С.~М.}{Средние значения и уточнения неравенств~/\!/~Вестник Воронежского Института МВД России.---Воронеж,~2005.---№5~(24).---С.~57--62.}
\bib{Ситник~В.~С., Ситник~С.~М.}{Исторические аспекты доказательства и обобщений некоторых неравенств~/\!/~В сб.:~Современные методы теории краевых задач.
Материалы Воронежской  весенней математической школы 'Понтрягинские чтения-XVI'.---Воронеж,~2006.---С.~170--171.}
\bib{Ситник~С.~М.}{О некоторых обобщениях неравенства Коши - Буняковского~/\!/~В сб.:~" Международная конференция 'Анализ и особенности', посвящённая 70--летию В.И. Арнольда ".
Тезисы докладов.---Москва, МИАН.---2007.---С.~92--94.}
\bib{Ситник~С.~М.}{Обобщения неравенства Коши --- Буняковского с
использованием средних~/\!/~В сб.:~ Международная конференция 'Дифференциальные уравнения и смежные вопросы', посвящённая памяти И.~Г.~Петровского. Тезисы докладов.---Москва,~МГУ.---2007.---С.~297--298.}
\bib{Ситник~С.~М., Анциферова~Г.~А.}{Некоторые обобщения неравенства Юнга~/\!/~Вестник Воронежского института МВД России.---Воронеж,~1999.---№2~(4).---С.~161--164.}
\bib{Ситник~С.~М.} {Сколько неравенств заключено в неравенстве Юнга?~~/\!/~
Актуальные проблемы обучения математике. Труды Всероссийской научно-практической конференции. Изд Орловского госуд. университета, Орёл, ~2007, ---С.~464--469.}

\bib{Dzrbashian~M.~M.}{ Harmonic Analysis and Boundary Value Problems in the Complex Domain.---Birkh\"auser, 1993.---257~p.}
\bib{Kiryakova~V.}{ Обобщённое дробное исчисление и его приложения в Анализе.---Докторская диссертация, София, Болгария, 2010.---532~c. (Личное сообщение).}
\bib{Kiryakova~V.}{ The multi-index Mittag-Leffler functions as an important class of special
functions of fractional calculus~/\!/~Computers and Mathematics with Applications.---2010.---Vol.~59, Issue~5.--- P.~ 1885--1895.}
\bib{Вирченко Н.~А., Рыбак В.~Я.}{  Основы дробного интегродифференцирования.---Киев, 2007.---362~с. (на украинском языке).}
\bib{Вирченко~Н., Гайдей~В.}{  Классические и обобщённые многопараметрические функции.---Киев, 2008.---227~с. (на украинском языке).}
\bib{Kilbas~A.~A., Saigo~M.}{ $H$ -- transforms. Theory and applications. Chapman and Hall,~CRC.---2004.}
\bib{Gorenflo~R.,  Mainardi~F.}{ Fractional calculus: Integral and differential equations of fractional order, in: A. Carpinteri, F. Mainardi (Eds.), Fractals and
Fractional Calculus in Continuum Mechanics, Springer, Wien, New York, 1997.---P.~ 223--278.}
\bib {E.T. Copson.}{On a Singular Boundary Value Problem
for an Equation of Hyperbolic Type// Arch. Ration.Mech. and Analysis 1 (1957), No. 1, 349 - 356.}
\bib {E.T. Copson.}{ Partial Differential Equations . Cambridge University Press, 1975.}
\bib {E.T. Copson, A. Erdelyi.}{ On a Partial Differential Equation
with Two Singular Lines// Arch. Ration.Mech. and Analysis 2 (1958), No. 1, 76--86.}
\bib{Ситник~С.~М. }{Унитарные операторы преобразования, связанные с операторами Харди.
Третья международная конференция
"Функциональные пространства. Дифференциальные операторы.
Общая топология. Проблемы математического образования",
посвящённая 85--летию
члена--корреспондента РАН профессора
Льва Дмитриевича Кудрявцева.
Тезисы докладов.
Москва, МФТИ, 2008, С. 324--327.}
\bib{Ситник~С.~М. }{О теоремах вложения для пространств С.Л. СОболева и И.А. Киприянова./
В сб.: Дифференциальные уравнения, Функциональные пространства, Теория приближений.
Тезисы Докладов Международной конференции,
посвященной 100-летию со дня рождения Сергея Львовича Соболева.
Россия, Новосибирск,  5-12 октября 2008.
С. 359.}
\bib{Ситник~С.~М. }{Об одном обобщении операторов Харди
методами теории операторов преобразования./
В сб.: Международная конференция по математической физике
и её приложениям.
Самарский государственный университет.
Самара, Россия, 8--13 сентября 2008 г.
Тезисы докладов.
С. 191--192.
}
\bib{Ситник~С.~М. }{О некоторых задачах теории операторов преобразования/
Международная конференция
"Современные проблемы математики, механики и их приложений",
посвящённая 70--летию ректора МГУ академика
Виктора Антоновича Садовничего.
Материалы конференции, С. 49--50. МГУ, 2009.}
\bib{Ситник~С.~М. }{Обзор теории операторов преобразования.
В сб.: Современные методы теории краевых задач.
Материалы Воронежской весенней математической школы
"Понтрягинские чтения--XX",
посвящённой 70--летнему юбилею
академика Виктора Антоновича Садовничего.
Воронеж, ВГУ, 2009, С. 167--168.}
\bib{Ситник~С.~М. }{Fractional powers of the Bessel differential operator.
Аналитические методы анализа и дифференциальных уравнений (АМАДЕ).
Тезисы докладов международной конференции,
14-19 сентября 2009 года, Минск, Беларусь.
С. 149.}
\bib{Ситник~С.~М. }{Решение задачи об унитарном обобщении операторов
преобразования Сонина--Пуассона.---Научные ведомости
Белгородского государственного университета.---2010.---Выпуск~18, №~5 (76).---С.~135--153.}
\bib{Ситник~С.~М. }{Ограниченность операторов преобразования Бушмана--Эрдейи.
Труды 5-ой международной конференции
"Analytical Methods of Analysis and Differential Equations (AMADE)"
(Аналитические методы Анализа и дифференциальных уравнений).
Том 1: Математический Анализ.
Национальная Академия наук Белоруси, институт математики.
Минск, 2010. С. 120--125.}
\bib{Ситник~С.~М. }{Построение унитарных обобщений операторов преобразования
Сонина и Пуассона.
Тезисы докладов Всероссийского научного семинара
"Неклассические уравнения математической физики",
посвящённого 65--летию со дня рождения
профессора Владимира Николаевича Врагова.
Часть II, Якутск 2010.
С. 41--44.}
\bib{Ситник~С.~М. }{Операторы преобразования Бушмана--Эрдейи и их приложения.
Материалы международного Российско--Болгарского симпозиума
"Уравнения смешанного типа и родственные проблемы
анализа и информатики".
Нальчик--Хабез, КБР.
Научно--исследовательский институт
прикладной математики и автоматизации,
Нальчик, 2010, С. 223--224.}
\bib{Ситник~С.~М. }{Теория и приложения операторов Бушмана--Эрдейи.
Материалы 13 международной научной конференции
имени академика М.Кравчука.
Т.2. "Алгебра. Геометрия. Математический и численный анализ".
Киев, 2010, С. 249.}

\bib{Ситник~С.~М. }{Метод факторизации операторов преобразования
в теории дифференциальных уравнений.
Вестник Самарского Государственного Университета (СамГУ) — Естественнонаучная серия.
№ 8/1 (67), 2008, С.  237--248.}
\bib{Ситник~С.~М. }{Операторы преобразования для сингулярных дифференциальных уравнений
с оператором Бесселя.
Сборник научных трудов:
Краевые задачи для неклассических уравнений математической физики.
(ответственный редактор В.Н. Врагов).
1989, Новосибирск, с. 179-185. }
\bib{Kilbas A.A., Skoromnik O.V.}{Integral transforms
with the Legendre function of the first kind in the kernels on
$\mathcal{L}_{\protect\nu ,r}$- spaces // Integral Transforms and
Special Functions. 2009. Vol. 20, issue 9. P. 653--672.}
\bib{Килбас А.А., Скоромник О.В.}{Решение многомерного интегрального уравнения первого рода с с функцией Лежандра по пирамидальной области //  Доклады академии наук РАН. 2009. Т. 429, № 4. С. 442--446.}
\bib{Кузнецов~Н.В.}{Формулы следа и некоторые их приложения в теории чисел. Владивосток, Дальнаука, 2003.}

\bib{Сохин~А.~С.}{Об одном классе операторов преобразования~/\!/ Тр. физ.--тех. ин-та низких температурур АН УССР.---1969.---Вып.~1.---С.~117--125.}
\bib{Сохин~А.~С.}{Обратные задачи рассеяния для уравнений с особенностью~/\!/ Тр. физ.--тех. ин-та низких температурур АН УССР.---1971.---Вып.~2.---С.~182--233.}
\bib{Сохин~А.~С.}{Обратные задачи рассеяния для уравнений с особенностями специального вида~/\!/
Теория функций, функциональный анализ и их приложенияю.---Харьков.---1973.---\No~17.---C.~36--64.}
\bib{Сохин~А.~С.}{О преобразовании операторов  для уравнений с особенностью специального вида~/\!/
Вестник Харьковского университета.---1974.---\No~113.---C.~36--42.}
\bib{Ситник~С.М.}{О представлении в интегральном виде решений одного дифференциального уравнения с особенностями в коэффициентах.~/\!/ Владикавказский математический журнал.---2010.---Т.~12, вып.~4, 7~С.}
\bib{Карп\,Д.Б.}{Пространства с гипергеометрическими воспроизводящими ядрами
и дробные преобразования типа Фурье. ~/\!/ Диссертация канд. физ.-мат. наук, Владивосток, 2000.}
\bib{Осипов В. Ф.} { Почти периодические функции Бора-Френеля.
Изд-во Санкт-Петербургского унив., Санкт--Петербург, 1992.}
\bib{Абжандадзе З.Л., Осипов В. Ф.} {Преобразования Фурье--Френеля и некоторые его приложения. Изд-во Санкт-Петербургского унив., Санкт--Петербург, 2000.}
\bib{Ozaktas H. M., Kutay M. A., Zalevsky Z.}
{ The Fractional Fourier Transform with Applications in Optics and Signal Processing.
John Wiley and Sons, 2000.}
\bib{Качалов\,А.П., Курылёв\,Я.В.}{Метод операторов преобразования в обратной задаче рассеяния, одномерный Штарк--эффект.~/\!/ Записки научных семинаров ЛОМИ, т. 179. Математические вопросы теории распространения волн, 19. Под ред. В.М.\,Бабича. Ленинград, Наука, 1989. С. 73--87.}
\bib{Аршава\,Е.А.}{Об одном классе интегральных уравнений со специальной правой частью. ~/\!/ Аналитические методы Анализа и дифференциальных уравнений (AMADE). Труды 5--й международной конференции 14--19 сентября 2009 года, Минск, Беларусь.  Редакторы: А.А.\,Килбас и С.В.\,Рогозин. \ Том 1, Математический Анализ, Минск, 2010. С. 25--29.}
\bib{Гаршин\,С.В., Прядиев\,В.Л.}{Неулучшаемые условия существования и непрерывности производных второго порядка
у решения характеристической задачи для гиперболического уравнения с двумя независимыми переменными~/\!/
Чернозёмный альманах научных исследований. Серия "Фундаментальная математика".---2005---№1~(1).---С.~83--98.}
\bib{Гринив\,Р.О., Микитюк\,Я.В.}{О подобии возмущённых операторов умножения~/\!/Математические заметки. т. 70, вып. 1, 2001, С. 38--45.}
\bib{Hryniv\,R.O., Mykytyuk Ya.V.}{Transformation Operators for Sturm--Liouville
Operators with Singular Potentials. Dedicated to Professor V. A.Marchenko on the occasion of his 80th birthday.
~/\!/ Mathematical Physics, Analysis and Geometry 7, 2004, P. 119--149.}
\bib{Hryniv\,R.O., Mykytyuk Ya.V.}{Inverse spectral problems for Sturm–Liouville
operators with singular potentials.~/\!/ Inverse Problems 19, 2003, P. 665--684.}
\bib{ Albeverio S., Hryniv\,R.O., Mykytyuk Ya.V.}{Scattering theory for Schr\''{o}dinger
operators with bessel-type potentials. 2009, P. 1--33.}
\Adress{Воронеж, Воронежский институт МВД России.\\
\hspace*{12pt} Радиотехнический факультет, кафедра высшей математики.}
{\verb"mathsms@yandex.ru"}
\end{document}